\def\VERSION{22.07.2022}
\def\users{us}
\def\friends{for-friends} 

\documentclass[12pt,a4paper,final]{article} 
\usepackage{color}
\usepackage{amsmath}
\usepackage{amsthm}
\usepackage{amssymb}
\usepackage{epsfig}
\usepackage{psfrag}
\usepackage{graphicx}
\usepackage{textcomp}
\usepackage{pgf}
\usepackage{todonotes}
\numberwithin{equation}{section}
\numberwithin{figure}{section}
\usepackage{upgreek} 
\newcommand{\ITEM}[2]{\parbox[t]{.045\textwidth}{\textrm{#1}}%
         \hfill\parbox[t]{.95\textwidth}{#2}\pagebreak[3]\vspace{.5em}}
\usepackage{mathrsfs,cite}
\usepackage{ifthen}
\usepackage[normalem]{ulem}
\definecolor{cadmiumgreen}{rgb}{0.0, 0.42, 0.24}
\definecolor{camouflagegreen}{rgb}{0.47, 0.53, 0.42}
\definecolor{darkgreen}{rgb}{0.0, 0.2, 0.13}
\definecolor{junglegreen}{rgb}{0.16, 0.67, 0.53}
\definecolor{lasallegreen}{rgb}{0.03, 0.47, 0.19}
\usepackage{cancel}

\newtheorem{theorem}{Theorem}[section]
\numberwithin{theorem}{section}

\newtheorem{proposition}[theorem]{Proposition}

\theoremstyle{definition}

\newtheorem{remark}[theorem]{Remark}

\ifthenelse{\equal{\users}{world}}
{
\newcommand{\REM}[1]{}

	\newcommand{\DELETE}[1]{}

        \newcommand{\COMMENT}[1]{}
        \newcommand{\TCOMMENT}[1]{}

\newcommand{\EEE}{\color{black}}
}	
{

\usepackage[notcite,notref,color]{showkeys}
\usepackage{showkeys}
\definecolor{brown}{rgb}{0.6,0.2,0.2}
\newcommand{\REM}[1]{\marginpar{\bfseries\tiny{\color{blue}#1}}}

\newcommand{\COMMENT}[1]{{\color{red}\underline{#1}}}
 \newcommand{\DELETE}[1]{{\color{brown}\cancel{#1}\color{black}}}

 \newcommand{\TCOMMENT}[1]{{\color{blue}{ #1}}}

\newcommand{\EEE}{\color{black}}
\newcount\hour \newcount\minute
\hour=\time
\divide \hour by 60
\minute=\time
\loop \ifnum \minute > 59 \advance \minute by -60 \repeat
}

\makeatletter
\renewcommand*\env@cases[1][1.2]{%
  \let\@ifnextchar\new@ifnextchar
  \left\lbrace
  \def\arraystretch{#1}%
  \array{@{\,}c@{\ }l@{}}%
}
\makeatother
\newcommand{\dd}{\,\mathrm{d}}
\newcommand{\ol}{\overline}
\newcommand{\weak}{\rightharpoonup}

\newcommand{\sfR}{\mathsf R}

\newcommand{\calR}{\mathcal R}
\newcommand{\stst}{\mathrm{stst}}

\ifthenelse{\equal{\friends}{for-friends}}
{
}

\newcommand\DT[1]{\mathchoice
                 {{\buildrel{\hspace*{.1em}\text{\LARGE.}}\over{#1}}}
                 {{\buildrel{\hspace*{.1em}\text{\Large.}}\over{#1}}}
                 {{\buildrel{\hspace*{.1em}\text{\large.}}\over{#1}}}
                 {{\buildrel{\hspace*{.1em}\text{\large.}}\over{#1}}}}

\newcommand{\lineunder}[2]{\LU{\begin{array}[t]{c}\underbrace{#1}\vspace*{.5em}\end{array}}{\mbox{\footnotesize\rm #2}}}
\newcommand{\linesunder}[3]{\LSU{\begin{array}[t]{c}\underbrace{#1}\vspace*{.5em}\end{array}}{\mbox{\footnotesize\rm #2}}{\mbox{\footnotesize\rm#3}}}
\newcommand{\LU}[2]{\begin{array}[t]{c}#1\vspace*{-1em}\\_{#2}\end{array}}
\newcommand{\LSU}[3]{\begin{array}[t]{c}#1\vspace*{-1em}\\_{#2}\vspace*{-.3em}\\_{#3}\end{array}}
\renewcommand{\d}{\mathrm{d}}  
\newcommand{\eps}{\varepsilon}
\newcommand{\pl}{\partial}
\newcommand{\R}{\mathbb R}

\newcommand{\wt}{\widetilde}
\newcommand{\wb}{\overline}
\def\Domain{D}

\begin{document}
\begin{sloppypar}
\allowdisplaybreaks

\begin{center}
{\large\bfseries
Qualitative study of a 
geodynamical rate-and-state model\\[0.5em]
for elastoplastic shear flows in crustal faults 
 }\bigskip\bigskip

{\large\sc Alexander Mielke}\bigskip

{\it
Weierstra\ss-Institut f\"ur Angewandte Analysis und 
Stochastik,\\\hspace*{.8em}Mohrenstr.39, D-10117 Berlin, Germany
\\[0.3em]
Institut f\"ur Mathematik,  
Humboldt Universit\"at zu Berlin,\\
\hspace*{.8em}Rudower Chaussee 25, D-12489 Berlin, Germany.\\
\hspace*{.8em}E-mail: alexander.mielke@wias-berlin.de
\vspace{0.3em}}
\bigskip

{\large\sc  \ Tom\'a\v s Roub\'\i\v cek}
\bigskip

{\it 
Mathematical Institute, Charles University,\\
\hspace*{.8em}Sokolovsk\'a 83, CZ-186~75~Praha~8,  Czech Republic
\\[0.3em]
Institute of Thermomechanics, Czech Academy of Sciences,\\
\hspace*{.8em}Dolej\v skova 5, CZ-182~00~Praha~8, Czech Republic\\
\hspace*{.8em}E-mail: tomas.roubicek@mff.cuni.cz       
}\bigskip\bigskip

\bigskip\bigskip

\begin{minipage}[t]{36em}
{\small\baselineskip=11pt
{\it Abstract}:
The Dieterich-Ruina rate-and-state friction model is transferred
to a bulk variant and the state variable (aging) influencing
the dissipation mechanism is here combined also with a damage
influencing standardly the elastic response. As the aging
has a separate dynamics, the overall model does not have a
standard variational structure. A one-dimensional model is investigated
as far as the steady-state existence, localization of the cataclastic
core, and its time response, too. Computational experiments with a
damage-free variant show stick-slip behavior (i.e.\ seismic cycles
of tectonic faults) as well as stable slip under very large velocities.
\bigskip


\bigskip

}

\end{minipage}
\end{center}

\def\bbC{\mathbb{C}}

\section{Introduction}\label{sec-intro}

In the last decades the mathematical interest in geophysical problems was
steadily growing. While there is already a large body of work in atmospheric
and oceanographic fluid flows, the mathematics for geophysical models for solid
earth is much less developed. The latter concerns in particular the deformation
and motion of lithospheric plates in the upper crust, in particular
earthquakes.  The difficulties in these models is the complex behavior of rock
that behaves elastically like a solid in the case of seismic waves on short
time scales but behaves like a viscoplastic fluid when considered over
centuries. However, very slow motion of long periods are crucial for building
up internal stresses that are then released in short rupture events triggering
earthquakes. Only recently, a new class of periodic motions in the Earth crust
was detected by evaluating GPS measurements, namely the so-called
``\emph{episodic tremor and slip}'' (cf.\ \cite{KTOT12ETSS,Bart20LTVE}): Here
all motions are so slow that no seismic waves are emitted, but there exist two
distinct regimes, one involving inelastic motions and one involving slow smooth
slip. These events are observed in so-called subduction zones and have periods
in the range of a few years while the overall shear velocity rate is in the range of
millimeter per year.  

In addition to these temporal time scales there are also several spatial scales
involved. For instance, between tectonic plates there form weak
regions called faults that are relatively narrow but may accumulate 
relatively large deformations, in particular in rapid shearing events.  We
refer to \cite{PKRO16ERNS, HeGevD18IRSD, RHCRKO19SGSM, PHGA19SAFG, NGBPW20D3DR}
for some recent efforts in geodynamical modeling towards a better understanding
of these phenomena. On the mathematical side the work started less than a
decade ago and is still comparably small, see \cite{RoSoVo13MRLF, PiSaKo13VFRS,
  Pipping2019, HeKoPo20FHMI, EiHoMi22LHSV, EiHoLa21?WSUE}. Moreover, there is a
dichotomy with respect to bulk interface models, where most of the nonlinear
effects are localized in the interface (e.g. by a so-called rate-and-state
dependent friction law), and pure bulk models where typically only existence
results for solutions are obtained but no qualitative behavior of the solutions
can be deduced.\medskip

With this work we want to initiate a mathematical study where pure bulk models
are considered but still interesting qualitative features can be deduced. In
this first study we will confine ourselves to a simplified ``stratified''
setting where only shear deformations are considered that depend on a
one-dimensional variable $x \in (-H,H)$ representing the transverse direction
to a straight fault or damage zone between two compact rocks representing two
plates that move with respect to each other, see Figure \ref{fig:geometry}.
The continuum model is given in terms of
\\
\textbullet\ the shear velocity $v= v(t,x) \in \R$,
\\
\textbullet\ the elastic strain $\eps = \eps(t,x)$,
\\
\textbullet\ the plastic strain $p= p(t,x)$,
\\
\textbullet\ the internal damage variable $\alpha=\alpha(t,x)$, and
\\
\textbullet\ the internal aging variable $\theta=\theta(t,x)$.

The model to be studied in its simplest form is the following system of five 
partial differential equations posed for $(t,x)\in
(0,\infty)\times (-H,H)$ (see \eqref{evol} for the more general case treated below): 
\begin{subequations}
\label{eq:I.evol}
\begin{align}
\label{eq:I.evol.a}
&\varrho \DT v = \big( \bbC(\alpha) \eps\big)_x , &&\DT\eps + \DT p = v_x,
\\
\label{eq:I.evol.b}
&\pl_{\DT p} R(\DT p,\theta)  \ni \bbC(\alpha) \eps + \eta\DT p_{xx},
 \hspace{-1em}
&& \DT\alpha =- \frac12\bbC'(\alpha) \eps^2 + \beta(1{-}\alpha) + \gamma
   \alpha_{xx} ,\quad \mbox{}
\\
\label{eq:I.evol.c}
\mbox{}\qquad&\DT\theta = 1-\theta/\theta_\infty - \lambda |\DT p|\theta +
\kappa \theta_{xx} ,  \hspace{-2em}
\\
\intertext{with the dot-notation $(\cdot)\!\DT{^{}}$ and the notation $(\cdot)_x$
for the partial derivatives in time and in space, respectively. 
We complete it with boundary conditions} 
\label{eq:I.evol.d}
& 
v(t,\pm H) = \pm v_\infty(t), \ \  p(t,\pm H)=0,   && 
\alpha(t,\pm H)=1, \ \
\theta(t,\pm H) = \theta_\infty. 
\end{align} 
\end{subequations}

Here $\beta,\ \gamma,\ \eta,\ \kappa$, and $\lambda$ are positive constants,
whereas $\alpha \mapsto \bbC(\alpha)>0$ and
$(\pi,\theta) \mapsto R(\pi,\theta)>0$ are general smooth constitutive
functions. In particular, the state of damage $\alpha$ may decrease the elastic
stiffness $\bbC(\alpha)$, and even more importantly the yield stress
$\mu(\pi,\theta)$ may depend on the plastic rate $\pi=\DT p$ as well as on the
aging variable $\theta$. Thus, we are able to mimic the commonly used
Dieterich-Ruina {\it rate-and-state friction} law
\cite{Diet07ARSD,Ruin83SISVFL} where now the aging variable can be interpreted
as the ``state'' while the dependence on $\pi=\DT p$ gives the rate dependence.

Here $R(\;\!\cdot\!\;,\theta):\R\to \R$ is the plastic dissipation potential
depending on the aging variable $\theta$,
i.e.\ it is convex and satisfies $R(\pi,\theta) \geq 0 = R(0,\theta)$. The
plastic yield stress (or dry friction coefficient) is encoded by assuming
$R(\pi,\theta) = \mu(0,\theta)|\pi| + \mathscr{O}(\pi^2)$. Hence, we obtain a set-valued
convex subdifferential, which we assume to have the form $\pl_\pi R(\pi,\theta) =
\mu(\pi,\theta) \mathop{\mathrm{Sign}}(\pi) +\mathscr{O}(\pi)$, where 
``\,Sign'' is the set-valued sign function,  see \eqref{eq:def.Sign}. Thus, 
the first equation in \eqref{eq:I.evol.b}, involving the nonsmooth
convex function $R(\cdot,\theta)$, is an inclusion 
and gives rise to a {\it free boundary}, namely between regions with 
the purely elastic regime with $\pi = \DT
p\equiv 0$ where $\text{Sign}(\DT p) = [-1,1]$ and the plastic regime where
$\pi = \DT p \neq 0$ and $\text{Sign}(\DT p)= \{-1\} $ or $\{+1\}$. \medskip

Our paper is organized as follows: In Section \ref{se:Setup} we provide the
background from geodynamics introducing the rate-and-state friction models with
a given interface and our distributed-parameter model which is slightly more
general than \eqref{eq:I.evol}. In particular, Section \ref{sec-steady} 
discusses the steady-state equation where $\DT v=\DT \alpha=\DT \theta=0$
while the plastic flow rate $\pi=\DT p$ is independent of time.  
The full evolutionary model is then introduced in Section \ref{su:EvolMod}. 

The analysis of steady states is the content of Section \ref{se:AnaSteady}.
In Theorem \ref{th:ExistSteady} we provide an existence theorem for steady states
under quite natural assumptions and arbitrary shear velocities $v(\pm H)=\pm
v_\infty$. The proof relies on a Schauder fix-point 
argument and we cannot infer uniqueness, which is probably false in this
general setting. In Proposition \ref{prop2} we show that for steady states 
the limit $\eta\to 0^+$ in \eqref{eq:I.evol.b}  can be performed in such a way
that accumulation points are still steady states. 

In Section \ref{sec-evol} we discuss the full dynamic model, show its
thermodynamic consistency, and derive the natural a priori estimates. For our
main existence result we restrict to the case without damage, i.e.\ $\bbC$ is
independent of $\alpha$ and $\alpha\equiv 1$ solves \eqref{eq:I.evol.b}. The
result of Theorem \ref{th:EvolExist} is obtained by time
discretization and a staggered incremental scheme mimicking the solution of the
static problem in Theorem \ref{th:ExistSteady}. \EEE The analytical aspects are
nontrivial because of the non-variational character of the problem, the non-polynomial
friction law \eqref{DR1+} leading to usage of Orlicz spaces, and
the lack of compactness for the elastoplastic wave equation.

The final Section~\ref{se:NumSimul} is devoted to a numerical exploration of
some simplified models that show the typical behavior expected also for the
full model. The simplified model is obtained from \eqref{eq:I.evol} by
neglecting $\alpha$ as in Section \ref{sec-evol} and by further ignoring
inertia (i.e.\ setting $\varrho=0$ and choosing $\eta=0$), see Section
\ref{su:SimplifMod}:
\begin{equation}
  \label{eq:I.SimplMod}
  \frac{2H}{\bbC} \DT \sigma + \int_{-H}^H \varPi(\sigma,\theta) \dd x = 2
  v_\infty(t), \quad 
 \DT \theta = 1{-}\frac\theta{\theta_\infty} - \lambda \varPi(\sigma,\theta) +
 \kappa \theta_{xx} , 
\end{equation}
with $\theta(t,\pm H)=\theta_\infty$, where $\pi=\varPi(\sigma,\theta)=
\pl_\xi\calR^*(\sigma,\theta)$ is the unique 
solution of $\sigma \in \pl_\pi\calR(\pi,\theta)$. 

In Section \ref{su:NumSimSteady} we discuss the steady states
$(\theta_\stst,\pi_\stst)$ where $\pi_\stst=
\varPi(\sigma_\stst,\theta_\stst)$. We do a parameter study for varying
$\kappa$ and $v_\infty$ and obtain a monotone behavior with respect to
$v_\infty$, namely $\theta_\stst$ is decreasing and $\pi_\stst$ is increasing. 
We always observe spatial localization in the sense that
$\pi_\stst$ is supported on $[-h_*(v_\infty,\kappa), h_*(v_\infty,\kappa)]$ with
a free boundary positioned at the points $\pm h_*(v_\infty,\kappa)$ with 
$h_*(v_\infty,\kappa) \lneqq H$ and $h_*(v_\infty,\kappa)\approx 0.55
\sqrt\kappa$ for $\kappa, v_\infty \to 0^+$.

The pure existence of steady states does not say anything about stability in
the dynamic model \eqref{eq:I.SimplMod}. In Section \ref{su:NumSimODE} we
provide a two-dimensional ODE model where there is a unique steady state that
is unstable for small positive $v_\infty$ and convergence of general solutions
to periodic motions. Similarly, Section \ref{su:NumSimPDE} shows simulations
for system \eqref{eq:I.SimplMod} which shows convergence towards
$(\theta_\stst,\pi_\stst)$ if $v_\infty$ is large but predicts convergence
towards time-periodic solutions that also have a clearly defined plastic zone
smaller than $(-H,H)$, see Figures \ref{fig:SM.converge} and \ref{fig:SM.osc}.

A surprising effect is that the width $2h$ of the core of the fault (the active
cataclastic zone) does not tend to be 0 if the plasticity gradient is ignored
by setting $\eta=0$, and even not if the aging gradient is ignored by setting
$\kappa=0$. In Proposition \ref{prop3} we show that under natural assumptions
on the rate-and-state friction law one obtains a linear dependence 
$h=h_*(v_\infty,0)= |v_\infty|/\pi_*$  for shear velocities with
$|v_\infty|< H\pi_*$, where $\pi_*$ is uniquely determined by the friction law
and the aging law.  

Another noteworthy effect is that the length scale of the aging qualitatively
influences the character of response, varying in between the stick-slip and the
sliding regimes. In particular, for very large shear velocities $v_\infty$
(which are not relevant in usual geophysical faults in the lithosphere) the
fault goes into a continuous sliding mode and no earthquakes occur. Actually,
this is a recognized attribute of this friction model which in
\cite{Baum96DFDL} has been compared to the observation of our ``everyday life
when one often manages to get rid of door-squeaking by a fast opening''. In
contrast under very slow shear velocities, the friction threshold is not
reached for large time spans after a relaxation. Only when enough shear stress
has build up, the threshold can be overcome. But then not only stresses are
released but also the aging variable is reduced which leads to a much larger
stress release than needed. Hence, another long waiting time is needed until
next ``earthquake'' will start.

\section{Setup of the geodynamical model}
\label{se:Setup}

\subsection{Geodynamical background} 
\label{su:Geodynamics}

{\it Earth's crust} (together with lithosphere) is a rather solid rock bulk
surrounding the lower, more viscous parts of the planet. It is subjected by
damage typically along thin, usually flat weak surfaces, called {\it faults},
which exist within millions of years. The faults may exhibit slow sliding
(so-called aseismic slip) or fast {\it rupture} (causing {\it tectonic earthquakes}
and emitting seismic waves) followed by long period or reconstruction (healing)
in between particular earthquakes. The former phenomenon needs some extra
creep-type rheology modeled using a plastic strain variable
or some smoothing of the activated character of the frictional
resistance at very small rates (cf.\ Remark~\ref{rem-aseismic}) and will not be
scrutinized in this article, while the latter phenomenon needs some
friction-type rheology. Thus faults can be modeled as frictional contact
surfaces or as flat narrow stripes.

As for the frictional contact, the original Dieterich-Ruina
rate-and-state friction model \cite{Diet07ARSD,Ruin83SISVFL}
prescribes the tangential stress $\sigma_{\rm t}$ on the frictional interface as
\begin{align}\label{DR1}
\sigma_{\rm t}=\sigma_{\rm n}\Big(\!\!\!\!
\lineunder{\mu_0+a\,{\rm ln}\frac{v}{v_{\rm ref}}
  +b\,{\rm ln}\frac{v_{\rm ref}\theta}{d_{\rm c}}}
{= $\mu(v,\theta)$ = frictional resistance}\!\!\!\!\Big)
\end{align}
where the normal stress $\sigma_{\rm n}$ is considered to be given (=\,a
so-called Tresca friction model) and $v$ is (the norm of) the tangential
velocity jump along interface.  The (given) parameters $a$ and $b$ are the
direct-effect and the evolution friction parameters, respectively, $d_{\rm c}$
is the characteristic slip memory length, and $v_{\rm ref}$ reference
velocity. If $a{-}b>0$, we speak about velocity strengthening while, if
$a{-}b<0$, we speak about {\it velocity weakening} -- the latter case may lead
to instabilities and is used for earthquake modeling.  The friction
coefficient $\mu=\mu(v,\theta)$ depends in this model on the velocity magnitude
$v$ and an internal variable $\theta$ being interpreted as an {\it aging}
variable, sometimes also as damage. The evolution of $\theta$ is governed by a
specific flow rule typically of the form of an ordinary differential equation
at each spot of the fault, say:
\begin{align}
  \label{DR3}
  \DT\theta=f_0(\theta)-f_1(\theta)|v|\,
\end{align}
with some continuous nonnegative functions $f_0$ and $f_1$
More specifically, $f_0(\theta)=1$ and $f_1(\theta)=\theta/d_{\rm c}$
with $d_{\rm c}>0$ is most common, considered e.g.\ in 
\cite{Bizz11TVCP, BeTuGo08CRPB, BieTel01MCPF, DauCar08CMFG, DauCar10FFE,
  KaLaAm08SEMS, RDDDC09FDMR, Scho98EFL}; 
then for the static case  $v=0$, the aging variable $\theta$ grows linearly
in time and has indeed the meaning of an ``age'' as a time elapsed from 
the time when the fault ruptured in the past. The steady state $\DT\theta=0$
leads to $\theta=d_{\rm c}/|v|$ so that $\mu=\mu_0+(a{-}b)\,{\rm ln}|v/v_{\rm ref}|$.
Alternatively, one can consider the flow rule \eqref{DR3}
with some other $f_0$:
\begin{align}
   \label{DR3+}
f_0(\theta)=\max\Big(1-\frac\theta{\theta_\infty\!}\,,0\Big)\ \ \text{ and }\ \
f_1(\theta)=\frac\theta{d_{\rm c}}\,,
\end{align}
cf.\ \cite{PeRiZh95SHSP},
and then $\theta$ stays bounded and asymptotically approaches $\theta_\infty$
in the steady state if $v\to0$, namely
$\theta=d_{\rm c}\theta_\infty/(d_{\rm c}{+}\theta_\infty|v|)$. This suggests to
interpret $\theta$ rather as a certain hardening or ``gradual locking''
of the fault in the ``calm'' steady state $v=0$.

An obvious undesired attribute of \eqref{DR1} is, as already noted in 
\cite[p.108]{Diet07ARSD}, that, ``as $v$ or $\theta$ approach zero, 
eqn.\ \eqref{DR1} yields unacceptably small (or negative) values of sliding 
resistance'' $\mu$. Therefore, \eqref{DR1} obviously violates the
Clausius-Duhem entropy inequality, although being used in dozens of geophysical
articles relying that in specific applications
the solutions might not slide into these physically wrong regimes. 
Nevertheless, 
a regularization leading to $\mu>0$ and thus to a physically 
correct non-negative dissipation is used, too, typically 
as \cite{Diet87NTES}, cf.\ e.g.\ also \cite{PeRiZh95SHSP}:
\begin{align}\label{DR1+}
&\mu=\mu(v,\theta)=\mu_0+a\,{\rm ln}\Big(\frac{|v|}{v_{\rm ref}\!}\,{+}1\Big)
+b\,{\rm ln}\Big(\frac{v_{\rm ref}}{d_{\rm c}}\theta{+}1\Big)\,.
\end{align}
In what follows, we will therefore have in mind rather \eqref{DR1+} than 
\eqref{DR1}. For an analysis and numerics  of the rate-and-state friction
in the multidimensional visco-elastic context we refer to \cite{PiSaKo13VFRS,
  PKRO16ERNS, PatSof17ARSF, Pipping2019}.

Since the velocity occurs in the aging flow rule \eqref{DR3}, this nonisothermal
friction model however does not seem consistent with standard
thermodynamics as pointed out in \cite{Roub14NRSD} in the sense
that the evolution \eqref{DR3} does not come from any free energy.
On top of it, it has been known from the beginning of this rate-and-state
model that it does not fit well some experiments \cite{Ruin80FLIQ} and 
(rather speculative) modifications e.g.\ by using several aging variables 
(which naturally opens a space for fitting more experiments) have been 
devised, cf.\ \cite{Ruin83SISVFL}.

A rather formal attempt to overcome the mentioned thermodynamical inconsistency
has been done in \cite{PiSaKo13VFRS} by introducing two energy
potentials. Thermodynamically consistent models have been devised either by
using isothermal damage with healing \cite{RoSoVo13MRLF} or by nonisothermal
damage when temperature variation was interpreted approximately as a sliding
velocity magnitude $v$. The latter option uses the idea that the slip of the
lithospheric fault generates heat which increases temperature on the fault. In
geophysical literature, the heat produced during frictional sliding is believed
``to produce significant changes in temperature, thus the change of strength of
faults during seismic slip will be a function of ... also temperature'', cf.\
\cite[p.7260]{Ches94ETFC}. The usage of an (effective) interfacial temperature
discussed in \cite{DauCar08CMFG, DauCar10FFE} following ideas from
\cite{Lang08STZR}.  In \cite{BaBeCa99PASR, Ches94ETFC, Ches95RMWC, Scho02MEF}
the classical rate-and-state friction law is also made temperature dependent.
Experimentally, even melting of rocks due to frictional heating is sometimes
observed.

A simplified friction model $\mu(v)=\mu_0+a{\rm ln}(b|v|{+}1)$ or
$\mu(v)=\mu_0+(a{-}b){\rm ln}|v/v_{\rm ref}|$
is sometimes also considered under the name rate-dependent friction
\cite{Diet79MRFE,LyBZAg05VDRR,Ruin83SISVFL,TonLav18SSTE}
and was analyzed in \cite{Miel18TECI} as far as its stability.
In contrast, the above mentioned variant of temperature dependent 
friction can be called purely state dependent.

The friction model is sometimes ``translated'' into a bulk model
involving a plastic-like strain and the sliding-friction
coefficient $\mu$ then occurs as a threshold (a so-called yield stress)
in the plastic flow rule, cf.\ \cite[Sect.\,6]{Roub14NRSD},
or \cite{DauCar09SSIS,DauCar10FFE,HeGevD18IRSD,LyBZAg05VDRR,TonLav18SSTE}, known
also under the name a {\it shear-transformation-zone} (STZ) concept
referring to a (usually narrow) region in an amorphous solid that
undergoes plastification when the material is under a big mechanical load.
Instead of velocity dependence \eqref{DR1+}, one should play with dependence
on the strain rate, cf.\ \eqref{DR1++} below. These options can be ``translated''
into the bulk model by making the yield stress $\mu$ dependent, beside
the strain rate, also on an aging variable $\theta$,
or on an temperature, or on a damage, or on various combination of those.
Altogether, one thus get a wide menagerie of friction-type models.

Here we consider, as rather standard in geophysical modeling as \eqref{DR1+},
an isothermal variant and make $\mu$ dependent on strain rate and on aging. We
consider also damage (or phase-field) as usual in fracture mechanics  
to illustrate its a different position in the model.
The main phenomena are that aging evolution does not directly contribute
to energetics when influencing only dissipative ``friction'' $\mu$.
This is similar to a cam-clay model \cite{DaDeSo11QECC,DaDeSo12QECC}
where  the dissipative response is controlled through an internal variable
whose rate, however, does not explicitly contribute to energetics.
On the other hand, damage (or phase-field) influences the elastic response
through the elastic response in the stored energy and is also driven by the
resulting driving force from it. Also, 
we adopt the (realistic) assumption that the elastic strain
(as well as its rate) is small, which makes possible to let
$\mu$ dependent on the plastic strain rate rather than
elastic strain rate and to put it into the standard framework of
rate-dependent plasticity. The plasticity is consider without
any hardening which otherwise might dominate with big slips
on long time scales and would unacceptably corrupt the autonomy of
the model. In principle, damage may also influence friction
$\mu$ like in \cite{RoSoVo13MRLF,RouVal16RIPP} but we will not consider it.

\subsection{The one-dimensional steady-state model}\label{sec-steady}

\begin{figure}[ht]
\centering
\psfrag{cataclastic zone}{\hspace*{-.3em}\small\bf cataclastic zone}
\psfrag{compact rock}{\footnotesize\hspace*{0em}
\begin{minipage}[c]{5em}\baselineskip=8pt compact\\\hspace*{.7em}rock\end{minipage}}
\psfrag{damage zone}{\footnotesize damage zone $D$ (width = $2H$)}
\psfrag{2h}{\footnotesize $2h$}
\psfrag{velocity}{\footnotesize\hspace*{0em}
\begin{minipage}[c]{5em}\baselineskip=8pt prescribed\\velocity $v_\infty$\end{minipage}}
\includegraphics[width=28em]{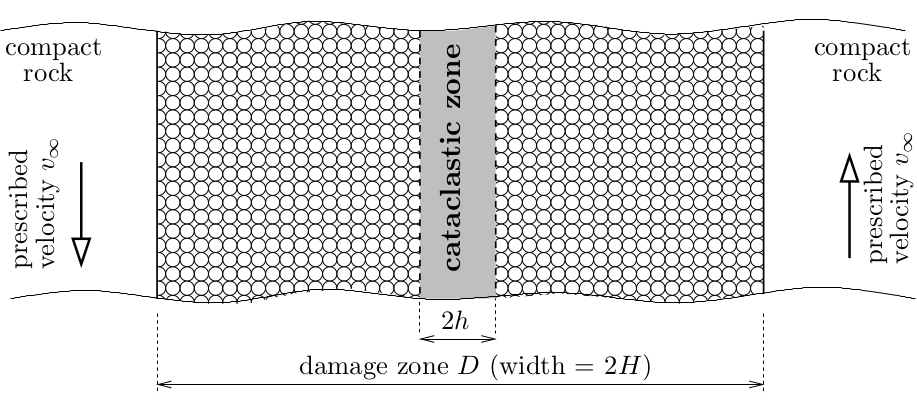}
\vspace*{-.1em}
\caption{\small\sl Schematic geometry: a cross-section through a fault.
}\label{fig:geometry}
\end{figure}

It is generally understood that fracture mechanics and in
particular fault mechanics is very complex and difficult to analyze.
Therefore, we focus to a very simplified situation: a flat
fault which is perfectly homogeneous in its tangential direction.
Thus all variables depend only on the position in the normal
direction and the problem reduces to be one dimensional,
cf.\ Figure~\ref{fig:geometry}.

We ask a question about existence of a steady state in the situations
where the sides of the fault move with a constant speed in opposite
directions. The model is thus expressed in rates rather than
displacements and plastic strains.
Such steady states are also called {\it aseismic slips} (sliding),
in contrast to seismic slips which are dynamical phenomena related
with a stick-slip motion and earthquakes. For the relation of the
aseismic slip (fault growth) and orientation of faults see
\cite{PHGA19SAFG}. The aseismic slip can be also understood as creep,
within which the Maxwellian viscoelastic rheology is manifested.

The variables of our steady-state model will thus be:
\\
\textbullet\ $\ \ v$ velocity (in m/s),
 \\
\textbullet\ $\ \ \pi$  plastic strain rate (in 1/s),
\\
\textbullet\ $\ \ \eps$ elastic strain (dimensionless),
 \\
\textbullet\ $\ \ \alpha$ damage (dimensionless, ranging over $[0,1]$), and
 \\
\textbullet\ $\ \ \theta$ aging (in seconds), and later also
\\
\textbullet\ $\ \ \sigma$ a stress (or, in one-dimensional case, rather a force in J/m=N).

\noindent
These first five variables are to satisfy the following system of five equations (inclusions): 
\begin{subequations}
 \label{eq}\begin{align}\label{eq1}
&(\bbC(\alpha)\eps)_x=0&&\text{(momentum equilibrium)}
\\ \label{eq2}
&\pi=v_x&&\text{(plastic shear rate)}
\\&\label{eq3}
\mu(\pi,\theta)\text{Sign}(\pi)\ni \bbC(\alpha)\eps+\eta\pi_{xx}\,,&&\text{(plastic flow rule)}
\\[-.2em]&\label{eq4}
\frac12\bbC'(\alpha)\eps^2+G_{\rm c}\frac{\alpha{-}1}{\ell^2}
=G_{\rm c}\ell^2\alpha_{xx}\,,&&\text{(damage flow rule)}
\\[-.2em]&\label{eq5}
|\pi|f_1(\theta)-f_0(\theta)=
\kappa\theta_{xx}\,,&&\text{(aging flow rule)}
\end{align}\end{subequations}
where $(\cdot)_x$ denotes the derivative (later also partial derivative) in
$x$. Actually, \eqref{eq3} contains a set-valued term
$\pl_\pi R(\pi,\theta) = \mu(\pi,\theta)\text{Sign}(\pi)$ and is thus an
inclusion rather than an equation. There, we have denoted by ``\,Sign'' in
set-valued sign function, i.e.
\begin{align}
\label{eq:def.Sign}
  \text{Sign}(\pi)=\begin{cases}1&\text{for }\pi>0,\\[-.2em]
  [-1,1]&\text{for }\pi=0.\\[-.2em]
  -1&\text{for }\pi<0.\end{cases}
\end{align}
This system arises as a steady state from an evolution model \eqref{evol} below. 
In particular, the equation \eqref{eq2} arises from the additive (Green-Naghdi's)
decomposition of the total strain into the elastic strain and the plastic
strain, cf.\ \eqref{evol2} below. Written in terms of rates and taking into account
that the rate of the elastic strain is zero in the steady state, we arrive at
\eqref{eq2}. In fact, the velocity $v$ here enters the rest of the system only
through the boundary condition \eqref{BC} below, in contrast to the full evolutionary
model later in Section~\ref{sec-evol} where velocity acts through the inertial
force. 

The data (or constitutive relations) in the model \eqref{eq} are:

$\ \ \ \mu=(\pi,\alpha)$ a yield stress (in the one-dimensional model in N=J/m)),

$\ \ \ \bbC=\bbC(\alpha)$ elastic modulus (smooth, nondecreasing, 
in N=J/m),

$\ \ \ f_0$ aging rate (dimensionless),

$\ \ \ f_1$ ``contra-aging'' coefficient (in seconds),

$\ \ \ G_{\rm c}$ fracture toughness (in a one-dimensional model in N=J/m),

$\ \ \ \eta>0$ a length scale coefficient for $\pi$
        (i.e.\ for the cataclastic zone, in W/m),

$\ \ \ \ell>0$ a length scale coefficient for the damage (in meters),

$\ \ \ \kappa>0$ a length scale coefficient for the aging (in m$^2$/s),

\noindent
while $f_0$ and $f_1$ are essentially borrowed from \eqref{DR3+}.
Actually, $v$ in \eqref{DR1+} has the meaning rather of a difference of
velocities across the contact interface than a velocity itself which would
not be Galilean invariant. In a variant of the bulk model, $\mu$ should
depend rather on a shear rate and, instead of the coefficient
$1/v_{\rm ref}^{}$, one should consider a $h/v_{\rm ref}^{}$ with $h$ a certain
characteristic width of the active slip area, likely to be identified with
the width of the cataclastic core zone, cf.~Figure~\ref{fig:geometry}. 
Thus, we consider
\begin{align}\label{DR1++}
&\mu=\mu(\pi,\theta)=\mu_0+a\,{\rm ln}\Big(\frac{h}{v_{\rm ref}}|\pi|{+}1\Big)
+b\,{\rm ln}\Big(\frac{v_{\rm ref}}{d_{\rm c}}\theta{+}1\Big)\,.
\end{align}

In comparison with \eqref{DR3}, the steady-state equation \eqref{eq5} 
contains the length-scale term $\kappa\theta_{xx}$. Also damage equation \eqref{eq5} 
contains a length-scale term $\ell^2\alpha_{xx}$ competing with  
the driving force $\frac12\bbC'(\alpha)\eps^2$ coming from the $\alpha$-dependence
in \eqref{eq1}. Note that the gradient term in \eqref{eq3} applies to plastic rate
and no gradient term involves directly the plastic strain, similarly as in
\cite{DaRoSt21NHFV,Roub22QHLS}. This eliminates spurious hardening-like
effects by large slips accumulated on faults in large time scales, which would
otherwise start dominating and corrupt the autonomous character of the model.

We have to complete the system \eqref{eq} by suitable boundary condition.
Specifically, we choose the boundary conditions
\begin{align}\label{BC}
v(\pm H)=\pm v_\infty,\ \ \ \ \
    \pi(\pm H)=0,\ \ \ \ \
    \alpha(\pm H)=1,\ \ \ \ \
    \theta(\pm H)=\theta_\infty
\end{align}
with $\theta_\infty$ from \eqref{DR3+}.
Let us mention that we use the mathematical convention that $\alpha=1$ means
undamaged material while $\alpha=0$ means maximally damaged material. 

From \eqref{eq1}, we can see that $\bbC(\alpha)\eps$ is constant on the damage
domain $\Domain=[-H,H]$, say $=\sigma$.
From this, we can express
\begin{align}\label{e=..alpha}
\eps(x)=\frac\sigma{\bbC(\alpha(x))}\qquad\text{ for all }\ x\in\Domain\,.
\end{align}
If $\bbC(\cdot)$ is increasing, one can conversely express $\alpha$ as a
function of $\eps$, but we will eliminate $\eps$ rather than $\alpha$.
Also the equation \eqref{eq2} can be eliminated because the velocity $v$ occurs
only in the first boundary condition in \eqref{BC}. This condition then
turns into an integral side constraint
$\int_\Domain\pi \dd x=\int_\Domain v_x \dd x=v(H)-v(-H)=2v_\infty$.
We can thus reduce \eqref{eq} to the system of three elliptic
ordinary-differential equations
\begin{subequations}\label{eq-three}\begin{align}
    &\label{eq3-three}
\mu(\pi,\theta){\rm Sign}(\pi)\ni\sigma+\eta\pi_{xx}\,,
\\&\label{eq4-three}
\frac{\bbC'(\alpha)}{2\bbC^2(\alpha)}\sigma^2+G_{\rm c}\frac{\alpha{-}1}{\ell^2}
=G_{\rm c}\ell^2\alpha_{xx}\,,
\\&\label{eq5-three}
|\pi|f_1(\theta)-f_0(\theta)=
\kappa\theta_{xx}\,
\end{align}\end{subequations}
with the integral and the boundary conditions
\begin{subequations}\label{BC-three}\begin{align}\label{BC-three1}
&\pi(\pm H)=0\ \ \ \text{ with }\ \ \int_\Domain\!\!\pi \dd x=2v_\infty,\ \ \ \ \ 
\\[-.6em]&\label{BC-three2}
\alpha(\pm H)=1,\ \ \ \ \ \\&\theta(\pm H)=\theta_\infty\,.
\end{align}\end{subequations}

It is noteworthy that \eqref{eq4-three} decouples from (\ref{eq-three}a,c)
which arises not from necessity but rather from our desire for simplicity and for
consistency with the standard rate-and-state friction as in
Section~\ref{sec-intro}: we assumed that $\mu$, $f_0$, and $f_1$ are
independent of $\alpha$. The system (\ref{eq-three}a,c)--(\ref{BC-three}a,c)
thus represents a nonstandard non-local two-point boundary-value problem for
the functions $(\pi,\theta)$ on $\Domain$ and one scalar variable
$\sigma$. When solved, the two-point boundary-value problem
\eqref{eq4-three}--\eqref{BC-three2} can be solved for $\alpha$. Then $\eps$ is
obtained from \eqref{e=..alpha}. Eventually, the velocity $v$ can be calculated
from (\ref{eq2}) when using also \eqref{BC-three1}.

\subsection{The evolutionary model}
\label{su:EvolMod}
\def\pp{p}

We will now investigate an evolution version of the steady-state model \eqref{eq},
which in particular explains how \eqref{eq} have arisen.
In addition to the variables needed in Section~\ref{sec-steady},
we now will exploit also:
\\
\textbullet\ $\ \ \pp$ plastic strain (dimensionless) and 
\\
\textbullet\ $\ \ \varrho$ mass density (in one-dimensional model kg/m).\\
An additional ingredient will be 
a dissipation potential $\zeta$ for damage, which is convex with 
subdifferential $\partial\zeta$ and has physical
dimension J/m.

\noindent
The evolution variant of \eqref{eq} then looks as:
\begin{subequations}\label{evol}\begin{align}\label{evol1}
&\varrho\DT v-(
\bbC(\alpha)\eps)_x=0\,,&&\text{(momentum equilibrium)}
\\\label{evol2}
&\DT\eps+\DT{\pp}=v_x\,,
&&\text{(additive decomposition)}
\\&\label{evol3}
\pl_\pi R(\DT\pp,\theta) \ni 
\bbC(\alpha)\eps+\eta\DT{\pp}_{xx}\,,&&\text{(plastic flow rule)}
\\&\label{evol4}
\partial\zeta(\DT\alpha)
+\frac12\bbC'(\alpha)\eps^2+G_{\rm c}\frac{\alpha{-}1}{\ell^2}
\ni G_{\rm c}\ell^2\alpha_{xx}\,,&&\text{(damage flow rule)}
\\&\label{evol5}
\DT\theta=f_0(\theta)-|\DT{\pp}|f_1(\theta)
+\kappa\theta_{xx}\,.&&\text{(aging flow rule)}
\end{align}\end{subequations}
It is to be completed with
boundary conditions as \eqref{BC} with possibly time dependent
boundary velocity $v_\infty=v_\infty(t)$, i.e.\ here
\begin{align}\label{BC-evol}
v(\pm H)=\pm v_\infty(t),\ \ \ \ \
    \pp(\pm H)=0,\ \ \ \ \
    \alpha(\pm H)=1,\ \ \ \ \
    \theta(\pm H)=\theta_\infty\,.
\end{align}
with $\theta_\infty$ constant in time.  The (Green-Naghdi's) additive
decomposition is written in rates, which just gives \eqref{evol2}.
Obviously, the steady-state variant of \eqref{evol} where all
time derivatives vanish yield just  \eqref{eq}.

The system (\ref{evol}a-d) has a rational physical background while \eqref{evol5}
expresses some extra phenomenology controlling the nonconservative part in
\eqref{evol3}. For $\varrho=0$, the system (\ref{evol}a--d) represents the
so-called Biot equation
$\partial_{\DT q}\mathcal{R}(q,\theta,\DT q)
+\partial_q\mathcal{E}(q,\theta)=0$ for the state $q=(u,\pp,\alpha)$ and
$\theta$ given with the total dissipation potential
$\mathcal{R}(q,\theta,\DT q)=\int_\Domain\zeta_{\rm
  tot}^{}(\theta,\alpha;\pi,\DT\alpha) \dd x$ and the stored energy
$\mathcal{E}(q,\theta)=\int_\Domain\psi(\eps,\alpha,\theta) \dd x$, while for
$\varrho>0$ it arises from the Hamilton variational principle generalized for
the dissipative systems with internal variables.

The underlying specific stored energy and the dissipation potential (in terms of
the rates of plastic strain $p$ and damage $\alpha$) behind
this model are
\begin{subequations}
\label{eq:energetics}
\begin{align}
&\varphi(\eps,\alpha)=\frac12\bbC(\alpha)\eps^2
+G_{\rm c}\Big(\frac{(1{-}\alpha)^2}{2\ell^2}
+\frac{\ell^2}2\alpha_x^2\Big)
\ \ \text{ and }\ \
\\&\zeta_{\rm tot}^{}(\theta;\DT{p},\DT\alpha)
=  R(\DT{p},\theta) + \zeta(\DT\alpha)+\frac\eta2{\DT p}_{x}^2\,,
\end{align}
\end{subequations}
where often $\bbC(\alpha)=(\ell^2/\ell_0^2+\alpha^2)C_0$ with some $\ell_0$.
The constants $\ell$ and $\ell_0$ are in meters while the
fracture toughness $G_{\rm c}$ is in J/m$^2$, cf.\
\cite[Eqn.\,(7.5.35)]{KruRou19MMCM}, or rather in J/m in our 1-dimensional model.
This is known as the Ambrosio-Tortorelli functional  \cite{AmbTor92AFDP}.

\section{Analysis of the steady state model} 
\label{se:AnaSteady}

Further on, we will use the standard notation for the function space. In
particular, $C(\Domain)$ will be the space of continuous functions on $\Domain$
and $L^p(\Domain)$ will denote the Lebesgue space of measurable functions
on the domain $\Domain=[-H,H]$ whose $p$-power is integrable (or, when $p=\infty$,
which are bounded), and $W^{k,p}(\Domain)$ the Sobolev space of functions in
$L^p(\Domain)$ whose $k$-th distributional derivative belongs to $L^p(\Domain)$.
We abbreviate $H^k(\Domain)=W^{k,2}(\Domain)$. Besides, $H_0^1(\Domain)$ will denote
a subspace of $H^1(\Domain)$ of functions with zero values at $x=\pm H$. In
Section~\ref{sec-evol}, for the time interval $I=[0,T]$ and a Banach space $X$, 
we will also use the Bochner spaces $L^p(I;X)$ of Bochner-measurable functions
$I\to X$ whose norm in in $L^p(I)$, and the Sobolev-Bochner space $H^1(I;X)$ which
belong, together with their distributional time derivative, into $L^p(I;X)$.

\subsection{Existence of steady states} 

Let us recall the standard definition of a weak solution to the inclusion
\eqref{eq3} as a variational inequality
\begin{align}\label{eq3-weak}
  \int_\Domain\big( R(\wt\pi,\theta)-\sigma(\wt\pi{-}\pi)
  +\eta\pi_x(\tilde\pi{-}\pi)_x \big)  \dd x\ge\int_\Domain
  R(\pi,\theta) \dd x
\end{align}
to be satisfied for any $\wt\pi\in L^1(\Domain)$, where $\pl_\pi R(\pi,\theta)
=\mu(\pi,\theta)\mathrm{Sign}(\pi)$. 
We will prove existence of solutions due to even a stronger
concept of a classical (also called Carath\'eodory or strong) solution, namely
that $|\pi_{xx}|$ is integrable (actually in our case even bounded) and
\begin{align}\label{eq3-strong}
\forall\,\wt\pi\in\R:\quad 
R(\wt\pi,\theta)-\sigma(\wt\pi{-}\pi)+\eta\pi_{xx}(\wt\pi{-}\pi)
\ge R(\pi,\theta)
\end{align}
holds a.e.\ on $\Domain$.
As mentioned in Section~\ref{sec-intro}, the rate-and-state friction model lacks
standard thermodynamical consistency, which is reflected in the steady-state  
case by a lack of joint variational structure. Nevertheless, the two
equations \eqref{eq3} and \eqref{eq5} for $\pi$ and $\theta$, respectively, have
an individual variational structure governed by
the functionals 
\begin{equation}
  \label{eq:calA.calB}
  \mathcal A_\pi(\theta) := \int_\Domain|\pi|\varphi_1(\theta)- 
\varphi_0(\theta) +\frac\kappa2|\theta_x|^2 \dd x  
\quad \text{and} \quad 
\mathcal B_\theta(\pi) := \int_\Domain R(\pi,\theta)  +\frac\eta2|\pi_x|^2 \dd x , 
\end{equation}
where $\varphi_0$ and $\varphi_1$ are primitive functions to $f_0$ and $f_1$,
respectively. Then, the pair $(\theta,\pi)$ is a desired solution if and only
if $\theta $ minimizes $\mathcal A_\pi(\cdot)$ on
$\{\theta \in H^1(D) ;\ \theta(\pm H)=\theta_\infty\}$ and $\pi$ minimizes
$\mathcal B_\theta(\cdot )$ on
$\{\pi \in H^1_0(D);\ \int_D \pi \dd x =2v_\infty\}$. Since both functionals
$\mathcal A_\pi(\cdot)$ and $\mathcal B_\theta(\cdot )$ are strictly convex,
the solutions operators $\theta=S_{\mathcal A}(\pi)= \text{argmin}\mathcal
A_\pi$ and $\pi =S_{\mathcal B}(\theta) =  \text{argmin}\mathcal  B_\theta(\cdot )$
are well-defined. 
 The existence of steady states will be proved by a Schauder fixed-point
 theorem  applied to $S_{\mathcal A} \circ S_{\mathcal B}$. 

\begin{theorem}[Existence of steady states]
\label{th:ExistSteady}
Let the following assumptions hold: 
\begin{subequations}\label{ass}\begin{align}\nonumber
 &\mu:\R^2\to\R  \text{ continuous, }\mu(\cdot,\theta)\text{ non-decreasing
   on } [0,+\infty)
\\&\hspace{10em}
\text{and non-increasing on } (-\infty,0],\quad \inf\nolimits_\R\mu(0,\theta)>0,
 \label{ass1}\\
 &\bbC:\R\to\R \text{ continuously differentiable, }\
\bbC'([1,\infty))=0,\ \  \inf\nolimits_\R \bbC(\alpha)>0 ,
 \\\nonumber&f_0,f_1  \text{ continuous, non-negative, } 
 f'_1(\theta)>0, \ \ f_1(0)=0,
\\&\label{ass2}
\hspace{13.8em}  f'_0(\theta)< 0, \  \ f_0(\theta_\infty)=0,
 \\&\kappa>0,\ \ \ell>0, \ \ \eta>0  \,.
\end{align}
\end{subequations}
Then:\\
\ITEM{(i)}{For all $v_\infty \in \R$, problem \eqref{eq}--\eqref{BC} has a solution in the
classical sense (i.e.\ {\rm(\ref{eq}a,b,d,e)} hold everywhere 
and \eqref{eq3-strong} holds a.e.\ on $\Domain$) such that
$\eps\in W^{1,\infty}(\Domain)$, $v\in W^{3,\infty}(\Domain)$, and
$\pi,\alpha,\theta\in W^{2,\infty}(\Domain)$.}
\ITEM{(ii)}{Moreover, any solution
satisfied $0\le\theta\le\theta_\infty$ and 
$0\le\alpha\le1$ with $\alpha$ convex.}
\ITEM{(iii)}{If $v_\infty\ne0$, then $\sigma v_\infty>0$ with
$\sigma=\bbC(\alpha)\eps$ denoting the stress, and if also $\bbC'\le0$ with
$\bbC'(1)<0$, then $\alpha(x)<1$ except at $x=\pm H$.}
\ITEM{(iv)}{If $\bbC$, $f_0$, $f_1$, and
$\mu$ are smooth, then $\alpha,\theta\in W^{4,\infty}(\Domain)$.}
\end{theorem}
\begin{proof}
For a given $\wt\theta$, equation \eqref{eq3-three} with the nonlocal condition
in \eqref{BC-three} is equivalent to $\pi = S_\mathcal{B}(\wt\theta)=
\text{argmin} \mathcal B_{\wt\theta}(\cdot)$. 
The monotonicity of $\mu(\cdot,\wt\theta)$ assumed in \eqref{ass1} ensures the
uniform convexity of the functional $\mathcal B_\theta(\cdot)$. Therefore the
minimizer $\pi=S_\mathcal{B}(\wt\theta)$, which clearly exists by the direct
method in the calculus of variations, is uniquely determined. Moreover, it
depends depends continuously on $\wt\theta$ with respect to the weak topology
on $H^1(\Domain)$.  Thanks to \eqref{ass1}, for $v_\infty$ given,
$\mathcal B_\theta(\cdot)$ is coercive uniformly with respect to $\wt\theta$,
and therefore the minimizer $\pi= S_\mathcal{B}(\theta) $ can be a priori
bounded in $H^1(\Domain)$ independently on $\wt\theta$.

With a Lagrange multiplier $\sigma$ for the scalar-valued constraint
$\int_D \pi\dd x = 2v_\infty$, the Lagrangian for minimizing $\mathcal
B_{\wt\theta}$ reads
\begin{align}
  \label{eq:scrL}
\mathscr{L}(\pi,\sigma)=\int_\Domain R(\pi,\wt\theta)
+\frac\eta2\pi_x^2+\sigma\Big(\pi-\frac{v_\infty}H\Big) \dd x
\end{align}
and the optimality conditions $\pl_\pi\mathscr{L}(\pi,\sigma)\ni0$ and
$\pl_\sigma\mathscr{L}(\pi,\sigma)=0$ with ``$\pl$'' denoting the partial
subdifferentials (in the functional sense) give respectively the inclusion
\eqref{eq3-three} with $\wt\theta$ instead of $\theta$ and the integral
condition $\int_\Domain\pi \dd x=2v_\infty$ in \eqref{BC-three}. Also this
multiplier is determined uniquely and depends continuously on $\wt\theta$. From
\eqref{eq3-three} written as
$\sigma\in\mu(\pi,\wt\theta){\rm Sign}(\pi)-\eta\pi_{xx}\in H_0^1(\Domain)^*$,
we can see that also $\sigma\in\R$ is a priori bounded independently of
$\wt\theta$.

For a given $\pi$,
equation \eqref{eq5} is equivalent to $\theta=
S_\mathcal{A}(\pi)=\text{argmin}\mathcal A_\pi(\cdot)$. 
As $f_1$ is nondecreasing and $f_0$ is
nonincreasing, the functional $\mathcal A_\pi(\cdot)$ is convex, and it is to
be minimized on the affine manifold 
$\{\theta\in H^1(\Domain);\ \theta(\pm H)=\theta_\infty\}$, cf.\ the boundary
conditions \eqref{BC-three}. Therefore this boundary-value problem has a unique weak
solution $\theta\in H^1(\Domain)$,
which depends continuously on $\pi$ and
can be bounded independently of $\wt\theta$ when taking into
account the mentioned a priori bound for $\pi$.

Using $f_1(0)=0$, $f_0(\theta_\infty)=0$, and $\theta(\pm H)=\theta_\infty$,
the maximum principle implies $0\le\theta\le\theta_\infty$.

Altogether, we obtain a mapping
$\wt\theta \mapsto \theta  =S_\mathcal{A}\big( S_\mathcal{B}(\theta)\big)$
which is continuous with respect to the weak topology on $H^1(\Domain)$ and
valued in some bounded set (depending possibly on a given $v_\infty$).  By the
Schauder fixed-point theorem, this mapping has a fixed point $\theta$. This
thus determines also $\pi=S_\mathcal{B}(\theta)$ and $\sigma$.

Having $\sigma$ determined, we can find a unique weak solution
$\alpha\in H^1(\Domain)$ to the equation \eqref{eq4-three} with the boundary
conditions \eqref{BC-three2} and then, from \eqref{e=..alpha}, we also obtain
$\eps\in H^1(\Domain)$. From $v(x)=\int_{-H}^x\pi(\wt x)\,\d\wt x$, we also
obtain $v\in W^{2,2}(\Domain)$.

The quadruple $(\pi,\alpha,\theta,\sigma)$ solves
\eqref{eq-three}--\eqref{BC-three} in the weak sense. By comparison, we can
also see that $\pi_{xx},\alpha_{xx},\theta_{xx}\in L^\infty(\Domain)$, so that
$\pi,\alpha,\theta\in W^{2,\infty}(\Domain)$.

If $v_\infty\ne0$, then necessarily $\sigma\ne0$. If also $\bbC'\le0$ with
$\bbC'(1)<0$, the (convex) solution $\alpha$ to \eqref{eq4-three} must be
nontrivial, this $\alpha<1$ except the end points $x=\pm H$.

Then, from \eqref{e=..alpha} with $\sigma$ already fixed and $\bbC(\cdot)$
smooth, we obtain $\eps\in W^{2,\infty}(\Domain)$.  Eventually
$v\in W^{3,\infty}(\Domain)$ can be reconstructed from \eqref{eq2} with the
boundary conditions \eqref{BC}; here we used the constraint
$\int_D \pi \dd x = 2 v_\infty$.
\end{proof}

\begin{figure}
\begin{tikzpicture}
\draw[thick,->] (-1.2,0) -- (1.2,0) node[above]{$x$};
\draw[thick,->] (0,-0.4) -- (0,2.3) node[pos=0.8,left]{$f$};
\draw[] (-1,2)--(-1,-0.1) node[below]{$-1$};
\draw[] (1,2)--(1,-0.1) node[below]{$+1$};
\draw[very thick, color=blue] (-1,1)--(0.5555,0.2222)--(1,2);
\draw[very thick, color=red!80!blue, domain=-1:1,samples=161] 
                             plot (\x, { 0.4-0.5* cos (540*\x) } ) ;  
\end{tikzpicture}
\quad
\begin{tikzpicture}
\draw[thick,->] (-1.2,0) -- (1.2,0) node[above]{$x$};
\draw[thick,->] (0,-0.4) -- (0,2.3) node[pos=0.8,left]{$f_\mathrm{dr}$};
\draw[] (-1,2)--(-1,-0.1) node[below]{$-1$};
\draw[] (1,2)--(1,-0.1) node[below]{$+1$};
\draw[very thick, color=blue]
      (-1,0.2222)--(-0.25,1)--(0,2)--(0.25,1)--(1,0.2222);
\draw[very thick, color=red!80!blue, domain=-1:1, samples=161] 
                                plot (\x, { 0.4+0.5*cos(180*\x) } );
\end{tikzpicture}
\quad
\begin{tikzpicture}
\draw[thick,->] (-1.2,0) -- (1.2,0) node[above]{$x$};
\draw[thick,->] (0,-0.4) -- (0,2.2) node[pos=0.8,left]{$f_\mathrm{ir}\!$};
\draw[] (-1,2)--(-1,-0.1) node[below]{$-1$};
\draw[] (1,2)--(1,-0.1) node[below]{$+1$};
\draw[very thick, color=blue]
      (-1,2)--(-0.75,1)--(0,0.2222)--(0.75,1)--(1,2);
\draw[very thick, color=red!80!blue, domain=-1:1, samples=161] 
                                 plot (\x, { 0.4 - 0.5*cos(180*\x) } );
\end{tikzpicture}
\hfill
\begin{minipage}[b]{0.33\textwidth}\caption{\small\sl Two examples of functions $f$ and
    their decreasing and increasing
  rearrangements $f_\mathrm{dr}$ and  $f_\mathrm{ir}$. \label{fig:rearrange}}
\end{minipage}
\end{figure}
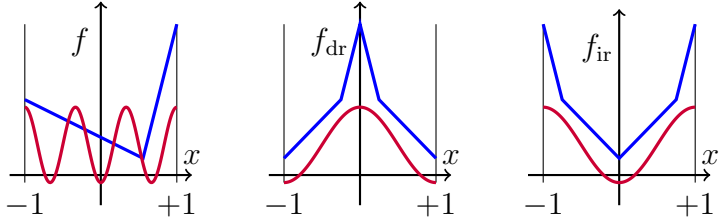

We discuss further qualitative properties of solution pairs $(\theta,\pi)$
that arise from the specific form of the steady state equations
\eqref{eq}--\eqref{BC}. As our above result does not imply uniqueness of
solutions, our next results states that there are solutions with symmetry and,
under a weak additional condition, these solutions are also monotone on
$[0,H]$.  For the latter we use the technique of rearrangements, which strongly
relies on the fact that we have no explicit $x$-dependence in our material
laws.  For general function $f \in L^1(D)$ we define its even decreasing and
even increasing rearrangements $f_\mathrm{dr}$ and $f_\mathrm{ir}$ via
\[
\{x\in D;\ f_\mathrm{dr}(x)>r\}= (-X(r),X(r)) \quad \text{where }
X(r):=\frac12\mathcal L^1\big( \{x\in D;\ f(x)>r\}\big) 
\] 
and $f_\mathrm{ir}(x) = f_\mathrm{dr}(H{-}|x|)$, see Figure
\ref{fig:rearrange}. 

The new condition \eqref{eq:mu.additive} for the following result is
satisfied in our adaptation \eqref{DR1++} of the classical Dieterich-Ruina
friction law \eqref{DR1}.

\begin{proposition}[Symmetric and monotone pairs]
\label{pr:SymMonotone}
Let the assumption \eqref{ass} of Theorem \ref{th:ExistSteady} hold. Then, for
all $v_\infty$ there exists an even solution pair $(\theta,\pi)$, i.e.\
$\theta$ and $\pi$ are even functions on $D=[-H,H]$.  If we additionally assume
\begin{equation}
  \label{eq:mu.additive}
  \mu(\pi,\theta) = \mu(\pi,0)+ B(\theta) \quad \text{with }
   B:\R\to [0,\infty) \text{ nondecreasing},
\end{equation}
then there exists an \emph{even, monotone pair} $(\theta,\pi)$, i.e.\ it is an
even pair such that additionally $[0,H]\ni x \mapsto \theta(x)$ is 
nondecreasing and $[0,H]\ni x \mapsto \pi(x)$ is nonincreasing.
\end{proposition}
\begin{proof} Throughout the proof we will restrict to the case $v_\infty>0$
  leading to $\sigma>0$ and $\pi\geq 0$. The case $v_\infty=0$ is trivial with
  $(\theta,\pi)\equiv (\theta_\infty,0)$, and $v_\infty<0$ follows similarly
  with $\sigma<0$ and $\pi\leq 0$.   

To obtain the evenness we simply restrict the existence theory
  developed in the proof of Theorem \ref{th:ExistSteady} to the closed
  subspaces of even functions. By the uniqueness of the minimizers of $\mathcal
  A_\pi$ and $\mathcal B_\theta$ it is clear that $S_\mathcal{A}$ and
  $S_\mathcal{B}$ map even functions to even functions. Hence, Schauder's
  fixed-point theorem produces an even solution. 

For showing the existence of monotone pairs we rely on classical results for
rearrangements, see e.g.\ \cite{Kawo85RCLS}, namely the Polya-Szeg\"o
inequality  
\begin{equation}
  \label{eq:PoSz}
  \int_D  (f_\mathrm{dr})_x ^2 \dd x =  \int_D  (f_\mathrm{ir})_x ^2 \dd x \leq
  \int_D f_x^2 \dd x
\end{equation}
and the Hardy-Littlewood inequality (cf.\ \cite[Ch.\,10]{HaLiPo34I})
\begin{equation}
\label{eq:HardyLittle}
\int_D f_\mathrm{dr}\,g_\mathrm{ir} \dd x= \int_D f_\mathrm{ir}\,g_\mathrm{dr} \dd
x \leq \int_D f\,g \dd  x \leq \int_D f_\mathrm{dr}\,g_\mathrm{dr} \dd x = \int_D
f_\mathrm{ir}\,g_\mathrm{ir} \dd x. 
\end{equation} 
While the upper estimate is classical and works for integration over $D=B_R(0)
\subset \R^d$ or $D=\R^d$, the lower estimate is special to $D\subset \R^1$,
see \cite[Eqn.\,(10.2.1)]{HaLiPo34I}. 

To exploit the theory of rearrangements we define the closed convex sets 
\begin{align*}
&\boldsymbol\Theta_\mathrm{ir}:=\big\{ \: \theta\in H^1(D); \ \theta(x)\in
   [0,\theta_\infty],   \ \theta(\pm H) =\theta_\infty, 
   \ \theta = \theta_\mathrm{ir} \: \big\} \quad \text{and}
\\
&\boldsymbol\Pi_\mathrm{dr}:=\big\{ \: \pi\in H^1(D); \ \pi(x) \geq 0, \ \pi(\pm
H)=0, \ \pi = \pi_\mathrm{dr}, \ \textstyle \int_D \pi\dd x = 2v_\infty \: \big\}
\end{align*}
and show below the mapping properties $S_\mathcal{A} : \boldsymbol\Pi_\mathrm{dr} \to
\boldsymbol\Theta_\mathrm{ir}$ and    $S_\mathcal{B} : \boldsymbol\Theta_\mathrm{ir}
\to \boldsymbol\Pi_\mathrm{dr} $. Thus, Schauder's fixed-point theorem can be
restricted to $S_\mathcal{A} \circ S_\mathcal{B} :
\boldsymbol\Theta_\mathrm{ir} \to \boldsymbol\Theta_\mathrm{ir}$ resulting in a
fixed point $\theta^* \in \boldsymbol\Theta_\mathrm{ir}$. With $\pi^*
=S_\mathcal{B}(  \theta^*)$, we obtain the desired even, monotone solution pair
$(\theta^*,\pi^*)$, namely $\theta^*=\theta^*_\mathrm{ir}$ and $\pi^*=\pi_\mathrm{dr}$. 

To establish $S_\mathcal{A} : \boldsymbol\Pi_\mathrm{dr} \to 
\boldsymbol\Theta_\mathrm{ir}$, we start with $\pi \in 
\boldsymbol\Pi_\mathrm{dr}$ and show $\mathcal A_\pi(\theta_\mathrm{dr}) 
\leq \mathcal A_\pi(\theta)$ for all $\theta \in H^1(D)$.  As $\theta=
S_\mathcal{A}(\pi)$ is the unique minimizer of $\mathcal A_\pi(\cdot)$, we
obtain $\theta= \theta_\mathrm{dr}$ as desired. 

To show $\mathcal A_\pi(\theta_\mathrm{dr}) \leq \mathcal A_\pi(\theta)$, we
exploit $|\pi|= \pi = \pi_\mathrm{dr}$ and the rearrangements estimates
\eqref{eq:PoSz} and \eqref{eq:HardyLittle} to obtain
\begin{align*}
&\int_D \theta_x^2 \dd x \!\overset{\text{\eqref{eq:PoSz}}}\geq\! \int_D \big(\theta_\mathrm{ir}\big){}_x^2 \dd x , \qquad 
\int_D \varphi_0(\theta) \dd x =  \int_D \varphi_0\big(\theta_\mathrm{ir} \big)
 \dd x ,\\ 
& \int_D |\pi|\,\varphi_1(\theta) \dd x 
  \!\overset{\text{\eqref{eq:HardyLittle}}}\geq\!   \int_D
  \pi_\mathrm{dr}\, \big( \varphi_1(\theta)\big)_\mathrm{dr}   \dd x = \int_D
  |\pi|\,\varphi_1\big(\theta_\mathrm{dr}\big)  \dd x .
\end{align*}
For the last identity we use $\big( \varphi_1(\theta)\big)_\mathrm{dr} =
\varphi_1(\theta_\mathrm{dr})$ which holds because of
$\varphi'_1=f_1(\theta) \geq 0$. Summing the three relations gives $\mathcal
A_\pi(\theta_\mathrm{dr}) \leq \mathcal A_\pi(\theta)$. 

Similarly, we derive $S_\mathcal{B} : \boldsymbol\Theta_\mathrm{ir} \to
\boldsymbol\Pi_\mathrm{dr}$ from $\mathcal B_{\theta}(\pi_\mathrm{dr}) \leq
\mathcal B_{\theta}(\pi)$ if $\theta \in \boldsymbol\Theta_\mathrm{ir}$. For
this we use assumption \eqref{eq:mu.additive}, which gives $R(\pi,\theta)=
R(\pi,0)+ B(\theta)|\pi|$, and the three relations 
\begin{align*}
&\int_D \pi_x^2 \dd x \!\overset{\text{\eqref{eq:PoSz}}}\geq\! \int_D 
   \big(\pi_\mathrm{dr}\big){}_x^2 \dd x , \qquad 
\int_D R(\pi,0) \dd x =  \int_D R(\pi_\mathrm{dr},0) \dd x ,\\ 
& \int_D |\pi|\,B(\theta) \dd x 
  \!\overset{\text{\eqref{eq:HardyLittle}}}\geq\!   \int_D
  \pi_\mathrm{dr}\, \big( B(\theta)\big)_\mathrm{dr}   \dd x = \int_D
  |\pi|\, B\big(\theta_\mathrm{dr}\big)  \dd x ,
\end{align*}
where we used that $B$ is nondecreasing. 

This finishes the proof of existence of even, monotone pairs. 
\end{proof}
\EEE

\begin{remark}[{\sl Aseismic-slip regime}]\label{rem-aseismic}
  Under very low shear velocities $|v_\infty|\ll 1$, real faults may go into
  so-called aseismic slip (also called aseismic creep), where one observes pure
  sliding like predicted by our steady state solutions constructed
  above. However, for our simplified evolutionary model introduced in Section
  \ref{se:NumSimul} (cf.\ \eqref{eq:SimpMod}) numerical simulations predict
  instability of the steady state and the development of stick-slip
  oscillations, see Section~\ref{su:NumSimPDE}. In the former case, stresses
  remain low and never challenge the plastic yield stress
  $\mu(0,\theta_\infty)$ at the core of the faults, a fact which is
  unfortunately not covered by our model.  One possible modification for
  modeling this effect would be to replace the set-valued Sign$(\cdot)$ in
  \eqref{eq3} by some monotone smooth approximation, e.g.\
  $\pi \mapsto \tanh(\pi/\delta)$ with $0<\delta\ll 1$.
\end{remark}

\subsection{Asymptotics of the plastic zone for $\eta\to0$ 
  and $\kappa\to 0$}
\label{sec-localization}
The gradient term in \eqref{eq3} and in \eqref{eq3-three} controls in a certain
way the width of the cataclastic zone where the slip is concentrated. There is
an expectation that, when suppressing it by $\eta\to0$, the slip zone will get
narrower. It is however a rather contra-intuitive effect that the zone
eventually does not degenerate to a completely flat interface like it would be
in so-called perfect plasticity where the plastic strain rate $\pi$ would be a
measure on $\Domain$.  Here, in the limit, $\pi$ only looses its
$W^{2,\infty}$-regularity as stated in Theorem~\ref{th:ExistSteady} for
$\eta>0$ but remains in $L^1(\Domain)$.

The definition of weak solutions \eqref{eq3-weak}
remains in its variational form or in its
strong form \eqref{eq3-strong} just putting $\eta=0$.
It should be emphasized that the boundary conditions $\pi(\pm H)=0$ are now
omitted. It will turn out that in the limit $\eta=0$  the plastic variable
$\pi$ becomes a pointwise function of $\theta$ and $\sigma$. 
By the strict convexity of $\pi\mapsto R(\pi,\theta)$ the set-valued mapping
$\pi\mapsto \pl_\pi(\pi,\theta)=\mu(\pi,\theta){\rm Sign}(\pi)$ is strictly
monotone (cf.\ \eqref{ass1}). Thus, $\pi$ in 
$\mu(\pi,\theta)\mathrm{Sign}(\pi)\in\sigma$ 
can be uniquely determined as a function of $\sigma$ and $\theta$. 
Specifically,
\begin{align}
  \label{pi=f(theta,sigma)}
\pi=\big[\mu(\cdot,\theta)\mathrm{Sign}(\cdot)\big]^{-1}(\sigma)=:
\varPi(\sigma,\theta)\,,
\end{align}
and the mapping $\varPi:\R^2\to\R$ is continuous. 

In this section, let us denote the solution obtained as a
Schauder fixed point in the proof of Theorem~\ref{th:ExistSteady} by
$(\eps_\eta,v_\eta,\pi_\eta,\alpha_\eta,\theta_\eta,\sigma_\eta)$.

\begin{proposition}[Convergence for $\eta\to0$]\label{prop2}
Let assumptions \eqref{ass} hold together with
\begin{subequations}
  \label{ass-M}
 \begin{align}\label{ass-M1}
   & \exists\, \varPhi:\R\to [0,\infty) \text{ continuous, superlinear }\forall\,
   (\pi,\theta): \quad R(\pi,\theta) \geq  \varPhi(\pi)   \ \ \text{ and}
   \\\label{ass-M2}
   &\big|\mu(\pi,\theta){-}\mu(\pi,\wt\theta)\big|\le o\big(|\theta{-}\wt\theta|\big)
   \ \text{ with some $o:\R^+\to\R^+$ continuous, $o(0)=0$}\,.
 \end{align}
\end{subequations}
There is a subsequence such that, for some
 $\pi\in  L^1(\Domain)$, $v\in W^{1,1}(\Domain)$,
  $\alpha\in W^{2,\infty}(\Domain)$, $\eps\in W^{1,\infty}(\Domain)$, 
  $\theta\in W^{1,\infty}(\Domain)$, and $\sigma\in\R$, it holds
\begin{subequations}
  \label{conv}
 \begin{align}
  \label{conv1}
\mbox{}\qquad\qquad &\eps_\eta\to\eps&&\text{weakly* in }\
W^{2,\infty}(\Domain),&&\qquad\qquad\mbox{}
  \\
  &v_\eta\to v&&\text{weakly in }\  W^{1,1}(\Domain),&&
  \\
  &\pi_\eta\to\pi&&\text{weakly in }\ L^1(\Domain),&&
  \\
  &\alpha_\eta\to\alpha&&\text{weakly* in }\ W^{2,\infty}(\Domain),&&
  \\
  &\theta_\eta\to\theta&&\text{strongly in }\ H^1(\Domain),&&
  \\
  &\sigma_\eta\to\sigma&&\text{in }\ \R,\ \ \text{ and }
  \\
  \label{conv.g}
  & \pi(x) =\varPi(\sigma,\theta(x))\hspace{-2em}&& \text{ for a.a. } x\in D.&&&&
 \end{align}
\end{subequations}
Moreover, $(\eps,v,\pi,\alpha,\theta,\sigma)$
is a classical solution to \eqref{eq}--\eqref{BC} in the sense that {\rm(\ref{eq}a,b,d,e)}
and \eqref{eq3-strong} with $\eta=0$ hold pointwise everywhere on $\Domain$.
More specifically, $\pi\in C(\Domain)$ and $v\in C^1(\Domain)$.
\end{proposition}
\begin{proof}
  From the proof of Theorem~\ref{th:ExistSteady}, we can see that the a priori
  bounds for 
\[
(\eps_\eta,v_\eta,\pi_\eta,\alpha_\eta,\theta_\eta,\sigma_\eta) \in 
  W^{2,\infty}(\Domain) {\times} W^{1,1}(\Domain)^2{\times} L^1(\Domain){\times}
  W^{2,\infty}(\Domain) {\times}  W^{2,1}(\Domain)^2 {\times} \R
\] 
are independent of $\eta>0$ and
$\|\pi_\eta\|_{H^1(\Domain)}=\mathscr{O}(1/\sqrt\eta)$. Moreover, from
$\pi_\eta= S_\mathcal{B}(\theta_\eta)$, we can easily see that even
$R(\pi_\eta,\theta_\eta)$ is bounded in $L^1(\Domain)$. Using
\eqref{ass-M1} we can apply the criterion of  de la Valle\'e
Poussin \cite{ValPou15IL} and obtain that $\{\pi_\eta\}_{\eta>0}$ is weakly
compact in $L^1(\Domain)$.

Then the limit passage in the weak solution to \eqref{eq}--\eqref{BC}
for $\eta\to0$ is quite easy. The only nontrivial point is the limit passage in
the variational inequality \eqref{eq3-weak}. We first use $\eta
(\pi_\eta)_x = \mathscr{O}(\sqrt\eta)$ in $L^2(D)$ and obtain, for all
$\wt \pi \in H^1(D)$, the relations  
\begin{align}
  \nonumber
  \!\!\int_\Domain\! &
  R(\wt\pi,\theta)-\sigma(\wt\pi{-}\pi) \dd x
  =\lim_{\eta\to0}\int_\Domain R(\wt\pi,\theta_\eta)- \sigma_\eta \, (\wt\pi{-}\pi_\eta)
  +\eta(\pi_\eta)_x\tilde\pi_x \dd x
  \\  \nonumber
  &\geq
  \limsup_{\eta\to0}\int_\Domain
  R(\wt\pi,\theta_\eta)- \sigma_\eta  \,(\wt\pi{-}\pi_\eta)
  +\eta\,(\pi_\eta)_x(\tilde\pi{-}\pi_\eta)_x \dd x 
  \overset{\text{\eqref{eq3-weak}}}\geq\liminf_{\eta\to0}\int_\Domain
  R(\pi_\eta,\theta_\eta) \dd x 
  \\ 
  &\geq \liminf_{\eta\to0}\int_\Domain  R(\pi_\eta,\theta) \dd x
  +\lim_{\eta\to0}\int_\Domain\!R(\pi_\eta,\theta_\eta){-}R(\pi_\eta,\theta) \dd x
  \ge\!\int_\Domain\!
  R(\pi,\theta) \dd x\ + \ 0 \, .  
   \label{eq3-weak+}
\end{align}
The liminf estimate follows because $R(\cdot,\theta)$ is convex and
continuous such that $\int_\Domain R(\cdot,\theta) \dd x $ is weakly lower
semicontinuous on $L^1(\Domain)$. The penultimate integral in
\eqref{eq3-weak+} converges to $0$ because 
$\theta_\eta\to\theta$ uniformly on $\Domain$ due to the compact embedding
$W^{2,1}(\Domain)\subset C(\Domain)$. Hence, 
$\lim_{\eta\to0}|\int_\Domain R(\pi_\eta,\theta_\eta)-R(\pi_\eta,\theta)| \dd x
\le\lim_{\eta\to0}\int_\Domain|\pi_\eta|o(\theta_\eta{-}\theta) \dd x\le
\lim_{\eta\to0}\|\pi_\eta\|_{L^1(\Domain)}^{}o(\|\theta_\eta{-}\theta\|_{L^\infty(\Domain)}^{})
=0$ where the function $o$ is from \eqref{ass-M2}.  

The variational inequality \eqref{eq3-weak+} does not contain any
$x$-derivatives any more and hence is equivalent to the pointwise inequality
$R(\wt\pi, \theta(x)) - \sigma (\wt\pi{-}\pi(x)) \geq R(\pi(x),\theta(x)) $
a.e.\ in $D$. But this is equivalent to
$\sigma \in \pl_\pi R(\wt\pi(x),\theta(x))$ and hence \eqref{conv.g} holds.

Since the mapping $\varPi:\R^2\to\R$ from \eqref{pi=f(theta,sigma)} is
continuous and since $\theta\in H^1(\Domain)\subset C(\Domain)$, we see that
$x\mapsto \pi(x)=\varPi(\sigma,\theta(x))$ is continuous as well, i.e.\
$\pi\in C(\Domain)$.
\end{proof}

We are now ready to study the limit $\kappa\to 0$ as well, which is really
surprising because we are losing all control over spatial derivatives and all
the modeling length scales induced by $\eta$ and $\kappa$ tend to $0$. In such
a situation the usual compactness arguments fail and fast spatial oscillations,
i.e.\ microstructures, may appear. Indeed we will see in Remark
\ref{re:ManySols} that there are indeed many complicated solutions without any
length scale. However, it is surprising that it is possible to show that
natural solutions exist, namely even, monotone pairs $(\theta,\pi)$. The idea
is to use for $\kappa>0$ and $\eta=0$ the even, monotone pairs
$(\theta^\kappa,\pi^\kappa)$ obtained from Proposition \ref{pr:SymMonotone} and
the subsequent limit $\eta\to 0$ in Proposition \ref{prop2}. The monotonicity
of the pairs $(\theta^\kappa,\pi^\kappa)$ allows us to deduce pointwise
convergence, which is good enough to pass to the limit $\kappa\to 0$ even in
nonlinear functions.
  
Under the additional assumptions \eqref{eq:wtmu.monotone}, which are satisfied
by our example treated in Section \ref{su:SimplifMod}, we then obtain the
typical behavior. There is a critical value $\pi_*>0$ such that for small
positive $v_\infty$ the cataclastic zone is $(-h,h)$ with $h=v_\infty/\pi_*$,
where $(\theta,\pi)$ assume constant values $(\theta_*,\pi_*)$ independent of
$v_\infty$, whereas for $x$ with $h<|x|<H$ we have
$(\theta,\pi)=(\theta_\infty,0)$, see \eqref{eq:PlateauSol}.

\begin{proposition}[The limit $\kappa\to 0$ for monotone pairs]
\label{prop3} 
Let the assumptions \eqref{ass}, \eqref{eq:mu.additive}, and \eqref{ass-M} hold
and let us consider a family $\big( (\theta^\kappa,\pi^\kappa)\big)_{\kappa>0}$ of
even, monotone solutions to \eqref{eq} with $\eta=0$ and $v_\infty>0$.  
Then:\\
\ITEM{(i)}{there exists a subsequence (not relabeled) and an even, monotone pair
$(\theta^0,\pi^0) \in L^\infty(D)\times  L^\infty(D) $ such that for $\kappa
\to 0$ we have the convergence}
\[
(\theta^\kappa(x),\pi^\kappa(x)) \ \to \ (\theta^0(x),\pi^0(x)) \quad
  \text{ for a.a. }x\in D 
\]
\ITEM{}{and that $(\theta^0,\pi^0)$ solves the minimization problems}
\begin{equation}
  \label{eq:calA.B.0}
  \mathcal A_{\pi^0}^0(\theta^0) \leq \mathcal A_{\pi^0}(\theta):=\!\int_D
|\pi^0| \varphi_0(\theta){-}\varphi_0(\theta) \dd x  
\  \text{ and } \ 
\mathcal B_{\theta^0}^0(\pi^0) \leq \mathcal B_{\theta^0}(\pi):=\!\int_D\!
R(\pi,\theta^0)\dd x 
\end{equation}
\ITEM{}{for all $(\theta,\pi)\in L^1(D){\times} L^1(D)$ with
$\int_D \pi \dd x = 2v_\infty$.}
\ITEM{(ii)}{Moreover, if we define $\theta=\Theta_f(\pi)$ to be the unique solution of
$f_0(\theta)=|\pi| f_1(\theta)$, set $\wt\mu:[0,\infty)\to
(0,\infty); \ \pi\mapsto \mu(\pi,\Theta_f(\pi))$, and  
assume that there exists $\pi_\circ>0$ such that}
\begin{equation}
  \label{eq:wtmu.monotone}
  \wt\mu \text{ is strictly decreasing on }[0,\pi_\circ]
           \quad \text{and} \quad 
  \wt\mu \text{ is strictly increasing on } [\pi_\circ,\infty),
\end{equation}
\ITEM{}{then there exists a unique $\pi_*>\pi_\circ$ such
$\int_0^{\pi_*}\wt\mu(\pi)\dd\pi= \pi_* \wt\mu(\pi_*)$ and the above solutions
$(\theta^0,\pi^0)$ are uniquely given by}
\begin{equation}
  \label{eq:PlateauSol}
  (\theta^0,\pi^0)(x)
= \begin{cases}  (\Theta_f(\pi_*),\pi_*) & \text{for } |x| <
            {v_\infty}/{\pi_*} \leq H, \\
  (\theta_\infty,0) & \text{for }
  {v_\infty}/{\pi_*} <|x| \leq H, \\
  \big(\Theta_f(v_\infty/H), v_\infty/H\big) & \text{for } v_\infty\geq \pi_* H.
  \end{cases}
\end{equation}
\ITEM{}{In particular, in this case the whole family $\big(
(\theta^\kappa,\pi^\kappa)\big)_{\kappa>0}$ converges pointwise.}
\end{proposition}
\begin{proof}
  By Proposition \ref{pr:SymMonotone} and Proposition \ref{prop2} we know that
  for all $\kappa>0$ even, monotone pairs $(\theta^\kappa,\pi^\kappa)$ exist
  and satisfy $\theta^\kappa \in W^1(D)$ and $\pi^\kappa \in C(D)$. Moreover,
  we have $\theta^\kappa(x)\in [0,\theta_\infty]$ and
  $\pi^\kappa (x) = \varPi(\sigma^\kappa, \theta^\kappa(x))$ for all $x\in D$. 

  \emph{Step 1. Superlinear a priori bound for $\pi^\kappa$:} We 
  again use the uniform superlinearity of the dissipation potential
  $R(\cdot,\theta)$ from  \eqref{ass-M1}.  As
  $\pi^\kappa$ is a minimizer of $\mathcal B_{\theta^\kappa} (\cdot)$ we obtain
  the uniform bound  $\int_D \varPhi(\pi_\kappa) \dd x \leq C_*< \infty$.
  Thus, we have weak compactness (by de la Valle\'e Poussin \cite{ValPou15IL})
  and along a subsequence (not relabeled) we have $\pi^\kappa \weak \pi^0$ and
  conclude $\int_D \pi^0 \dd x = 2 v_\infty$.  Moreover, using
  $\pi^\kappa = \pi_\mathrm{dr}^\kappa$ this implies the a priori bound
\begin{equation}
  \label{eq:Phi.pw.bdd}
  0 \leq \pi^\kappa (x) \leq R \quad \text{for } |x| \geq \frac{C_*}{\varPhi(R)} .
\end{equation}

\emph{Step 2. Pointwise convergence:} Exploiting the monotonicity and the a
priori bounds $\theta^\kappa\in [0,\theta_\infty]$ and \eqref{eq:Phi.pw.bdd},
we can apply the classical Helly's selection principle to obtain pointwise
convergence (everywhere in $D$). Along a subsequence (not relabeled) we have 
\[
\sigma^\kappa\to \sigma^0, \qquad (\theta^\kappa(x),\pi^\kappa(x)) \to
(\theta^0(x),\pi^0(x)) \text{ for all }x \in D. 
\]
Here the monotonicities are kept, i.e.\ $\theta^0=\theta_\mathrm{ir}$ and
$\pi^0= \pi^0_\mathrm{dr}$, but the continuity of the limits might be
lost. Moreover, $\pi^0(0)=\infty$ might be possible.

\emph{Step 3. Limit passage in the equations:} Since $\varPi$ is continuous,
the pointwise convergence yields the limit relation
\begin{equation}
  \label{eq:varPi.0}
  \pi^0(x)= \varPi(\sigma^0,\theta^0(x)) \quad \text{for all } x \in D.
\end{equation}
For the equation determining $\theta$ we can use the a priori estimate $ \kappa
\| \theta^\kappa\|_{L^2}^2 \leq C_*$ and pass to the limit in the weak form of
$(\kappa \,\theta^\kappa_x)_x+f_0(\theta^\kappa)=\pi^\kappa f_1(\theta^\kappa)$, i.e.\ in
the integral identity 
\[
\int_D \kappa \,\theta^\kappa_x \,\wt\theta_x - f_0(\theta^\kappa)\, \wt\theta +
\pi^\kappa f_1(\theta^\kappa)\, \wt\theta \dd x = 0 \quad \text{for all }
\wt\mu \in H^1_0(D). 
\]  
This provides the pointwise relation
\begin{equation}
  \label{eq:Thetaf.0}
   f_0(\theta^0(x)) = \pi^0(x)\, f_1(\theta^0(x)) \quad 
 \text{ for a.a.\ }x\in D.  
\end{equation}
From \eqref{eq:varPi.0} and \eqref{eq:Thetaf.0} we immediately see that
\eqref{eq:calA.B.0} holds.  

We next observe that $\theta=\Theta_f(\pi)$ is well-defined by the implicit
function theorem using \eqref{ass2}. Thus, the solutions satisfy
$\theta^0(x)=\Theta_f(\pi^0(x))$ for a.a.\ $x \in D$. Henceforth,
recalling
$\wt\mu(\pi)= \mu(\pi,\Theta_f(\pi))$,  the minimization
problem \eqref{eq:calA.B.0} is equivalent to $\sigma \in
\wt\mu(\pi)\mathrm{Sign}(\pi)$ and $\int_D \pi\dd x = 2v_\infty$. 
Defining the function $\sfR(\pi)=\int_0^\pi \wt\mu(s) \dd s $, this is
equivalent to the following problem: 
\[
\text{minimize }\ \pi \mapsto \int_D \sfR(\pi(x)) \dd x \quad \text{ subject to
} \ \pi\geq 0 \text{ and } \int_D \pi \dd x = 2v_\infty >0.
\]
However, this minimization problem is well understood via the convex hull
$\sfR^{**}$, see \cite[Ch.\,2]{Brai02GCB}. By our assumption
\eqref{eq:wtmu.monotone} we know that $\sfR^{**}$ has the form  
\begin{align}
  \label{eq:R**}
&\sfR^{**}(\pi)= \begin{cases} \sfR(\pi_*) \pi/\pi_* &
  \text{for } \pi \in [0,\pi_*], \\
 \sfR(\pi) & \text{for } \pi \geq \pi_*, \end{cases}. 
\end{align}
and satisfies $\sfR^{**}(\pi) \lneqq \sfR(\pi)$  for $\pi\in (0,\pi_*)$ and $
\sfR''(\pi)>0$ for $\pi \geq\pi_*$, see Figure \ref{fig:R**}. \EEE
\begin{figure}
\begin{tikzpicture}
\draw[->, thick] (0,0)--(4,0) node[below]{$\pi$};
\draw[->, thick] (0,0)--(0,3);
\draw (1.5,0.1)--(1.5,-0.15) node[below]{$\pi_\circ$};
\draw (2.2,0.1)--(2.2,-0.15) node[below]{$\pi_*$};
\draw[color =blue, very thick, domain=0.0:3.8] plot (\x, {0.4+ 0.4*(\x-1.5)^2})
 node[pos=0.8, left]{$\widetilde{\mu}(\pi)$};
\draw[color=blue!50!red, very thick, dashed] (0,0.62)--(2.2,0.62);
\end{tikzpicture}
\quad
\begin{tikzpicture}
\draw[->, thick] (0,0)--(4,0) node[below]{$\pi$};
\draw[->, thick] (0,0)--(0,3);
\draw (1.5,0.1)--(1.5,-0.15) node[below]{$\pi_\circ$};
\draw (2.2,0.1)--(2.2,-0.15) node[below]{$\pi_*$};
\draw[color =blue, very thick, domain=0.0:3.5] 
   plot (\x, {0.4*\x+ 0.1333*(1.5^3+(\x-1.5)^3)}) node[left]{$\sfR(\pi)$};
\draw[color=blue!50!red, very thick, dashed] (0,0)--(2.2,1.36) node[pos=0.8, below right]{$\sfR^{**}(\pi)$} ;
\end{tikzpicture}
\quad
\begin{minipage}[b]{0.29\textwidth}\caption{\small\sl The functions $\wt\mu$, $\sfR$, and
    $\sfR^{**}$.}  
\label{fig:R**}
\end{minipage}
\end{figure}

As our $\sfR$ is superlinear, a minimizer always exists. Moreover,
recalling that $v_\infty/H>0$ is the average value of $\pi:D\to \R$, the
minimizer is unique if and only if the tangent at $\pi=v_\infty/H$ is not in
the interior of an interval on which $\sfR^{**} $ is affine. In the open
interval $(0, v_\infty/H)$ the minimizers $\pi$ attain only the values $0$ and
$\pi_*$ on sets with the corresponding measures to fit the average.  However,
by constructing the even, nonincreasing rearrangement, we find a unique
minimizer, where only the value at the two jump points
$x = \pm h= v_\infty/\pi$ are free.

From these uniqueness results we also obtain the convergence of the full family
by the standard contradiction via compactness. With this, Proposition
\ref{prop3} is established. 
\end{proof}

The new condition \eqref{eq:wtmu.monotone} can be checked numerically for our
example specified in \eqref{eq:DataSimple} giving $\pi_* \approx
1.4923$ and $\pi_\circ=0.6193$. Indeed, to see the desired effect of a fixed
$\pi_*$ leading to a 
cataclastic zone of width
$2h=2v_\infty/\pi_*$, our condition \eqref{eq:wtmu.monotone} is sufficient, but
far from being necessary. What we really need is that $\sfR^{**}$ is affine in
an interval $[0,\pi_*]$, which automatically follows if $\sfR''(0^+)=\lim_{pi
  \searrow 0} \sfR''(\pi) < 0$. In fact, in general we can
consider the case $\mu(\pi,\theta)= \mu_0+A(\pi)+ B(\theta)$ and general $f_0$
and $f_1$. Using $\Theta_f(0)=\theta_\infty$ following from
$f_0(\theta_\infty)=0$, an explicit calculation gives
\[
\sfR''(0^+) =\wt\mu'(0^+) = \pl_\pi \mu(0^+,\theta_\infty) + \pl_\theta
\mu(0^+,\theta_\infty) \frac{f_1(\theta_\infty)}{f'_0(\theta_\infty)},
\]
which may be negative because of $f'_0(\theta_\infty)<0$.  
 
\begin{remark}[{\sl Nonuniqueness of solutions}]
\label{re:ManySols} 
We want to emphasize that the uniqueness result for $\kappa=\eta=0$ at the end
of Proposition \ref{prop3} concerns only even, monotone solutions. Because of
$\kappa=\eta=0$ there are indeed infinitely many solutions, as we can
``rearrange'' the function values of $(\theta,\pi)$ freely. In the case
$v_\infty< \pi_* H$, we can choose any open set $P\subset D$ with
$\int_D 1_P\dd x = 2v_\infty/\pi_*$ and the function
\[
\big(\theta(x),\pi(x)\big) = \begin{cases}
\big(\Theta_f(\pi_*),\pi_*\big)& \text{ for } x \in P,\\
\big(\theta_\infty,0\big)& \text{ for } x \in D\setminus P
\end{cases}
\]
is a solution of \eqref{eq:calA.B.0} as well. 
\end{remark} 
\EEE

\section{Analysis of the evolutionary model}\label{sec-evol}

We now consider the evolutionary model \eqref{evol}. 
The energetics \eqref{eq:energetics} behind this model can be revealed by
testing momentum balance \eqref{evol1} by $\wt v=v-w^\infty$ with
$w^\infty(t,x):=v_\infty(t) x /H$, the plastic flow rule \eqref{evol2} by
$\DT\pp$, and the damage rule \eqref{evol3} by $\DT\alpha$. 
Using the Dirichlet boundary condition for the velocity at $x=\pm H$, we
have $\wt v(\pm H)=0$, as needed. The first test gives, in particular, the term
\begin{align}\nonumber
\int_\Domain
\bbC(\alpha)\eps\Big(v_x-\frac{v_\infty}H\Big) \dd x=\int_\Domain\bbC(\alpha)\eps
(\DT\eps{+}\DT\pp) \dd x-\frac{v_\infty}H\int_\Domain\bbC(\alpha)\eps \dd x\qquad\qquad
\\=\frac{\d}{\d t}\int_\Domain\frac12\bbC(\alpha)\eps^2\dd x 
+\int_\Domain\bbC(\alpha)\eps\DT\pp-\frac12\bbC'(\alpha)\eps^2\DT\alpha \dd x
-\frac{v_\infty}H\int_\Domain\bbC(\alpha)\eps \dd x\,,
\end{align}
where also \eqref{evol2} has been used. This test of the inertial form gives
\begin{align}\nonumber
\int_\Domain\varrho\DT v\Big(v-v_\infty\frac{x}H\Big) \dd x=
\frac{\d}{\d t}\int_\Domain\frac\varrho2v^2 \dd x
-v_\infty\int_\Domain\varrho\DT v \frac{x}H \dd x\,.
\end{align}
Combining it with the tests of \eqref{evol2} by $\DT\pp$ and of \eqref{evol3} by
$\DT\alpha$ which give
\begin{subequations}\begin{align}\label{subst}
&\!\!\int_\Domain\bbC(\alpha)\eps\DT\pp \dd x=
\int_\Domain\mu(\DT\pp,\theta)|\DT\pp|+\eta\DT\pp_x^2 \dd x\ \ \text{ and}
\\&
\!\!\int_\Domain\!\!-\frac12\bbC'(\alpha)\eps^2\DT\alpha \dd x=
\!\int_\Domain\!\!\DT\alpha\pl\zeta(\DT\alpha)+\Big(\frac12\bbC'(\alpha)\eps^2\!
+G_{\rm c}\frac{\alpha{-}1}{\ell^2}\Big)\DT\alpha \dd x
+\frac{\d}{\d t}\int_\Domain\frac12G_{\rm c}\ell^2\alpha_x^2 \dd x,
\end{align}\end{subequations}
we altogether obtain the energy balance
\begin{align}\nonumber
&\frac{\d}{\d t}\int_\Domain\!\!\!\linesunder{\frac\varrho2v^2+\varphi(\eps,\alpha)+\frac12G_{\rm c}\ell^2\alpha_x^2}{kinetic and stored}{energies}\!\!\!\!\d x\!
\\[-1em]&\hspace{11em}
+\int_\Domain\!\!\!\lineunder{\mu(\DT\pp,\theta)|\DT\pp|+\DT\alpha\pl\zeta(\DT\alpha)+\eta\DT\pp_{x_{_{_{}}}}^2}{dissipation rate}\!\!\!\!\d x=\!\!\!\!\!\linesunder{\langle\tau,(v_{\infty_{_{}}},-v_\infty)\rangle}{power of}{external load}\!\!\!,
\label{engr}\end{align}
where $\tau\in\R^2$ is the traction on the boundary (i.e.\ here two
forces at $x=\pm H$) defined as a functional
$\langle\tau,(z(H),z(-H))\rangle
=\int_\Domain\varrho\DT v z+\bbC(\alpha)\eps z_x \dd x$ for any
$z\in H^1(\Domain)$, cf.\ e.g.\ \cite[Sect.6.2]{KruRou19MMCM}.

Further on, we will be interested in an initial-value problem. For this, we
prescribe some initial conditions, i.e.\ 
\begin{align}
  \label{IC}
v(\cdot,0)=v_0\,,\ \ \ \eps(\cdot,0)=\eps_0\,,\ \ \
\alpha(\cdot,0)=\alpha_0\,,\ \ \text{ and }\ \ 
\theta(\cdot,0)=\theta_0\,.
\end{align}

A definition of the weak solutions of particular equations/inclusions
in \eqref{evol} can be cast by standard way, using convexity of the
involved functionals. Let us specify, rather for illustration, the weak 
formulation for the inclusion \eqref{evol3} exploiting that 
$\mu(\DT\pp,\theta){\rm Sign}(\DT\pp)$, i.e.\
$\mu(\pi,\theta){\rm Sign}(\pi) = \pl_\pi R(\DT\pp,\theta)$ where
$R(\pi ,\theta)$  is convex in
the variable $\pi=\DT\pp$. This leads to the variational inequality
\begin{align}\label{eq3-weak+++}
  \int_0^T\!\!\int_\Domain
  R(\wt\pi,\theta)-\bbC(\alpha)\eps(\wt\pi{-}\DT\pp)
  -\eta\DT\pp_x(\tilde\pi{-}\DT\pp)_x
   \dd x\d t\ge\int_0^T\!\!\int_\Domain
  R(\DT\pp,\theta) \dd x\d t
\end{align}
to be valid for any $\wt\pi\in L^\infty(I{\times}\Omega)$.

Beside the previous assumptions, we now also assume
\begin{align}
v_0\in L^2(\Domain)\,,\ \ \ \eps_0\in L^2(\Domain)\,,\ \ \
\alpha_0\in H^1(\Domain)\,,\ \ \
\theta_0\in H^1(\Domain)\,.
\label{ass-IC}\end{align}
The definition of weak solutions to \eqref{evol} with \eqref{BC-evol}
and \eqref{IC} is standard and we will not write it explicitly; the
variational inequality \eqref{eq3-weak} is to hold 
integrated over $I$. Furthermore, we also exploit the superlinear growth of
$R(\cdot,\theta)$ from  \eqref{ass-M1}, namely
\begin{align}
   \label{mu-growth}
   \mu(\pi,\theta)|\pi| \geq R(\pi,\theta) \geq \varPhi(\pi) ,
\end{align}
which is a standard estimate for $\wt\mu \in \pl\psi(\pi)$, namely $\pi\wt\mu =
\psi(\pi)+\psi^*(\wt\mu)\geq \psi(\pi)$ as $\psi^*\geq 0$. 
Note that the standard model \eqref{DR1+} complies
with assumption \eqref{ass-M1}. 

Relying formally on the tests leading to \eqref{engr}, after integration
in time on the interval $[0,t]$ when using also the by-part
integration, we obtain 
\begin{align}\nonumber
&\int_\Domain\frac\varrho2v^2(t)+\varphi(\eps(t),\alpha(t))+\frac12G_{\rm c}\ell^2
\alpha_x^2(t) \dd x
+\int_0^t\!\!\int_\Domain
\mu(\DT\pp,\theta)|\DT\pp|+\DT\alpha\pl\zeta(\DT\alpha)+\eta\DT\pp_x^2 \dd x\d t
\\[-.3em]&\qquad\qquad\nonumber
=\int_\Domain\frac\varrho2v_0^2+\varphi(\eps_0,\alpha_0)
+\frac12G_{\rm c}\ell^2[\alpha_0]_x^2 \dd x
+\int_0^t\!\!\int_\Domain\varrho\DT v w^\infty+\bbC(\alpha)\eps_xw^\infty_x \dd x\d t
\\&\qquad\qquad\nonumber=\int_\Domain\frac\varrho2v_0^2+\varphi(\eps_0,\alpha_0)
+\frac12G_{\rm c}\ell^2[\alpha_0]_x^2+\varrho v(t)\big(v_\infty(t){-}v_\infty(0)\big)\frac{x}H \dd x
\\[-.3em]
&\hspace*{17em}+\!\int_0^t\!\!\int_\Domain\!\bbC(\alpha)\eps_x \frac{v_\infty}H
-\varrho  v\DT v_\infty\frac{x}H \dd x\d t.
\label{engr1}\end{align}

Moreover, the aging equation \eqref{evol5} has to be tested separately by using
the test function $\theta{-}\theta_\infty$, which has zero traces
for $x=\pm H$.  Integrating the result over $[0,t]$ leads to 
\begin{align}
 \nonumber
\int_\Domain \frac12\,\theta^2(t) \dd x+\int_0^t\!\!\int_\Domain\kappa\theta_x^2\dd x\dd t
&=\int_\Domain(\theta(t){-}\theta_0)\theta_\infty \dd x
\\[-.5em]&\quad+\int_0^t\!\!\int_\Domain|\DT\pp|f_1(\theta)(\theta{-}\theta_\infty)
-f_0(\theta)(\theta{-}\theta_\infty)\dd x\,.
\label{engr2}\end{align}

When summing \eqref{engr1} and \eqref{engr2}, we can use the H\"older and a
(generalized) Young inequality to estimate the resulting right-hand side.
Actually, the only nontrivial term is
$|\DT\pp|f_1(\theta)(\theta{-}\theta_\infty)$ in
\eqref{engr2} and it can be estimated as
\begin{align}
  \nonumber
  \int_\Domain|\DT\pp|f_1(\theta)(\theta{-}\theta_\infty) \dd x
  &\leq \int_\Domain \frac12\varPhi\big(|\DT\pp|\big) 
  +\frac12 \varPhi^*\big(2f_1(\theta)(\theta{-}\theta_\infty)\big) \dd x
  \\ \label{engr3}
  &\overset{\text{\eqref{mu-growth}}}\leq \int_\Domain\frac12\mu(\DT\pp,\theta)|\DT\pp|
  +\frac12\varPhi^*\big(2f_1(\theta)(\theta{-}\theta_\infty)\big) \dd x\,,
\end{align}
where $\varPhi^*$ is the Fenchel-Legendre conjugate of $\varPhi$, i.e.\
$\varPhi^*(s)=\sup_{\pi \in \R} \big(\pi s - \varPhi(\pi) \big) $.   \EEE

The term $\frac12\mu(\DT\pp,\theta)|\DT\pp|$ in \eqref{engr3} can then be
absorbed in the left-hand side of \eqref{engr1} while 
$\frac12 \varPhi^*(2f_1(\theta)(\theta{-}\theta_\infty))$ is a priori 
bounded since $0\le\theta\le\theta_\infty$. Eventually, the last term in
\eqref{engr1} can be  estimated as $\varrho(1{+}|v|^2)|\DT v_\infty|$.

Assuming $v_\infty\in W^{1,1}(I)$ and using Gronwall's inequality, from the
left-hand sides of  \eqref{engr1} and \eqref{engr2} we can read the a priori
estimates
\begin{subequations}
 \label{est}
 \begin{align}
  \label{est1}
  &
  \|v\|_{L^\infty(I;L^2(\Domain))}^{}\le C,
  \\&\|\eps\|_{L^\infty(I;L^2(\Domain))}^{}\le C,
  \\&\|\pp\|_{H^1(I;H^1(\Domain))}\le C,
  \\&\|\alpha\|_{L^\infty(I;H^1(\Domain))\,\cap\,H^1(I;L^2(\Domain))}^{}\le C,
  \\&\|\theta\|_{L^\infty(I;L^2(\Domain))\,\cap\,L^2(I;H^1(\Domain))}^{}\le C.
\label{est4}
 \end{align}
\end{subequations}
By comparison, we will get also an information about
$\DT v=(\bbC(\alpha)\eps)_x/\varrho\in L^\infty(I;H^1(\Domain)^*)$, about
$\DT\eps=v_x-\DT\pp\in L^2(I;H^1(\Domain)^*)$, and also
about $\DT\theta=f_0(\theta)-|\DT{\pp}|f_1(\theta)
+\kappa\theta_{xx}\in L^2(I;H^1(\Domain)^*)$.

The rigorous existence proof of weak solutions is however very nontrivial and
seems even impossible for the full dynamical model \eqref{evol} with damage.
Some modifications by involving some additional dissipative terms or some
higher-order conservative terms seem necessary, cf.\
\cite[Sect.7.5]{KruRou19MMCM} or also \cite{RoSoVo13MRLF} for the model
without aging. Consistently also with the computational experiments in
Section~\ref{se:NumSimul} below, we thus present the rigorous proof only
for a model without damage, i.e.\ for $\bbC>0$ constant.

\begin{theorem}[Damage-free
case -- existence and regularity of solutions]\label{th:EvolExist}
Let {\rm(\ref{ass}a,c,d)} with $\mu$ smooth, \eqref{ass-IC}, and \eqref{mu-growth}
hold, and $\varrho>0$ be a constant and $v_\infty\in W^{1,1}(I)$. Then:\\
\ITEM{(i)}{There is a weak solution $(v,\eps,\pp,\theta)
\in L^\infty(I;L^2(\Domain))^2\times 
H^1(I;H^1(\Domain))\times (L^\infty(I;L^2(\Domain))\cap L^2(I;H^1(\Domain)))$
to the initial-boundary-value problem for the system {\rm(\ref{evol}a-c,e)}
with the boundary conditions \eqref{BC-evol} and the initial conditions \eqref{IC}.}
\ITEM{(ii)}{If $\sup_{0\le\theta\le\theta_\infty}^{}\mu(\cdot,\theta)$ does not have
a growth more than $\mathscr{O}(|\pi|^s)$, then these solutions are,
in fact, regular in the sense that $\pp\in W^{1,s}(I;H^2(\Domain))$ and,
if $s\ge2$, also $\theta\in H^1(I;L^2(\Domain))
\cap L^\infty(I;H^1(\Domain))\cap L^2(I;H^2(\Domain))$ and also each such weak
solution satisfies the energy balance \eqref{engr} without $\alpha$-terms integrated
over a time interval $[0,t]$ with any $t\in I$.}
\end{theorem}

Let us note that the $\mathscr{O}(|\pi|^s)$-growth condition in the point (ii)
surely covers the model \eqref{DR1++} for any $1\le s<\infty$.

\begin{proof}[Sketch of the proof]
  Actually, the above formal procedure is to be made first for a suitable
  approximation whose solutions exist by some specific arguments, and then to
  pass to the limit. Imitating the split for the static problem used in the
  proof of Theorem~\ref{th:ExistSteady}, we choose a staggered time
  discretization. We take an equidistant partition of the time interval $I$ by
  using the time step $\tau>0$, assuming $T/\tau$ integer and considering a
  sequence of such $\tau$'s converging to 0. Then, recalling
  $\pl_\pi R(\pi,\theta)=\mu(\pi,\theta){\rm Sign}(\pi)$, we consider a
  recursive boundary-value problem for the system
\begin{subequations}
  \label{disc}
 \begin{align}
   \label{disc1}
&\varrho\frac{v_\tau^k-v_\tau^{k-1}\!\!}\tau
-(\bbC\eps_\tau^k)_x=0\,,
\\&\frac{\eps_\tau^k-\eps_\tau^{k-1}\!\!}\tau=(v_\tau^k)_x-\pi_\tau^k\,,
\\&\mu(\pi_\tau^k,\theta_\tau^{k-1})\xi_\tau^k=\bbC\eps_\tau^k+\eta(\pi_\tau^k)_{xx}
\ \ \ \text{ with }\ \xi_\tau^k\in{\rm Sign}(\pi_\tau^k)\,,
\\&\frac{\theta_\tau^k-\theta_\tau^{k-1}\!\!}\tau=f_0(\theta_\tau^k)
-|\pi_\tau^k|f_1(\theta_\tau^k)+\kappa(\theta_\tau^k)_{xx}\,\label{disc4}
\end{align}\end{subequations}
to be solved for $k=1,2,...,T/\tau$ starting for $k=1$ from the
initial conditions $v_\tau^0=v_0$, $\eps_\tau^0=\eps_0$, and $\theta_\tau^0=\theta_0$.
The boundary conditions for \eqref{disc} are like in \eqref{BC} but now with
time-varying velocity $v_\infty$, i.e.\
\begin{align}\label{BC-disc}
v_\tau^k(\pm H)=\pm v_\infty^k=:\int_{(k-1)\tau}^{k\tau}\!\!\!\frac{v_\infty(t)}\tau \dd t,
\ \ \ \ \ \ \ \pi_\tau^k(\pm H)=0,\ \ \ \ \ \ \ \theta_\tau^k(\pm H)=\theta_\infty\,.
\end{align}

The system (\ref{disc}a-c) has a variational structure with a convex coercive
potential
\begin{align}
(v,\eps,\pi)\mapsto\int_\Domain\!
\varrho\frac{(v{-}v_\tau^{k-1})^2\!}{2\tau}
+\bbC\eps(v_x{-}\pi)+\bbC\,\frac{(\eps{-}\eps_\tau^{k-1})^2\!}{2\tau}
+R(\pi,\theta_\tau^{k-1})+\frac\eta2\pi_x^2 \dd x\,.
\end{align}
For a sufficiently small $\tau>0$, this potential is convex and coercive
on $L^2(\Domain)^2\times H^1(\Domain)$. Minimization of this functional on an
affine manifold respecting the boundary conditions $v(\pm H)=\pm v_\infty^k$,
$\pi(\pm H)=0$, and $\theta(\pm H)=\theta_\infty$ gives  by 
the standard direct-method argument existence of an (even unique) minimizer,
let us denote
it by $(v_\tau^k,\eps_\tau^k,\pi_\tau^k)\in L^2(\Domain)^2\times H^1(\Domain)$.
This minimizer satisfies (\ref{disc}a,b) in the weak sense
and also the inclusion $\pl_\pi R(\pi_\tau^k,\theta_\tau^{k-1})\ni
\bbC\eps_\tau^k+\eta(\pi_\tau^k)_{xx}$. Therefore, there exists
$\xi_\tau^k\in{\rm Sign}(\pi_\tau^k)\subset H^1(\Domain)^*$ such
that $\mu(\pi_\tau^k,\theta_\tau^{k-1})\xi_\tau^k=\bbC\eps_\tau^k+\eta(\pi_\tau^k)_{xx}$
in the weak sense. Then we can solve \eqref{disc4} by
minimization of the convex functional
\begin{align}
\theta\mapsto\int_\Domain\!\frac{(\theta-\theta_\tau^{k-1})\!}{2\tau}
+|\pi_\tau^k|\varphi_1(\theta)-\varphi_0(\theta)+\frac\kappa2\theta_x^2 \dd x\,,
\end{align}
where $\varphi_i$ are the primitive functions to $f_i$, $i=0,1$. This functional 
is coercive on a linear manifold of the space $H^1(\Domain)$ respecting
the boundary condition \eqref{BC}. Let us denote its unique minimizer by
$\theta_\tau^k$.

We introduce the piecewise affine continuous and the piecewise constant
interpolants. Having $\{v_\tau^k\}_{k=0}^{T/\tau}$, we define
\begin{align}\label{def-of-interpolants}
&\wb v_\tau(t):=v_\tau^k,\ \ \ \underline v_\tau(t):=v_\tau^{k-1},
\ \text{ and }\ v_\tau(t):=\Big(\frac t\tau{-}k{+}1\Big)v_\tau^k
\!+\Big(k{-}\frac t\tau\Big)v_\tau^{k-1}
\end{align}
for $(k{-}1)\tau<t\le k\tau$ with $k=0,1,...,T/\tau$. Analogously,
we define also $\wb\eps_\tau$, or $\underline\theta_\tau$, etc.
This allows us to write the system \eqref{disc} in a ``compact'' form:
\begin{subequations}\label{disc+}\begin{align}\label{disc1+}
&\varrho\DT v_\tau -(\bbC\wb\eps_\tau)_x=0\,,
\\&\DT\eps_\tau=(\wb v_\tau)_x-\wb\pi_\tau\,,\label{disc2+}
\\&\mu(\wb\pi_\tau,\underline\theta_\tau)\wb\xi_\tau=
\bbC\wb\eps_\tau+\eta(\wb\pi_\tau)_{xx}
\ \ \ \text{ with }\ \wb\xi_\tau\in{\rm Sign}(\wb\pi_\tau)\,,\label{disc3+}
\\&\DT\theta_\tau=f_0(\wb\theta_\tau)
-|\wb\pi_\tau|f_1(\wb\theta_\tau)+\kappa(\wb\theta_\tau)_{xx}\,.\label{disc4+}
\end{align}\end{subequations}

By modifying appropriately the procedure which led to the a priori estimates
(\ref{est}a-c,e), we obtain here
\begin{subequations}\label{est++}\begin{align}
&\|\wb v_\tau\|_{L^\infty(I;L^2(\Domain))}^{}\le C\,,
\\&\|\wb\eps_\tau\|_{L^\infty(I;L^2(\Domain))}^{}\le C\,,
\\&\|\wb\pi_\tau\|_{L^2(I;H^1(\Domain))}^{}\le C\,,
\\&\label{est++theta}
\|\wb\theta_\tau\|_{L^\infty(I{\times}\Domain)\cap L^2(I;H^1(\Domain))}^{}\le C\,,\ \
\text{ and here also}
\\&\|\wb\xi_\tau\|_{L^\infty(I{\times}\Domain)\cap L^2(I;H^1(\Domain)^*)}^{}\le C\,.
\end{align}\end{subequations}
All these estimates hold also for the piecewise affine interpolants, and
\eqref{est++theta} holds also for $\underline\theta_\tau$. The last estimate is
obtained by comparison from
$\wb\xi_\tau=
(\bbC\wb\eps_\tau{+}\eta(\wb\pi_\tau)_{xx})/\mu(\wb\pi_\tau,\underline\theta_\tau)$
when testing it by functions bounded in $L^2(I;H^1(\Domain))$ and using the
smoothness of $1/\mu(\wb\pi_\tau,\underline\theta_\tau)$.

Then, by the Banach selection principle, we obtain subsequences indexed, for
simplicity, again by $\tau$) weakly* converging in the topologies indicated in
\eqref{est++}, and we pass to a limit for $\tau\to0$ and are to show
that such limit (let us denote it by $(v,\eps,\pi,\theta,\xi)$) solve the
continuous problem with $\pi=\DT\pp$. For this, one uses the Aubin-Lions
compactness theorem adapted for the time-discretization method as in
\cite[Sect.\,8.2]{Roub13NPDE}. Thus we can rely on that
\begin{align}\label{theta-strong}
\wb\theta_\tau\to\theta\quad\text{ strongly in }\ L^c(I{\times}\Domain)\
\text{ for any $1\le c<\infty$}.
\end{align}
The limit passage in the linear hyperbolic equation \eqref{evol1} is due to a
weak convergence of both $v$ and $\eps$ and also the limit passage in the
linear equation \eqref{evol2} is easy via weak convergence. Yet, there is one
peculiarity in the limit passage in the nonlinearity in \eqref{evol3} for which
a strong convergence of $\eps$ is needed, but we do not have any information
about space gradient of $\eps$. The other peculiarity is a need of the strong
convergence of $\DT\pp$ which is needed for \eqref{evol5}, but we do
not have any information about $\DT\pi$, so that mere compactness arguments
cannot be used. This can be obtained from the momentum equation \eqref{evol1}
and from \eqref{evol3} when using the strong monotonicity of the operators in
\eqref{evol1} and \eqref{evol3} simultaneously. As for \eqref{evol3}, note that
$\mu(\pi,\theta){\rm Sign}(\pi)=\pl_\pi R(\pi,\theta)$ and that
$R(\cdot,\theta)$ is convex, to that $\pl_\pi R(\cdot,\theta)$ is monotone. In
particular, for any $\wb\xi_\tau\in{\rm Sign}(\wb\pi_\tau)$ and
$\xi\in{\rm Sign}(\pi)$, we have
$\int_0^t\langle\wb\xi_\tau{-}\xi,\wb\pi_\tau{-}\pi\rangle \dd t\ge0$, where
$\langle\cdot,\cdot\rangle$ denotes the duality pairing between
$H^1(\Domain)^*$ and $H^1(\Domain)$.

The usage of this monotonicity of the set-valued mapping
$\pl_\pi R(\cdot,\theta)$ should be done carefully. The time-discrete approximation
of \eqref{eq3-weak+++} gives some $\wb\pi_\tau\in L^2(I;H^1(\Domain))$
and $\wb\xi_\tau\in L^2(I;H^1(\Domain)^*)$ satisfying \eqref{disc3+}
together with the boundary conditions $\pp(\pm H)=0$ in the weak form.
From the mentioned monotonicity and by using \eqref{disc1+} and \eqref{disc3+} tested
by $\wb v_\tau{-}v$ and $\wb\pi_\tau{-}\pi$ and integrated over a time interval
$[0,t]$ and the domain $\Domain$, we obtain
\begin{align}\nonumber
&\int_\Domain\frac\varrho2(v_\tau(t){-}v(t))^2
+\frac12\bbC(\eps_\tau(t){-}\eps(t))^2 \dd x+ \int_0^t\!\!\int_\Domain
\eta(\wb\pi_\tau{-}\pi)_x^2 \dd x\d t
\\[-.3em]&\nonumber\le
\int_0^t\!\bigg(\big\langle\varrho\DT v_\tau{-}\varrho\DT v,
\wb v_\tau{-}v\big\rangle+
\big\langle\DT\eps_\tau{-}\DT\eps,\bbC\wb\eps_\tau{-}\bbC\eps\big\rangle
+\big\langle\mu(\wb\pi_\tau,\underline\theta_\tau)\wb\xi_\tau
-\mu(\pi,\underline\theta_\tau)\xi,\wb\pi_\tau{-}\pi\big\rangle
\\[-.3em]&\hspace{5em}\nonumber
+\big\langle\varrho\DT v,\wb v_\tau{-}v_\tau\big\rangle
+\big\langle\DT\eps,\bbC\wb\eps_\tau{-}\bbC\eps_\tau\big\rangle
+\!\int_\Domain\eta(\wb\pi_\tau{-}\pi)_x^2 \dd x\bigg) \dd t
\\[-.3em]&\nonumber=-\int_0^t\!\bigg(\big\langle\varrho\DT v,
\wb v_\tau{-}v\big\rangle+\big\langle\DT\eps_\tau{-}\DT\eps,\bbC\eps\big\rangle
+\big\langle\mu(\pi,\underline\theta_\tau)\xi,\wb\pi_\tau{-}\pi\big\rangle
\\[-.5em]&\hspace{5em}
-\big\langle\varrho\DT v,\wb v_\tau{-}v_\tau\big\rangle
-\big\langle\DT\eps,\bbC\wb\eps_\tau{-}\bbC\eps_\tau\big\rangle
+\!\int_\Domain\!\eta\pi_x(\wb\pi_\tau{-}\pi)_x \dd x\bigg) \dd t\to0\,,
\label{->0}\end{align}
where $\langle\cdot,\cdot\rangle$ again denotes the duality pairing between
$H^1(\Domain)^*$ and $H^1(\Domain)$. The meaning of
$\langle\mu(\wb\pi_\tau,\underline\theta_\tau)\wb\xi_\tau
,\wb\pi_\tau{-}\pi\rangle$ for $\wb\xi_\tau$ valued in $H^1(\Domain)^*$ is
rather $\langle\wb\xi_\tau,\mu(\wb\pi_\tau,
\underline\theta_\tau)(\wb\pi_\tau{-}\pi)\rangle$, relying that
$\mu(\wb\pi_\tau,\underline\theta_\tau)(\wb\pi_\tau{-}\pi)$ is valued in
$H^1(\Domain)$; here we need $\mu$ smooth so that
$(\mu(\wb\pi_\tau,\underline\theta_\tau)(\wb\pi_\tau{-}\pi))_x
=\mu(\wb\pi_\tau,\underline\theta_\tau)(\wb\pi_\tau{-}\pi)_x
+(\mu_\pi'(\wb\pi_\tau,\underline\theta_\tau)(\wb\pi_\tau)_x+
\mu_\theta'(\wb\pi_\tau,\underline\theta_\tau)(\underline\theta_\tau)_x)
(\wb\pi_\tau{-}\pi)$ is valued in $L^2(\Domain)$. Similarly, it applies also
for $\langle\mu(\pi,\underline\theta_\tau)\xi,\wb\pi_\tau{-}\pi\rangle$. For
the inequality in \eqref{->0} see \cite[Remark~8.11]{Roub13NPDE}. For the
equality in \eqref{->0}, we used \eqref{disc2+} together with its limit
obtained by the weak convergence, i.e.\ $\DT\eps=v_x-\pi$, and also
(\ref{disc+}a,c) for the identity
\begin{align}\nonumber
&\big\langle\DT\eps_\tau{-}\DT\eps,\bbC\wb\eps_\tau{-}\bbC\eps\big\rangle
=\big\langle(\wb v_\tau{-}v)_x,\bbC\wb\eps_\tau\big\rangle
-\big\langle\wb\pi_\tau{-}\pi,\bbC\wb\eps_\tau\big\rangle
-\big\langle\DT\eps_\tau{-}\DT\eps,\bbC\eps\big\rangle
\\&\nonumber\ \ \ =-\big\langle\varrho\DT v_\tau,\wb v_\tau{-}v\big\rangle
-\big\langle\mu(\wb\pi_\tau,\underline\theta_\tau)\wb\xi_\tau,\wb\pi_\tau{-}\pi\big\rangle-\!\int_\Domain\!\eta(\wb\pi_\tau)_x(\wb\pi_\tau{-}\pi)_x\d x
-\big\langle\DT\eps_\tau{-}\DT\eps,\bbC\eps\big\rangle\,.
\end{align}
It is important, that \eqref{->0} holds for any $\xi\in{\rm Sign}(\pi)$ and, at
this moment, we do not assume that $\xi$ comes as a limit from the (sub)sequence
$\{\wb\xi_\tau\}_{\tau>0}$. 

To the convergence in \eqref{->0}, we used that
$\DT v\in L^2(I;H^1(\Domain)^*)$ while $\wb v_\tau{-}v\to0$ weakly 
$L^2(I;H^1(\Domain))$, and that
$\DT\eps_\tau{-}\DT\eps\to0$ weakly in $L^2(I;H^1(\Domain)^*)$, and eventually that
$\mu(\pi,\underline\theta_\tau)$ converges (to a limit which is not important here)
strongly in $L^c(I{\times}\Domain)$ due to \eqref{theta-strong} while
$\wb\pi_\tau{-}\pi\to0$ weakly in $L^2(I;H^1(\Domain))$ so that also
$\mu(\pi,\underline\theta_\tau)(\wb\pi_\tau{-}\pi)\to0$ weakly in
$L^2(I;H^1(\Domain))$. Therefore, considering \eqref{->0} integrated over $I$,
we obtain
\begin{subequations}\label{conv-v-eps-pi}\begin{align}
&&&\wb v_\tau\to v&&\text{strongly in }\ L^2(I{\times}\Domain)\,,&&&&\\
&&&\wb\eps_\tau\to\eps&&\text{strongly in }\ L^2(I{\times}\Domain)\,,\\
&&&\wb\pi_\tau\to\pi&&\text{strongly in }\ L^2(I;H^1(\Domain))\,.\label{conv-pi}
\end{align}\end{subequations}
In fact, by interpolation, (\ref{conv-v-eps-pi}a,b) holds even in
$L^c(I;L^2(\Domain))$ for any $1\le c<\infty$.
For \eqref{conv-pi}, we used the strong convergence of gradients of $\DT\pp_k$
and the fixed boundary conditions, so that we do not need to rely on
the monotonicity of $\pl_\pi R(\cdot,\theta)$ which may not be strong.

Having the strong convergence \eqref{conv-v-eps-pi} at disposal, the limit
passage is then easy, showing that the previously obtained weak limit
$(v,\eps,\pi,\theta)$ is a weak solution to the system \eqref{evol}.
In particular, from the inclusion in \eqref{disc3+} one obtains
$\xi\in{\rm Sign}(\pi)$ by using maximal monotonicity of the graph of
the set-valued mapping
${\rm Sign}:L^2(I;H^1(\Domain))\rightrightarrows L^2(I;H^1(\Domain)^*)$
and the strong convergence \eqref{conv-pi}. Thus (i) is proved.

As to (ii), if $\mu(\pi,\theta)\le \mathscr{O}(|\pi|^s)$, then
$\eta\pi_{xx}\in\bbC\eps-\mu(\pi,\theta){\rm Sign}(\pi)$ is bounded
in $L^s(I;L^2(\Domain)$ so that $\pi\in L^s(I;H^2(\Domain))$.

If $s\ge2$,
the procedure which led to the energy balance \eqref{engr} considered here without
$\alpha$-terms but integrated over a time interval $[0,t]$ was indeed rigorous.
This is because $v\in L^2(I;H^1(\Domain))$, as can be seen by comparison from
\eqref{evol2}, is in duality with $\varrho\DT v\in L^2(I;H^1(\Domain)^*)$ and with
$(\bbC\eps)_x\in L^2(I;H^1(\Domain)^*)$, so that testing the momentum equation
\eqref{evol1} and the related by-part integration is legitimate. Similar arguments
concern also the aging rule \eqref{evol5}.
Since $\eta\pi_{xx}\in L^2(I{\times}\Domain)$ if $s\ge2$, also the
test of the plastic rate equation \eqref{evol3} by $\pi\in L^s(I{\times}\Domain)$
is legitimate together with the related by-part integrations.

In this case when $s\ge2$, also \eqref{disc4+} can be tested by $\DT\theta_\tau$,
which gives the regularity
$\theta\in H^1(I;L^2(\Domain))\cap L^\infty(I;H^1(\Domain))$.
By comparison $\kappa\theta_{xx}=\DT\theta+|\pi|f_1(\theta)-f_0(\theta)\in
L^2(I{\times}\Domain)$, we obtain also $\theta\in L^2(I;H^2(\Domain))$.
\end{proof}

\begin{remark}[{\sl Stability and time-periodic solutions}] In geodynamics the
  phenomenon called \emph{episodic tremor and slip} describes time-periodic
  motions in subduction zones where shorter periods of plastic slips alternate
  with longer periods with slow slip events. Hence, it would be interesting to
  complement our existence result for ``transient events'' governed by the
  above initial-value problem by a theory for time-periodic solutions. The
  aim would be show that there is a period $t_*>0$ and a solution of the
  system \eqref{evol} with the boundary conditions \eqref{BC} satisfying
  $(\DT v, \DT \eps, \DT \alpha, \DT\theta) \not\equiv 0$ and 
\begin{align}\label{PC}
v(\cdot,t_*)=v(\cdot,0)\,,\ \ \ \eps(\cdot,t_*)=\eps(\cdot,0)\,,\ \ \
\alpha(\cdot,t_*)=\alpha(\cdot,0)\,,\ \text{ and }\ \ 
\theta(\cdot,t_*)=\theta(\cdot,0)
\end{align}
instead of \eqref{IC}. 
Of course, a general question is that of stability of the steady state
solutions $(\pi,\theta)$ obtained in Section \ref{se:AnaSteady} or potentially
of such time-periodic solutions as described here. As we will see in the
following section, one indication of the existence of time-periodic solutions
is the loss of stability of the steady state solution. But because of the
complexity of the model, these questions are beyond the scope of this
paper.
\end{remark}

\begin{remark}[{\sl Asymptotics for $\eta\to0$ and $\kappa\to0$}]
  Unlike to the case for steady solutions for \eqref{eq} as in
  Section~\ref{sec-localization}, it is not possible in the evolutionary model
  \eqref{evol} to pass to the limit for $\eta\to0$. In particular, a limit
  passage in the term $\bbC\eps^\eta(\wt\pi^\eta{-}\DT\pp^\eta)$ occurring in
  \eqref{eq3-weak+++} seems to be out of reach. The substitution
  \eqref{subst} by a convex term in $\DT\pp$ could not help, being not weakly
  upper-semicontinuous. If also \eqref{ass-M1} holds, then  like in
  Propositions \ref{prop2} and \ref{prop3}, we can at least obtain some
  uniform bounds, in particular for the plastic strain rate $\pi=\DT\pp$ in the
  Orlicz space $L_\varPhi(I{\times}\Domain)$ with $\varPhi$ from
  \eqref{ass-M1}, i.e.\ $\int_I\int_D \varPhi(\pi(t,x))\dd x \dd t
  <\infty$. Yet, the limit passage for $\eta\to0$, even while keeping
  $\kappa>0$ fixed, remains intractable.
\end{remark}

\section{Illustrative numerical simulations}
\label{se:NumSimul}

We illustrate the response of the evolutionary model in Section~\ref{sec-evol}
by a simplified model derived in Section
\ref{su:SimplifMod}. This model still has exactly the same steady states
as the full model, such that all the theory of Section \ref{se:AnaSteady} applies
to it, when ignoring statements about the damage variable $\alpha$. We expect
that the simplified model is still relevant as far as
usually observed dynamical features concern. Moreover,  it also displays the 
effect of the free boundary occurring between the elastic zone and the plastic
zone. In Section \ref{su:NumSimSteady} we show by numerical simulations
that the steady states localize for $v_\infty\to 0$ in such a way
that $\pi_\stst$ has support (i.e.\ the so-called cataclastic zone)
in $[-h_*(v_\infty,\kappa),h_*(v_\infty,\kappa)]$ with
$h_*(v_\infty,\kappa)\sim\sqrt\kappa$ for $\kappa \to 0^+$. Moreover, we show
that, when keeping $v_\infty\neq 0$ fixed but sufficiently small, we obtain a support
with $h_*(v_\infty,\kappa) \to v_\infty/\pi_*$ for $\kappa \to 0^+$.  

In Section \ref{su:NumSimODE} \EEE we study an ODE model for
scalars $\theta(t)$ and
$\sigma(t)$ which displays the effect of oscillatory behavior for $|v_\infty| <
v_\text{crit} $ while solutions converge to the unique steady state for
$|v_\infty| > v_\text{crit} $. Finally Section \ref{su:NumSimPDE} 
presents
simulations for the simplified evolutionary model. In particular, we observe
again that for small nontrivial values of $|v_\infty|$ we have oscillatory
behavior, where the plastic zone is spatially and temporarily localized in the
sense that the support of $\pi(t,\cdot)$ is compactly contained in $D=[-H,H]$
for all $t \in [0,T_\text{per}]$ and that $\pi(t,x)=0$ for all $x\in D$ and
all $t\in [t_1,t_2]$ for a nontrivial interval $[t_1,t_2]\subset [0,T_\text{per}]$.  
For $|v_\infty| $ large, we find convergence into a steady state with a
nontrivial plastic (cataclastic) zone. 
All the following results are derived from numerical experiments only.

\subsection{The simplified model without damage}
\label{su:SimplifMod}

To display the main features of our rate-and-state friction model we reduce
the full evolutionary model \eqref{evol} by making the following simplifications: 
\\
\textbullet\ we neglect inertial effects (i.e.\ we set $\varrho=0$ in
\eqref{evol1}), thus\\
\mbox{}\quad making the system quasistatic but still keeping a
rate-and-state dependent plasticity;
\\
\textbullet\ we choose $\eta=0$ for the length-scale parameter in \eqref{evol3}\\
\mbox{}\quad 
as analyzed in Section~\ref{sec-localization} for the steady-state solutions;  
\\
\textbullet\ we neglect all damage effects through $\alpha$ and
omit \eqref{evol4} as we did in Theorem~\ref{th:EvolExist}.

\noindent
Because of $\varrho=0$, the momentum balance leads to a spatially constant
stress $\sigma(t)=\bbC \eps$. As now $\bbC$ is constant, also $\eps(t)$ is
spatially constant. Integrating \eqref{evol2} over $x\in D =[{-}H,H]$ and using
the boundary condition for $v$ from \eqref{BC-evol} gives the following coupled system
for $\sigma$, $\pi=\DT p$, and $\theta$:
\begin{subequations}
\label{eq:SimpMod}
\begin{align}
\label{eq:SimpMod.a}
&\frac{2H}{\bbC} \DT \sigma + \int_\Domain \pi \dd x = 2 v_\infty(t),
\\
\label{eq:SimpMod.b}
&\mu(\pi,\theta) \mathrm{Sign}(\pi) \ni \sigma,
\\
& \DT\theta = f_0(\theta)-|\pi|f_1(\theta) + \kappa \theta_{xx} , \quad
\theta(t,\pm H)=\theta_\infty.
\label{eq:SimpMod.c}
\end{align}
\end{subequations}

Throughout this section we assume that $\mu$ has the form
\[
\mu(\pi,\theta) = \mu_0 + A(\pi) + B(\theta)\ \ \text{ with }\ \ 
A(\pi),B(\theta)\geq 0\ \ \text{ and }\ \ A(-\pi)=A(\pi)\,;
\]
cf.\ also \eqref{DR1++}.  Assuming further $A'(\pi)>0$ for $\pi>0$ we can solve
\eqref{eq:SimpMod.b} in the form
\begin{equation}
  \label{eq:def.Pi}
  \pi=\varPi(\sigma, \theta)\ \ \text{ with }\ \ \varPi(\sigma,\theta) = \left\{
\begin{array}{cl} 
0 & \text{for }|\sigma| \leq  \mu_0{+}B(\theta), \\
A^{-1}\big(\sigma{-}\mu_0{-}B(\theta)\big)& \text{for }\sigma>\mu_0{+}B(\theta), \\
\!\!-A^{-1}\big(|\sigma{-}\mu_0{-}B(\theta)|\big)& \text{for }\sigma<-\mu_0{-}B(\theta).\\
\end{array} \right. 
\end{equation}
Thus, we obtain our final coupled system of a scalar ODE for $\sigma$  with a
non-locally coupled scalar parabolic PDE for $\theta$, namely 
\begin{subequations}
\label{eq:SM}
\begin{align}
& \label{eq:SM.a}
\DT \sigma = \frac{\bbC}H\, v_\infty(t) - \frac{\bbC}{2H}
 \int_\Domain \varPi(\sigma,\theta) \dd x, 
\\
& \label{eq:SM.b} 
\DT\theta = f_0(\theta)-|\varPi(\sigma,\theta)|f_1(\theta) + \kappa \theta_{xx} , \quad
\theta(t,\pm H)=\theta_\infty.
\end{align}
\end{subequations}
Here the nonsmoothness due to the plastic behavior is realized by the nonsmooth
function $\pi=\varPi(\sigma,\theta)$ defined in \eqref{eq:def.Pi}. 

For all the following simulation we choose the 
following parameters and functions:
\begin{equation}
  \label{eq:DataSimple}
  \begin{aligned}
    &H=1,\ \ \bbC=1,\ \ \theta_\infty=10,\ \ \mu_0=1,\ \
    f_0(\theta)=1-\theta/\theta_\infty,
    \\
    &f_1(\theta)=10\, \theta,\ \ A(\pi)=\ln(|\pi|{+}1),\ \
    B(\theta)=\ln(4\theta{+}1).
  \end{aligned}
\end{equation}
Subsequently, we will only vary the coefficient $\kappa>0$ and the shear velocity
$v_\infty$. 

\subsection{Steady states}
\label{su:NumSimSteady}

\begin{figure}[h]
\centerline{\scriptsize\sf Stationary profiles $\theta_\stst $ of the aging variable}
\begin{tabular}{cccc}
\includegraphics[width=0.23\textwidth]{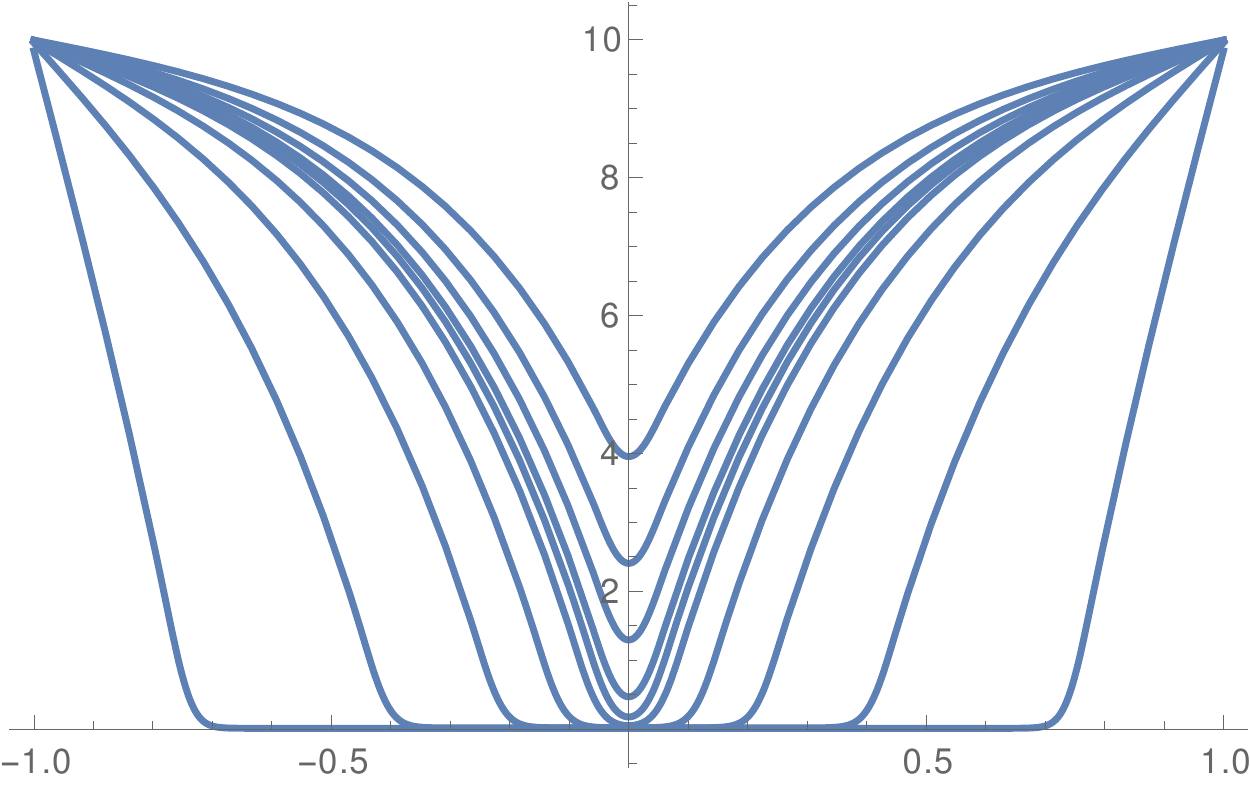}& 
\includegraphics[width=0.23\textwidth]{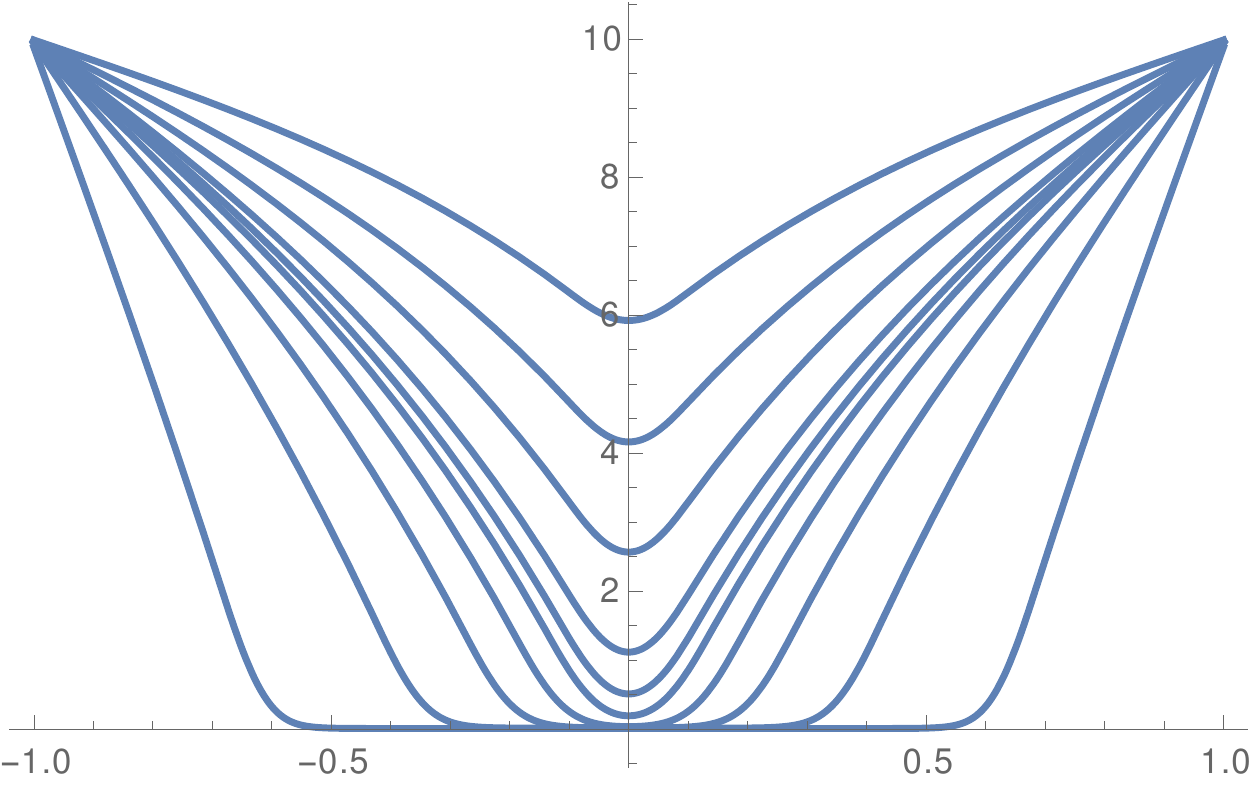}& 
\includegraphics[width=0.23\textwidth]{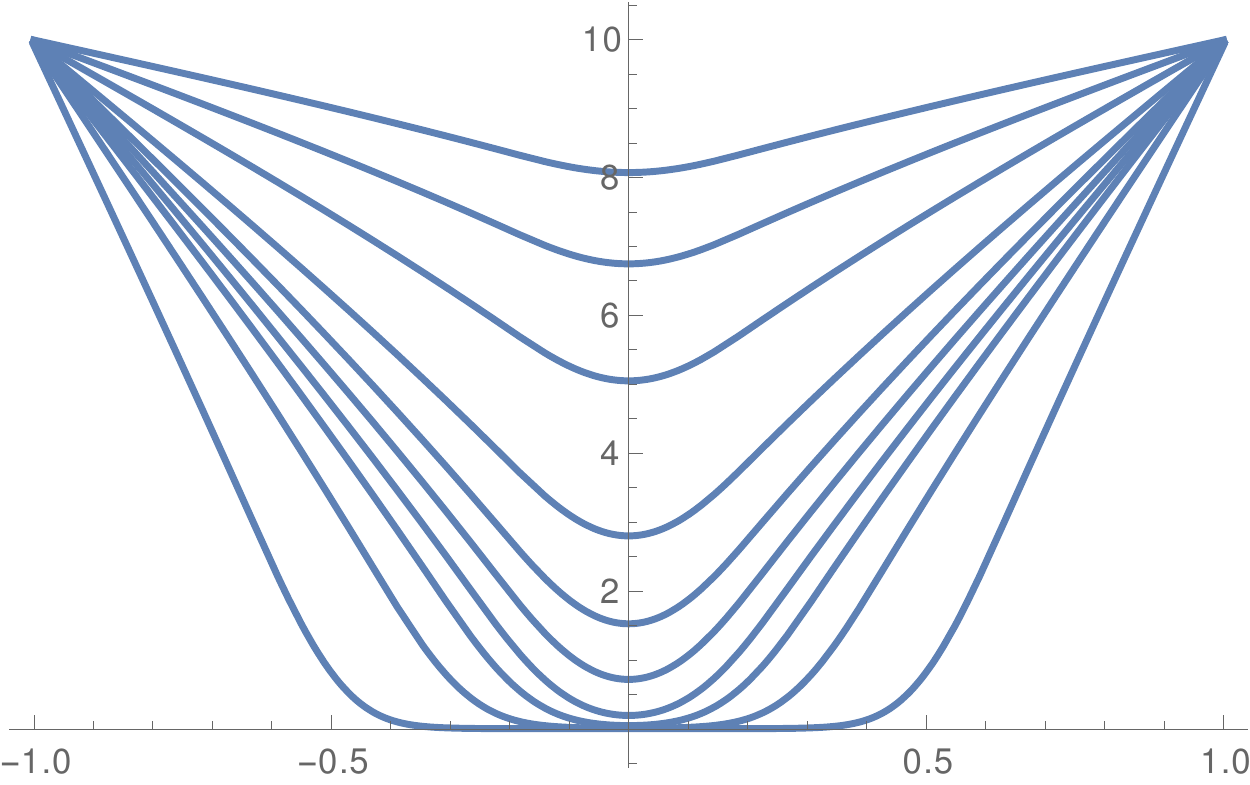}& 
\includegraphics[width=0.23\textwidth]{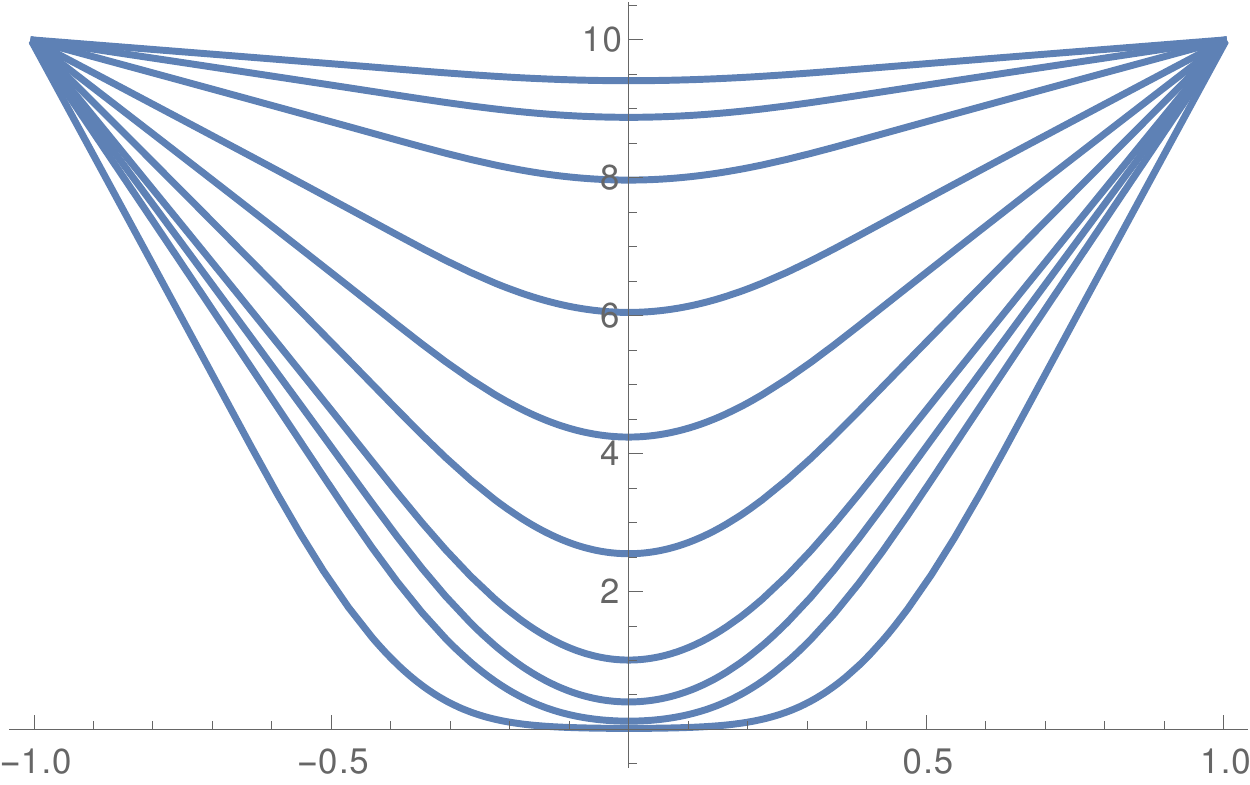}
\end{tabular}
\\[0.5em]
\centerline{\scriptsize\sf Stationary profiles $\pi_\stst $ of the plastic
  strain rate\hspace{9em}}\vspace*{-0.5em}
 \begin{tabular}{cccc}
 \includegraphics[width=0.23\textwidth,height=0.12\textwidth]{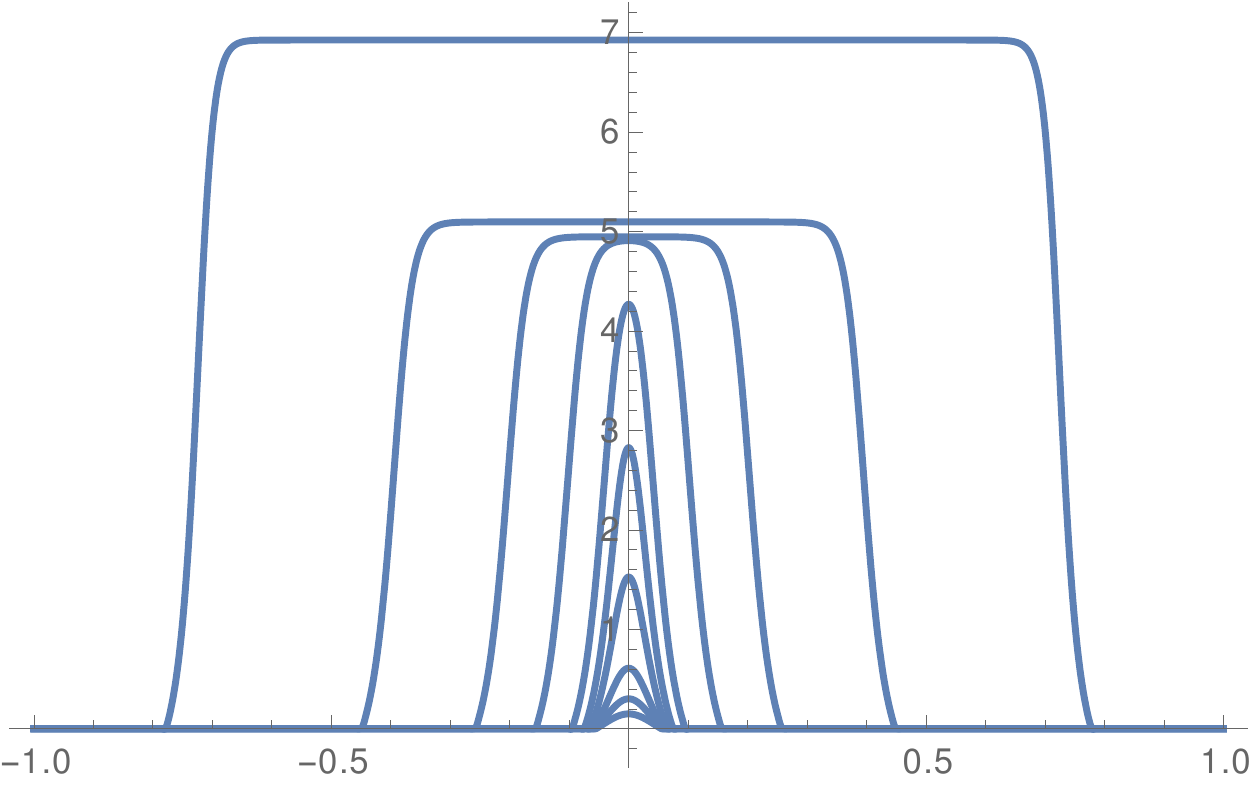}& 
\includegraphics[width=0.23\textwidth,height=0.146\textwidth]{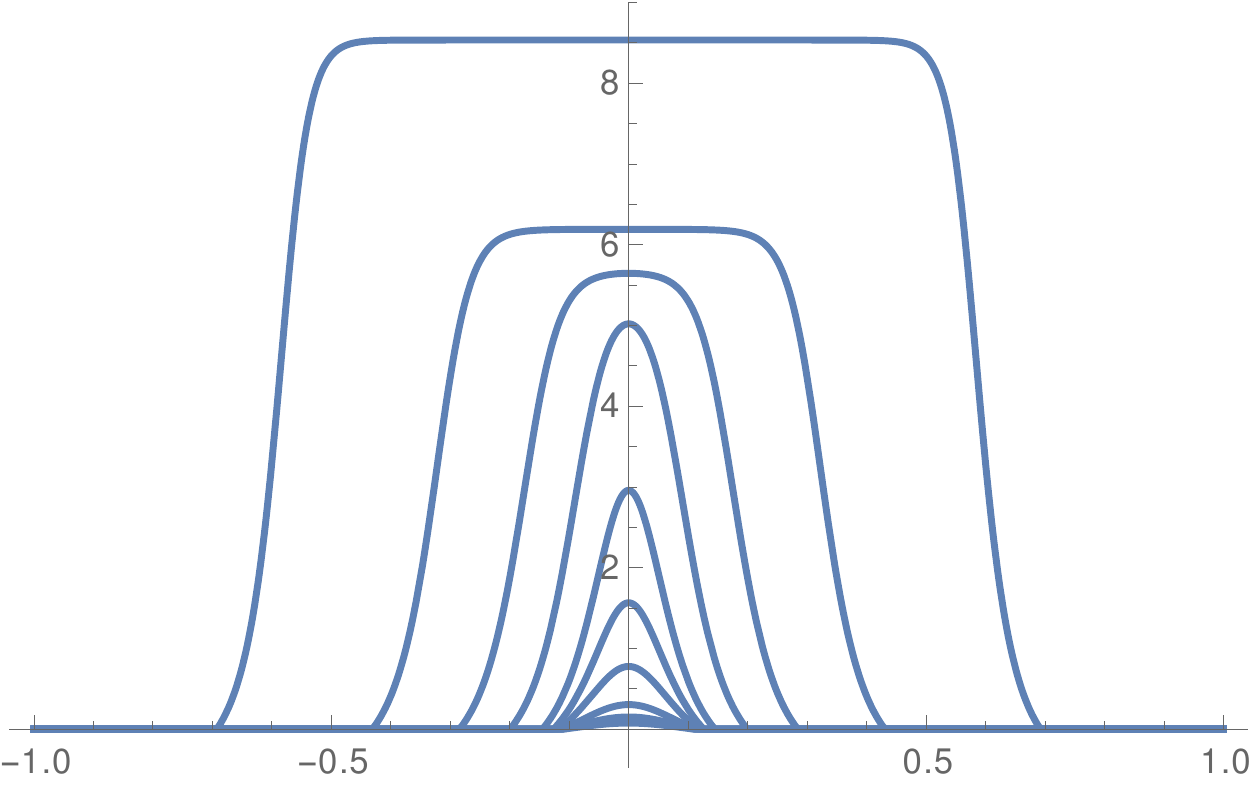}& 
\includegraphics[width=0.23\textwidth,height=0.196\textwidth]{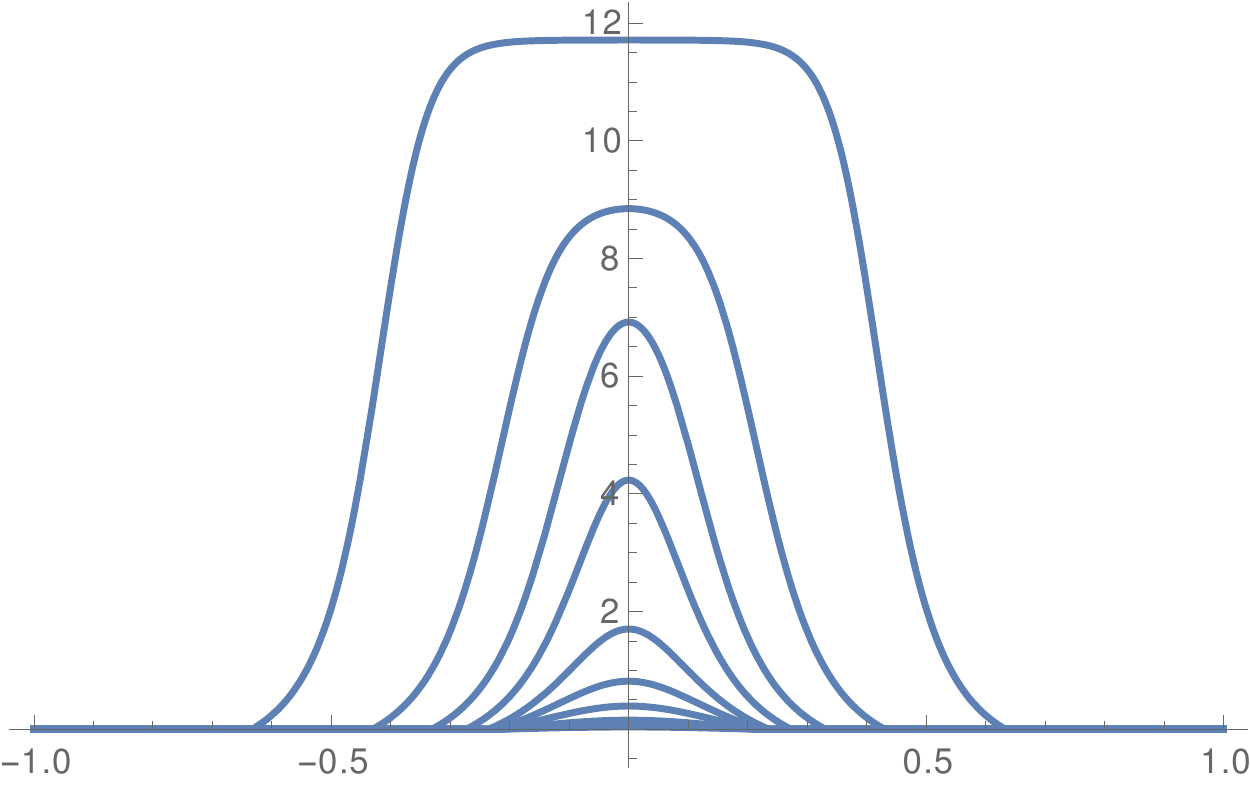}&
\includegraphics[width=0.23\textwidth,height=0.196\textwidth]{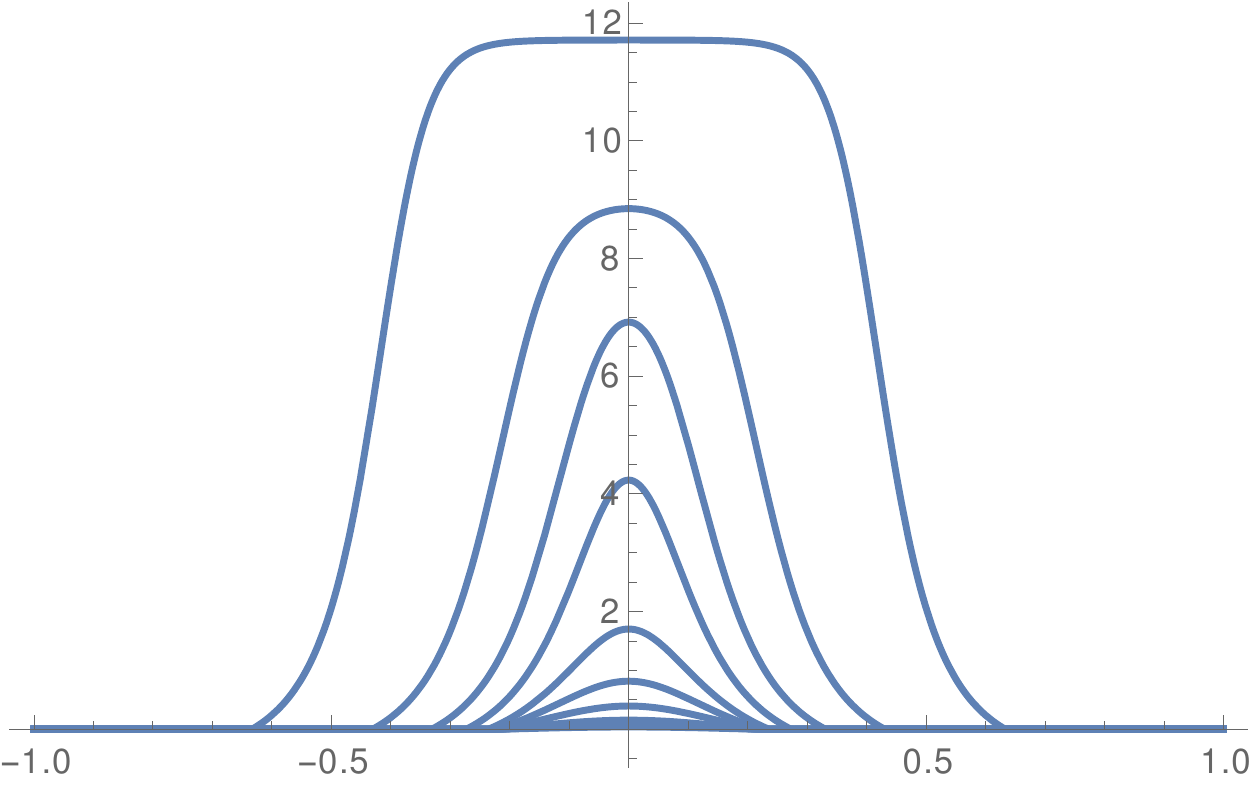}
\\
{\small$\kappa=0.01$}&{\small$\kappa=0.04$}&{\small$\kappa=0.16$}&{\small$\kappa=0.64$}
\end{tabular}
\caption{\small\sl Each picture shows ten curves that
correspond to the shear velocities $v_\infty\in \{0.005,\,0.01,\,0.02,\,0.05,\,0.1,\,0.2,\,0.5,\,1.0,\,2.0,\,5.0\}$, respectively. The rows shows
  $\theta_\stst $ (decreasing with $v_\infty$) and the lower rows shows
  $\pi_\stst $ (growing $v_\infty$).}
\label{fig:De-In-Crease}
\end{figure}

We first discuss the steady states for \eqref{eq:SM}, which are indeed a
special case of the steady states obtained in Proposition \ref{prop2}. 
Numerically, we always found exactly one steady state
$\theta_\stst = \Theta(v_\infty,\kappa)$, but were unable to prove its
uniqueness rigorously. When varying the parameters $v_\infty$ and $\kappa$ we
can easily observe clear trends for $(\theta_\stst ,\pi_\stst )$, where the
associated plastic flow rate is given by
$\pi_\stst =P(\sigma_\stst ,\theta_\stst )$, see Figure
\ref{fig:De-In-Crease}. We first observe that for fixed $\kappa$ the functions
$\theta_\stst $ and $\pi_\stst $ depend
monotonically  on $v_\infty$ in the
expected way, namely $\theta_\stst $ decreases with the shear velocity
$v_\infty$, while $\pi_\stst $ increases, which fits to the relation
$2v_\infty= \int_D \pi_\stst (v_\infty, \kappa; x) \dd x$.

Moreover, for $v_\infty\to 0^+$ the scaled plastic rate
$\pi_\stst/v_\infty$ converges to a
nontrivial limit with localized support, while $\theta_\stst $ converges
uniformly to $\theta_\infty$. For larger and larger $v_\infty$ the plastic zone
occupies more and more of the domain $D=[{-}1,1]$ and $\theta_\stst $ is 
very small in most of the plastic zone, namely $\theta \approx
\Theta_f(\pi)=\theta_\infty/(1{+}10 \pi \theta_\infty)\approx 1/(10 \pi)$.

\begin{figure}[h]
\centerline{\scriptsize\sf Rescaled stationary profiles $\pi_\stst/v_\infty $ of the
  plastic strain rate}
\begin{tabular}{cccc}
\includegraphics[width=0.23\textwidth,height=0.2\textwidth]{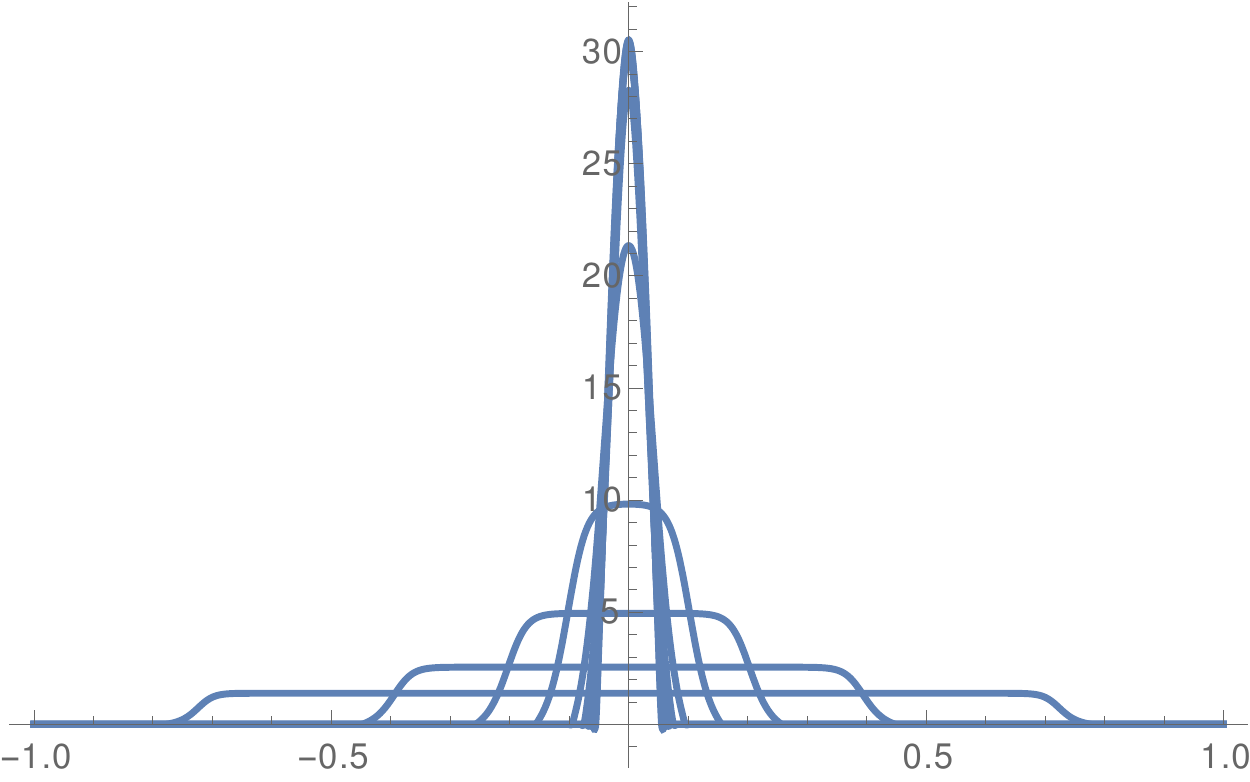}& 
\includegraphics[width=0.23\textwidth,height=0.2\textwidth]{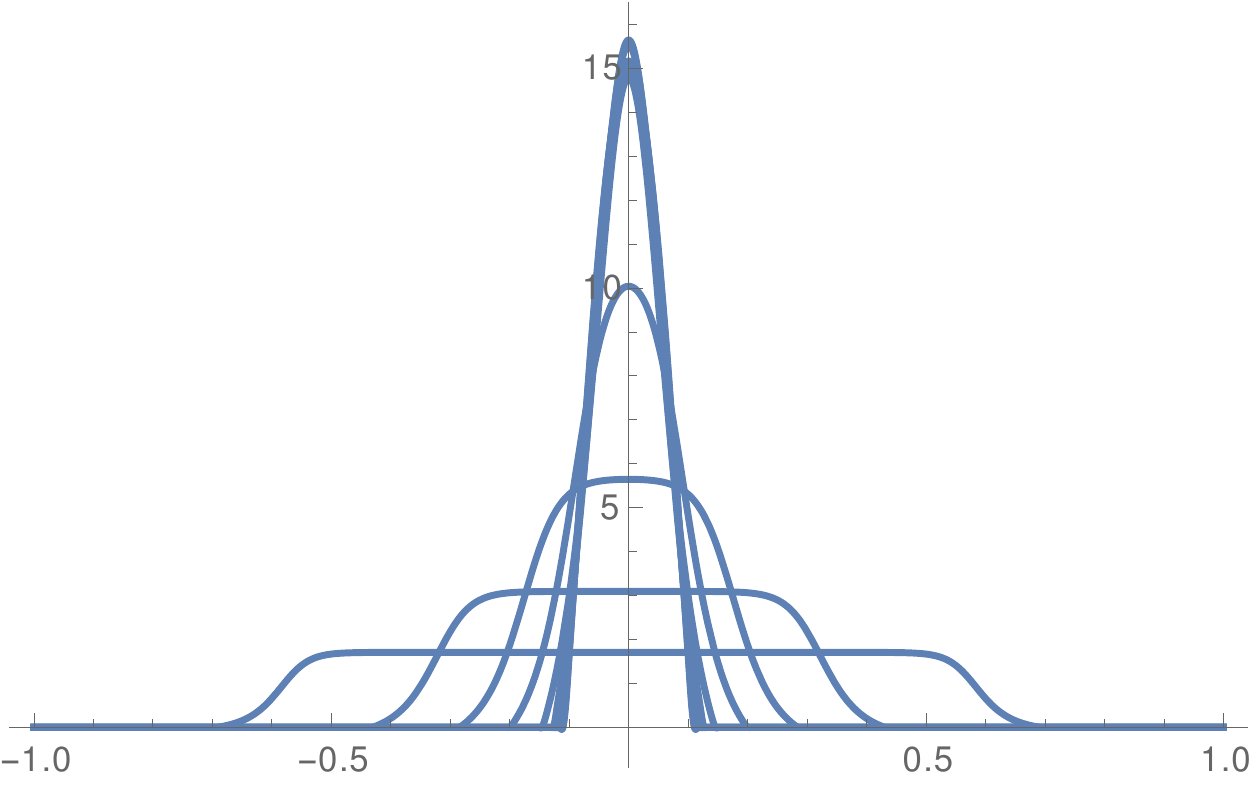}& 
\includegraphics[width=0.23\textwidth,height=0.2\textwidth]{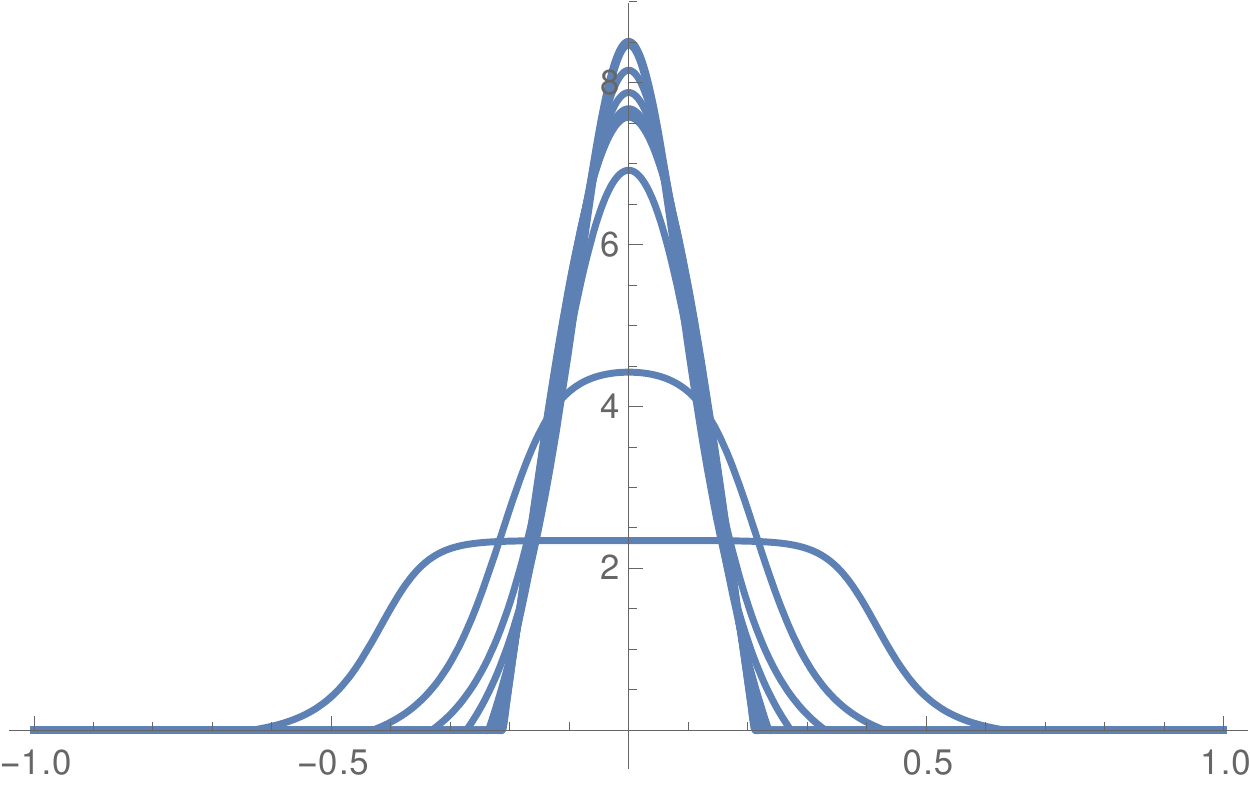}&
\includegraphics[width=0.23\textwidth,height=0.2\textwidth]{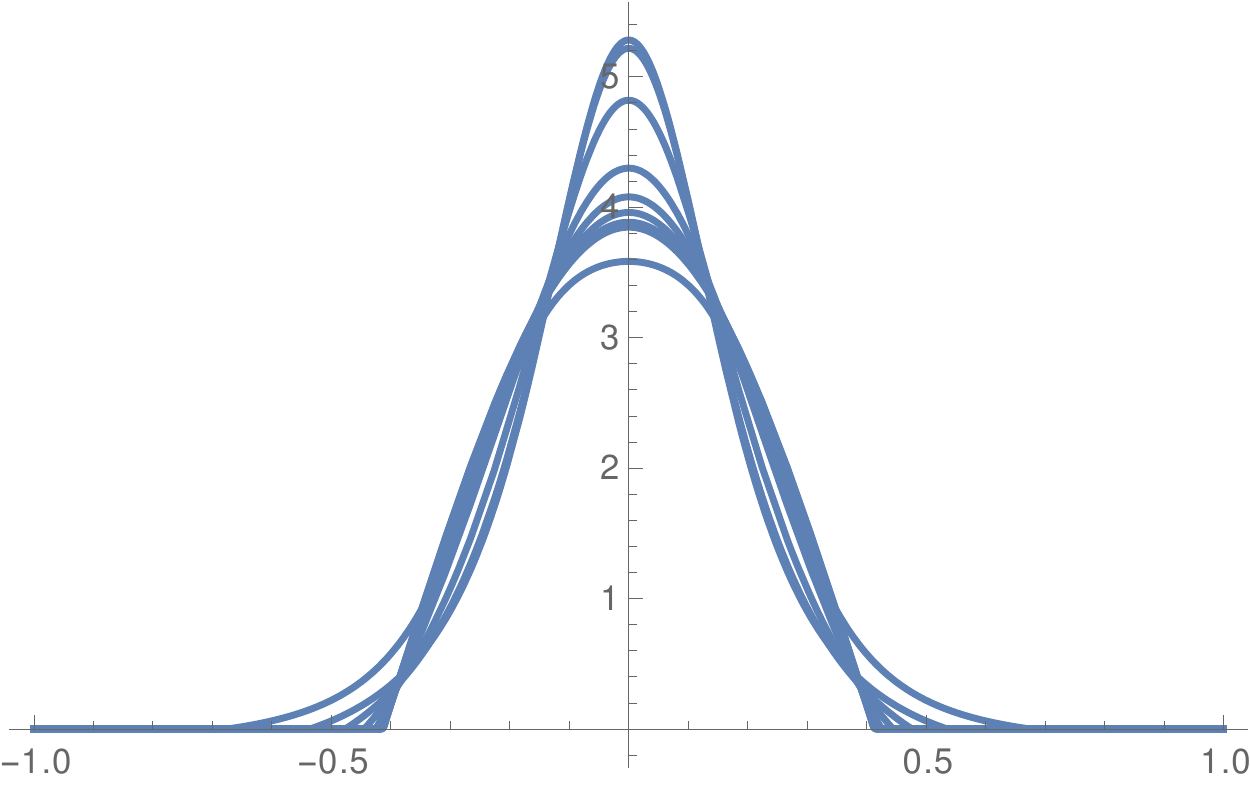}
\\
{\small$\kappa=0.01$}&{\small$\kappa=0.04$}&{\small$\kappa=0.16$}&{\small$\kappa=0.64$}
\end{tabular}
\caption{\small\sl The figures display the rescaled plastic strain rates
  $\pi_\stst/v_\infty$ for shear velocities $v_\infty\in
  \{0.005,\,0.01,\,0.02,\,0.05,\,0.1,\,0.2,\,0.5,\,1.0,\,2.0,\,5.0\}$,
  respectively. For $v_\infty\to 0$ one sees convergence to a limit shape with
  minimal support $[-h_*(\kappa), h_*(\kappa)]$ where
  $\kappa(0.01)\approx 0.055$, $\kappa(0.04)\approx 0.11$,
  $\kappa(0.16)\approx 0.21$, and $\kappa(0.64)\approx 0.41$. Effectively, we
  can see a free boundary between active cataclastic core zone and the rest of
  the fault.}
\label{fig:Support}
\end{figure}

When reducing the size of $\kappa$ we also see that the size of the plastic
zone shrinks. For small $v_\infty$ it can be seen that the support of
$\pi_\stst  $ is $[{-}h_*(v_\infty,\kappa),h_*(v_\infty,\kappa)]$ with
$h_*(v_\infty,\kappa) \sim \sqrt\kappa$, see Figure \ref{fig:Support}. 

\begin{figure}
\begin{tabular}{@{}ccc@{}}
\includegraphics[width=0.31\textwidth]{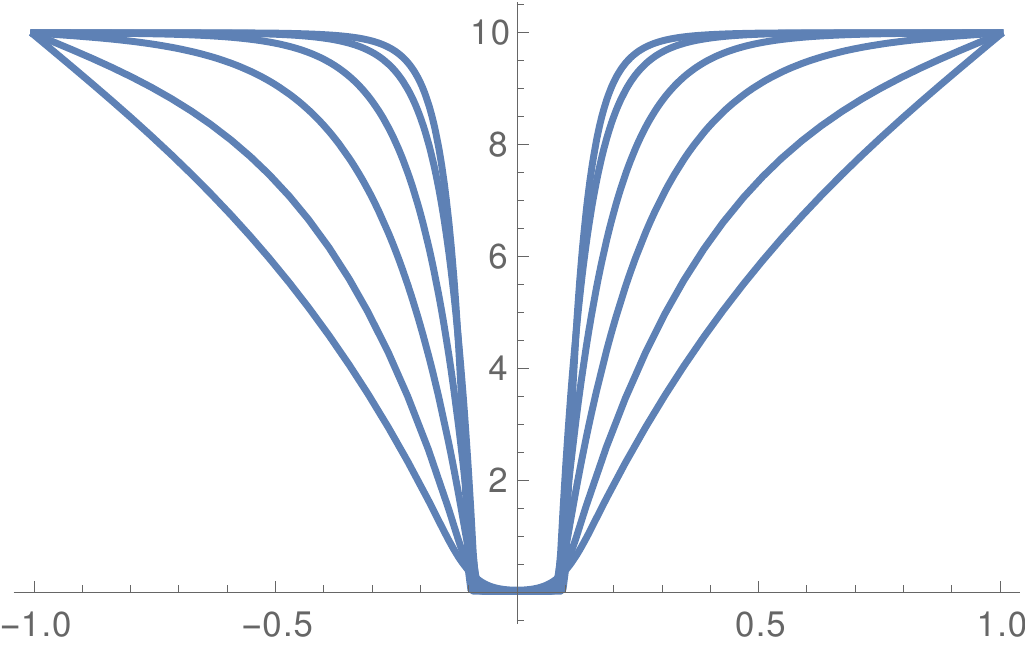}&
\includegraphics[width=0.31\textwidth]{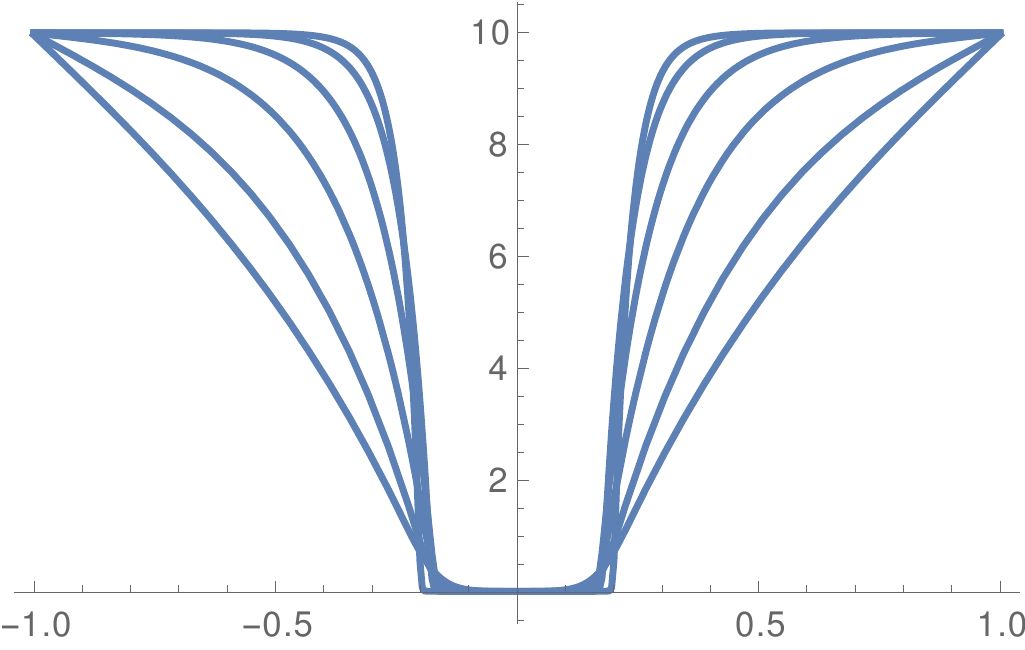}&
\includegraphics[width=0.31\textwidth]{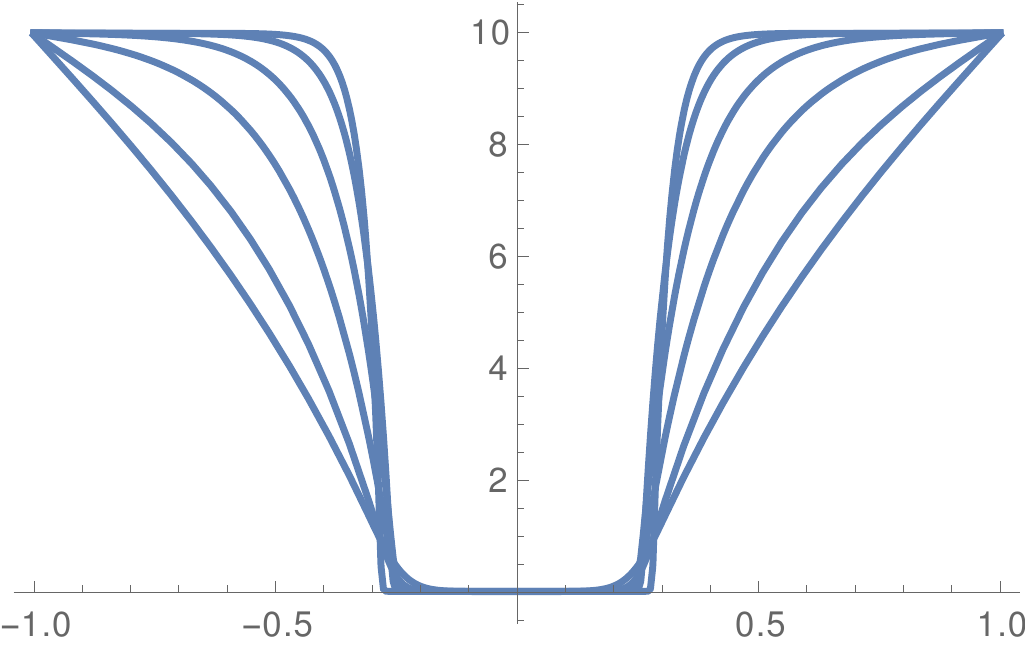}
\\
$v_\infty=0.4$ &  $v_\infty=0.8$ &  $v_\infty= 1.2$ 
\\
\includegraphics[width=0.25\textwidth]{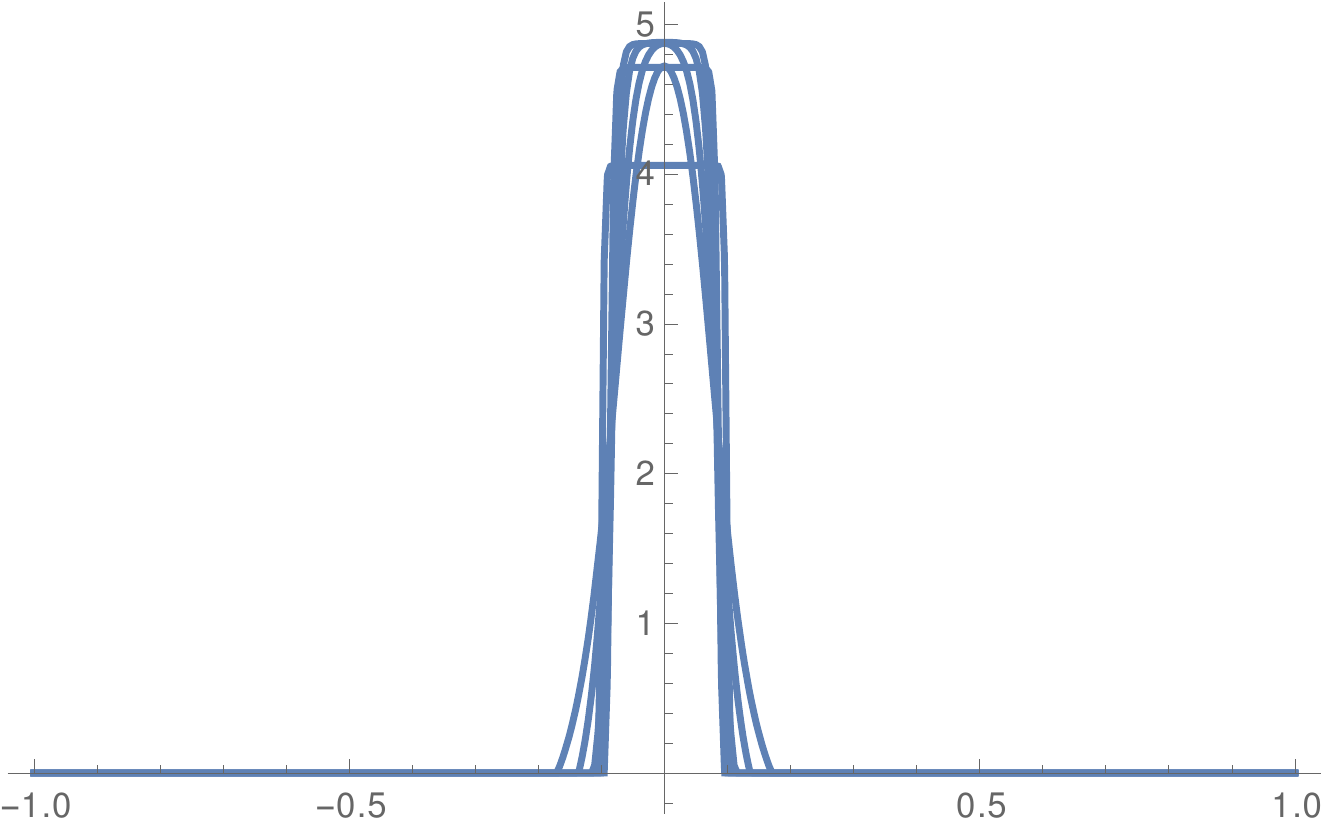}&
\includegraphics[width=0.25\textwidth]{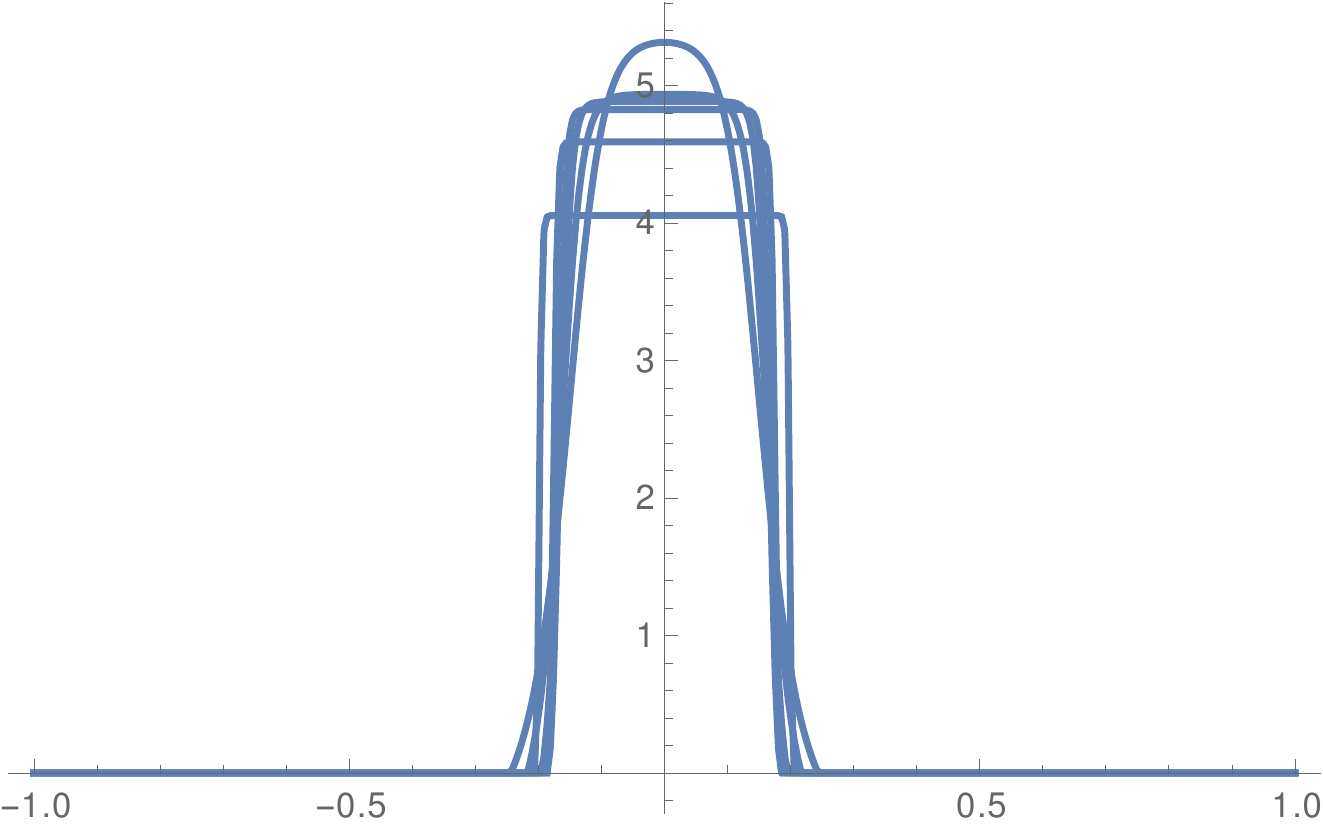}&
\includegraphics[width=0.25\textwidth]{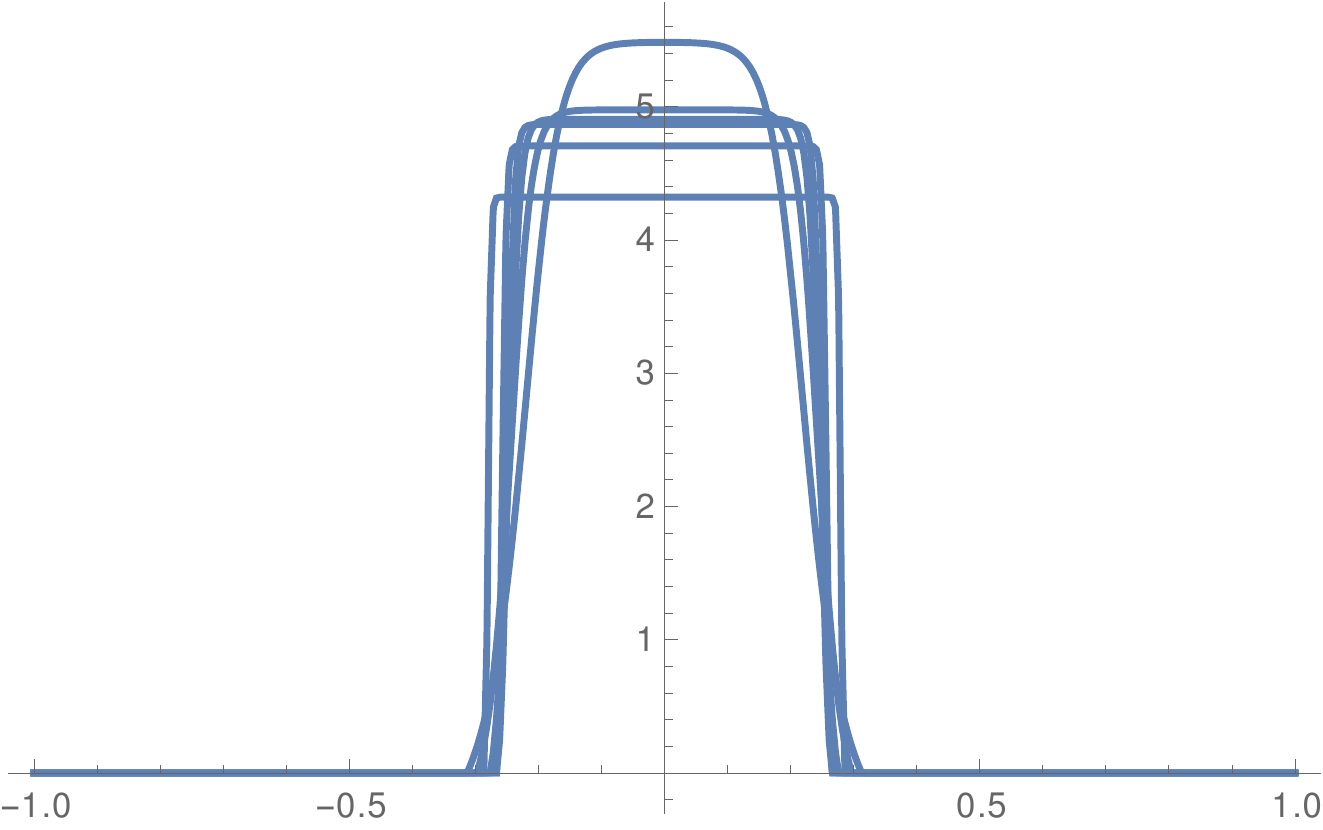}
\end{tabular}
\caption{\small\sl A study for the limit $\kappa\to 0^+$ of the steady state solutions
  $(\theta_\stst,\pi_\stst)$. For $v_\infty  \in \{0.4,0.8,1.2\}$ the profiles
  are plotted for $\kappa \in \{0.03, 0.01, 0.003, 0.001, 0.0003, 0.0001
  \}$. Convergence to rectangular profiles is observed.}
\label{fig:Lim.kappa0}
\end{figure}

Finally, we want to study the case corresponding to Proposition
\ref{prop3}, where $v_\infty$ is kept fixed and the limit $\kappa \to 0$ is
performed. In Figure \ref{fig:Lim.kappa0} we show plots of the steady states
$(\theta_\stst^\kappa,\pi_\stst^\kappa)$ for three different values of
$v_\infty$ for a sequence of decreasing $\kappa$. We clearly see the predicted
development of convergence against towards the limit
$(\theta_\stst^0,\pi_\stst^0)$ taking only two different values. Moreover, the
values are roughly independent of $v_\infty$, where the active plastic zone
$(-h,h)$ behaves like $h=v_\infty/\pi_*$, as proved in Proposition \ref{prop3}.

\subsection{An ODE model showing oscillations in time }
\label{su:NumSimODE}

Oscillatory behavior is most easily seen in a simple finite dimensional model,
consisting only of $\sigma(t)$ and $\ol\theta(t)$, where we may consider
$\ol\theta(t)$ as the average of $\theta(t,x)$ 
over the critical plasticity region where $\pi(t) = P(\sigma(t),\theta(t))$ is
positive. We also refer to the analysis of a spring-slider model in
\cite{Miel18TECI} as well as the geophysical paper \cite{AbeKat12CECS}.

Thus, our simplified model \eqref{eq:SM} is even more simplified to the ODE
system
\begin{equation}
  \label{eq:lumped}
  \frac{2H}{\bbC} \DT \sigma = 2v_\infty - 2h \,\varPi(\sigma,\ol\theta) \quad
\text{and} \quad 
\DT{\ol\theta} = 1 - \frac{\ol\theta}{\theta_\infty}
   - 10 \varPi(\sigma,\ol\theta)\,\ol\theta.  
\end{equation}
Here $h \in {]0,H[}$ represents the width of the plastic zone, which has to be
adapted accordingly.  We may consider \eqref{eq:lumped} as an evolutionary 
lumped-parameter system, which  in geophysical literature is often referred to
as a {\it 1-degree-of-freedom slider} and is considered as a basic test of
every new friction model.

\begin{figure}[h]
\begin{tabular}{ccc}
\includegraphics[width=0.37\textwidth]{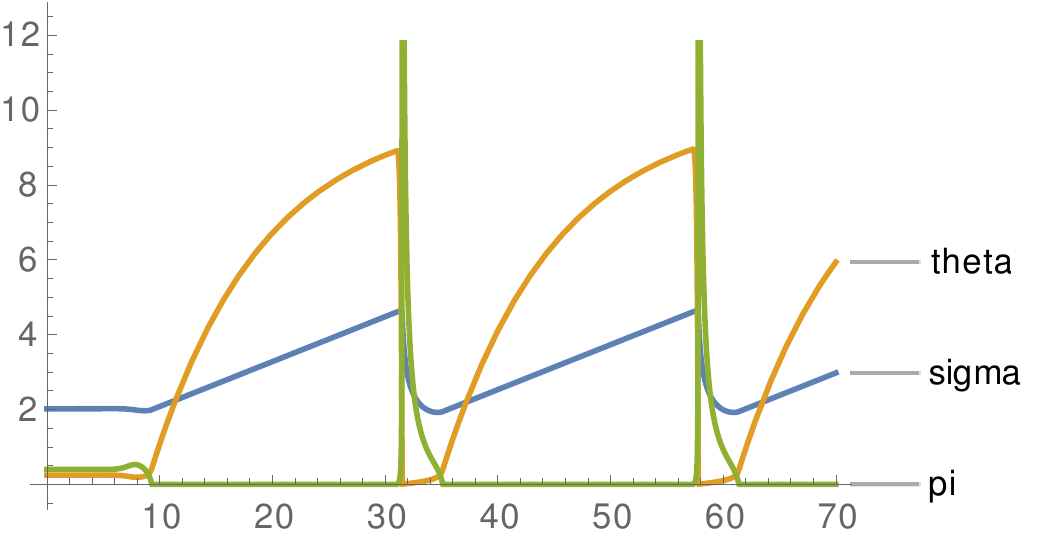}&
\includegraphics[width=0.27\textwidth,trim=0 0 66 0,clip=true]{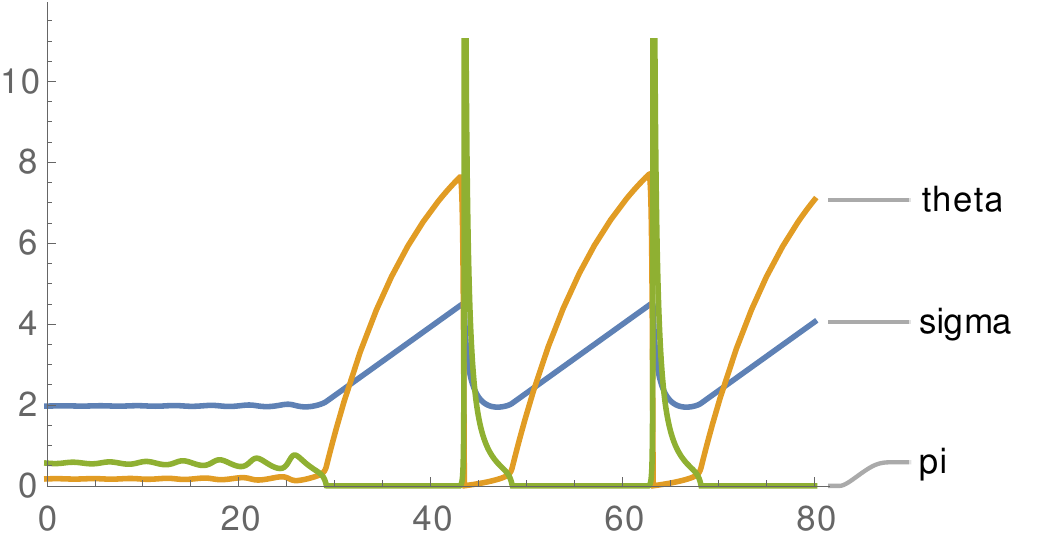}&
\includegraphics[width=0.27\textwidth,trim=0 0 90 0,clip=true]{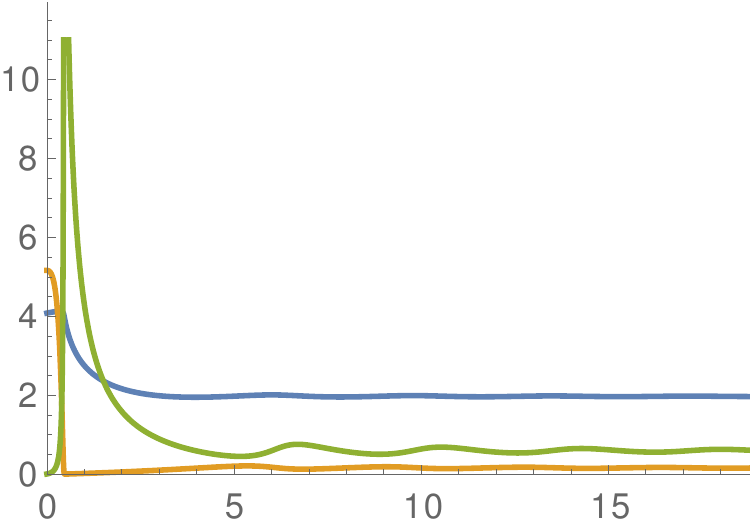}
\\
{\small$v_\infty=0.12$}&{\small$v_\infty=0.17$}&{\small$v_\infty=0.18$}
\end{tabular}
\caption{\small\sl Solutions $(\ol\theta(t),\sigma(t))$ together with
  $\pi(t)=P(\sigma(t),\ol\theta(t))$ for $h=0.3$ and three different values of $v_\infty$. In the first two cases the solutions
  start very close to the unstable steady state. In the third case the
  solution starts far away but soon returns to the stable fixed point.}
\label{fig:ODE}
\end{figure}

The nice feature of this ODE model is that the steady states
can be calculated explicitly, and even a stability analysis can be
performed. Indeed there is exactly one steady state, namely
\[
\ol\theta_\stst = \frac{\theta_\infty }{1{+}10 
  (v_\infty/h) \theta_\infty} \quad \text{ and } \quad 
\sigma_\stst =\mu_0+A\Big(\frac{v_\infty}h\Big)+B(\ol\theta_\stst ). 
\]
Instead of performing a rigorous analysis, we simply display the solution
behavior of this ODE by a few numerical results. We find that for small
positive $v_\infty$ we obtain oscillatory behavior, while for larger $v_\infty$
the solutions converge to the steady state, see Figure \ref{fig:ODE}. Indeed,
the oscillations can be interpreted physically in terms of geophysical
processes as {\it seismic cycles}. 

During the oscillatory behavior there is a large part of the interval where
there is no plastic slip (i.e.\ $\pi(t)=0$).  In these intervals the stress is
growing linearly with a slope that is proportional to $v_\infty$, and the aging
variable $\ol\theta$ is relaxing exponentially back to its equilibrium value
$\theta_\infty$. However, if the stress reaches a critical value, then the
plastic strain rate is triggered, which leads to reduction of the aging
variable. This leads to a simultaneous weakening of the plastic yields stress
$\mu(\pi,\ol\theta)$ such that $\pi$ can grow even more. As a result the stress
is drastically reduced in a rather short time interval, and $\ol\theta$ is
reduced almost down to $0$ (refreshing). If the inertial term would be
included, then this fast rupture-like processes could emit elastic waves, i.e.\
{\it earthquakes}. Because of the stress release
the plastic strain rate reduces to $0$, and the process starts
again by a slow aging and building up the stress.

\begin{figure}[h]
\centering \begin{tikzpicture}
\node[above] at (0,0){\includegraphics[width=0.3\textwidth]{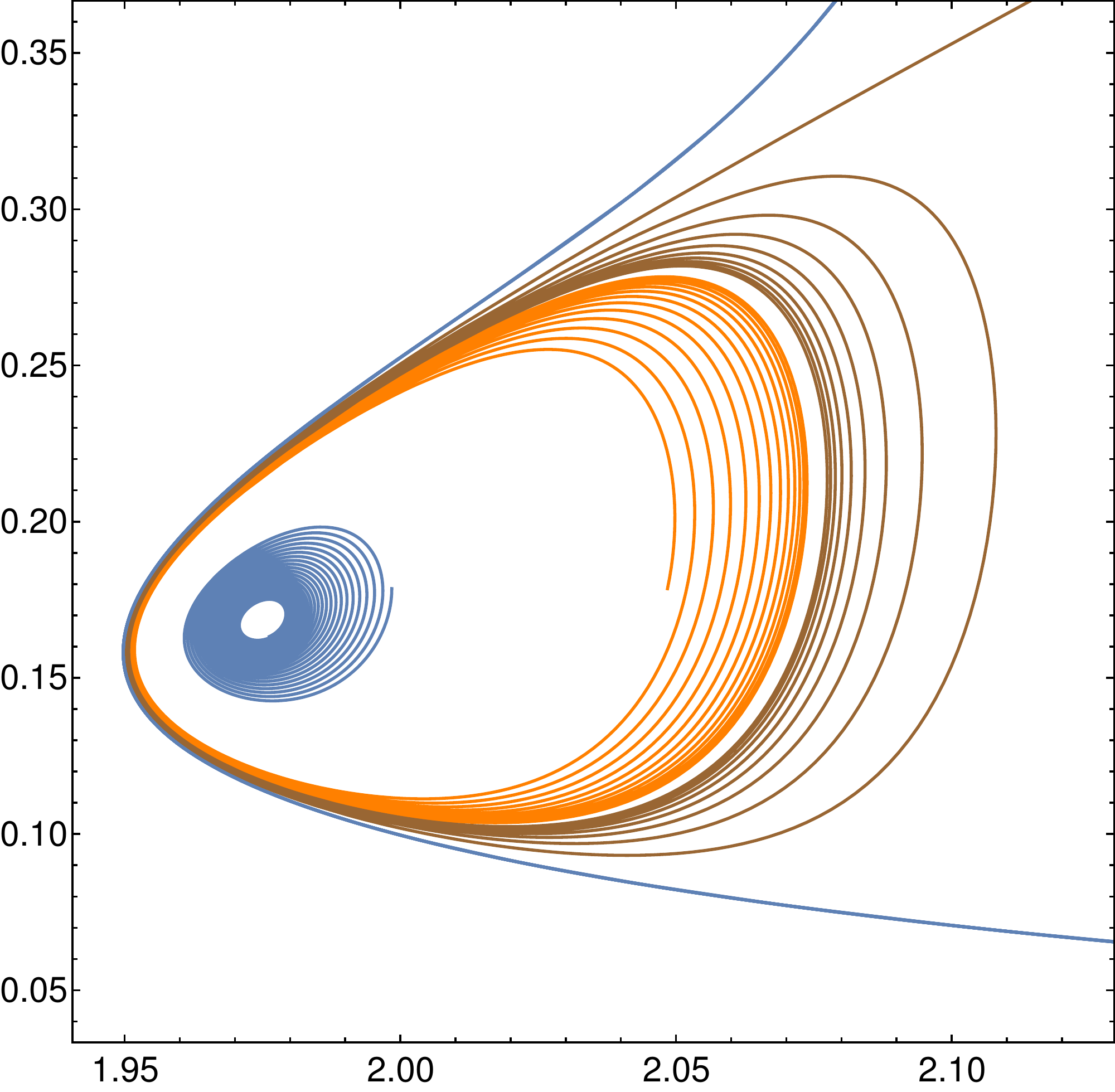}};
\node[above] at (7,0){\includegraphics[width=0.37\textwidth]{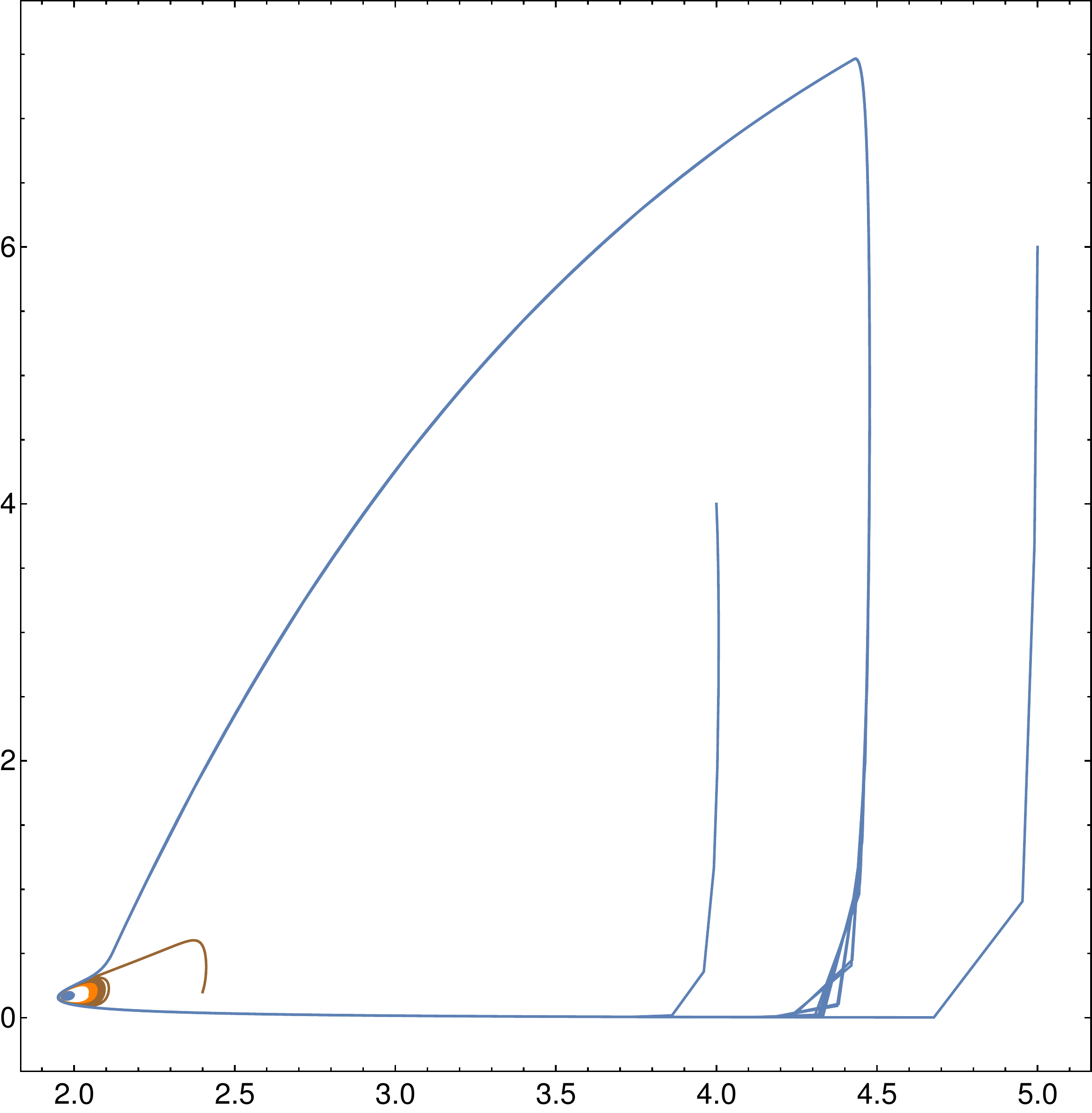}}; 
\draw (4.3,0.6) rectangle (4.7,1);
\draw (4.3,0.6) --(2.45,0.35);
\draw (4.3,1.0) --(2.45,4.8);
\node[right] at (10.2,0.35){$\sigma$};
\node[left] at (4,6){$\ol\theta$};
\end{tikzpicture}\par

\caption{\small\sl The $(\sigma,\ol\theta)$ phase plane for $h=0.3$ and $v_\infty=0.175$,
  where all trajectories rotate clockwise around the fixed point
  $(\sigma_\stst , \ol\theta_\stst )\approx (1.973,0.168)$. 
There are two periodic solutions. The outer one is stable and is approached by
the blue trajectories from inside and outside. The unstable periodic orbit lies
between the orange and the brown trajectory.}
\label{fig:TwoPeriodicOrbits}
\end{figure}
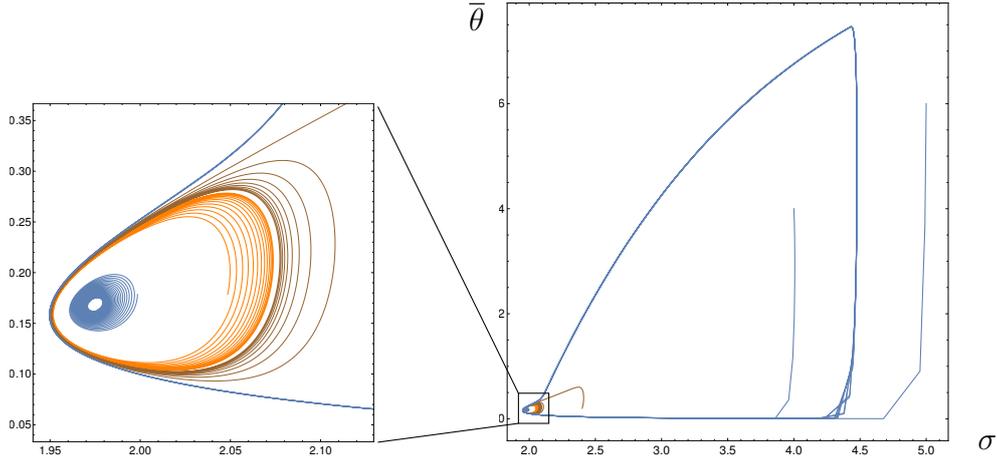

In fact, choosing $h=0.3$ a closer analysis of the system shows that the steady
states are stable if and only if $v>v_\infty^{(1)}\approx 0.17462$. However,
stable oscillations are already seen for $ v< v_\infty^{(2)} \approx
0.175452$. A careful analysis of the trajectories in the phase plane for 
$(\ol\theta,\sigma)$ reveals that for $v_\infty \in (v_\infty^{(1)},
v_\infty^{(2)})$ there are two periodic solutions, as smaller unstable one that
encircles the stable fixed point and a larger stable one that encircles the
unstable one, see Figure \ref{fig:TwoPeriodicOrbits}. Thus, in the small
parameter interval $(v_\infty^{(1)}, v_\infty^{(2)})$ we have coexistence of a
stable fixed point and a stable periodic orbit. 

\subsection{Convergence to steady states versus oscillations for \eqref{eq:SM}}
\label{su:NumSimPDE} 

The behavior of the evolutionary coupled system \eqref{eq:SM} coupling the
parabolic PDE for the aging variable $\theta(t,x)$ to the ODE for the stress
$\sigma(t)$ displays roughly a similar behavior as the lumped ODE system
\eqref{eq:lumped}. For large $|v_\infty|$ one observes convergence into the
steady states analyzed in Section \ref{se:AnaSteady} and displayed numerically
in Section \ref{su:NumSimSteady}. For small nontrivial values of $v_\infty$ one
observes oscillatory behavior. Of course, the new feature is the spatial
distribution of the plastic rate $\pi(t,x)=\varPi(\sigma(t),\theta(t,x))$ and
the aging variable $\theta(t,x)$. In most cases one observes that $\pi(t,x)$
has a nontrivial support in the sense that the support of $\pi(t,\cdot)$ is
compactly contained in $({-}H,H)$. Moreover, in the oscillatory case, we also
observe that there are large parts of the periodicity interval, in which there
is no plastic flow at all (i.e.\ $\pi=\DT p=0$), but there is aging and
slow building up of stress. Then, in sudden plastic bursts there is a
strong plastic flow that leads to stress release and refreshing, i.e.\
reduction of $\theta$ almost down to $0$ inside the cataclastic zone.

\begin{figure}[h]
\centering 
\begin{tabular}{cc@{\qquad}cc}
\includegraphics[width=0.27\textwidth]{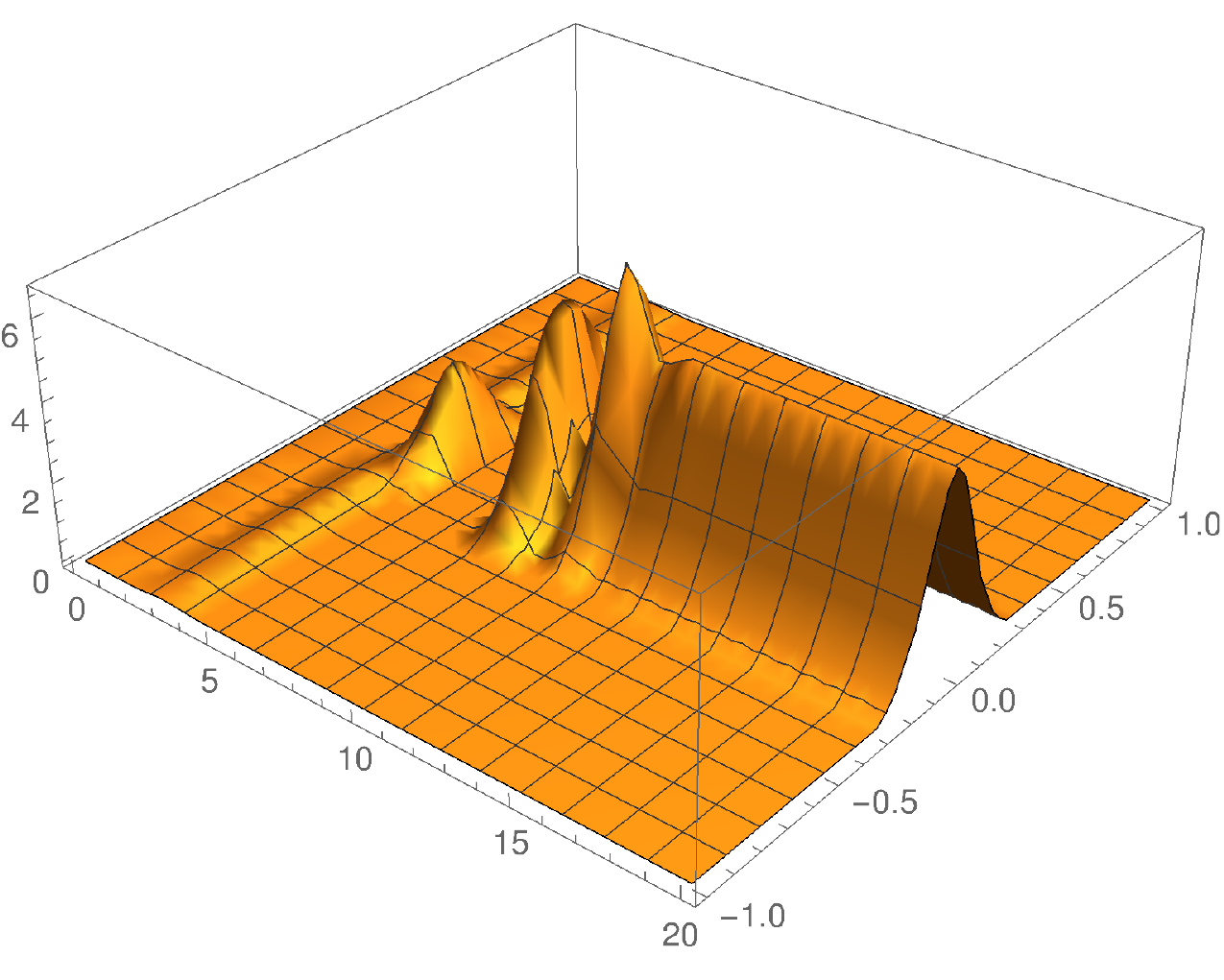}\hspace*{-1.3em}&
\includegraphics[width=0.22\textwidth]{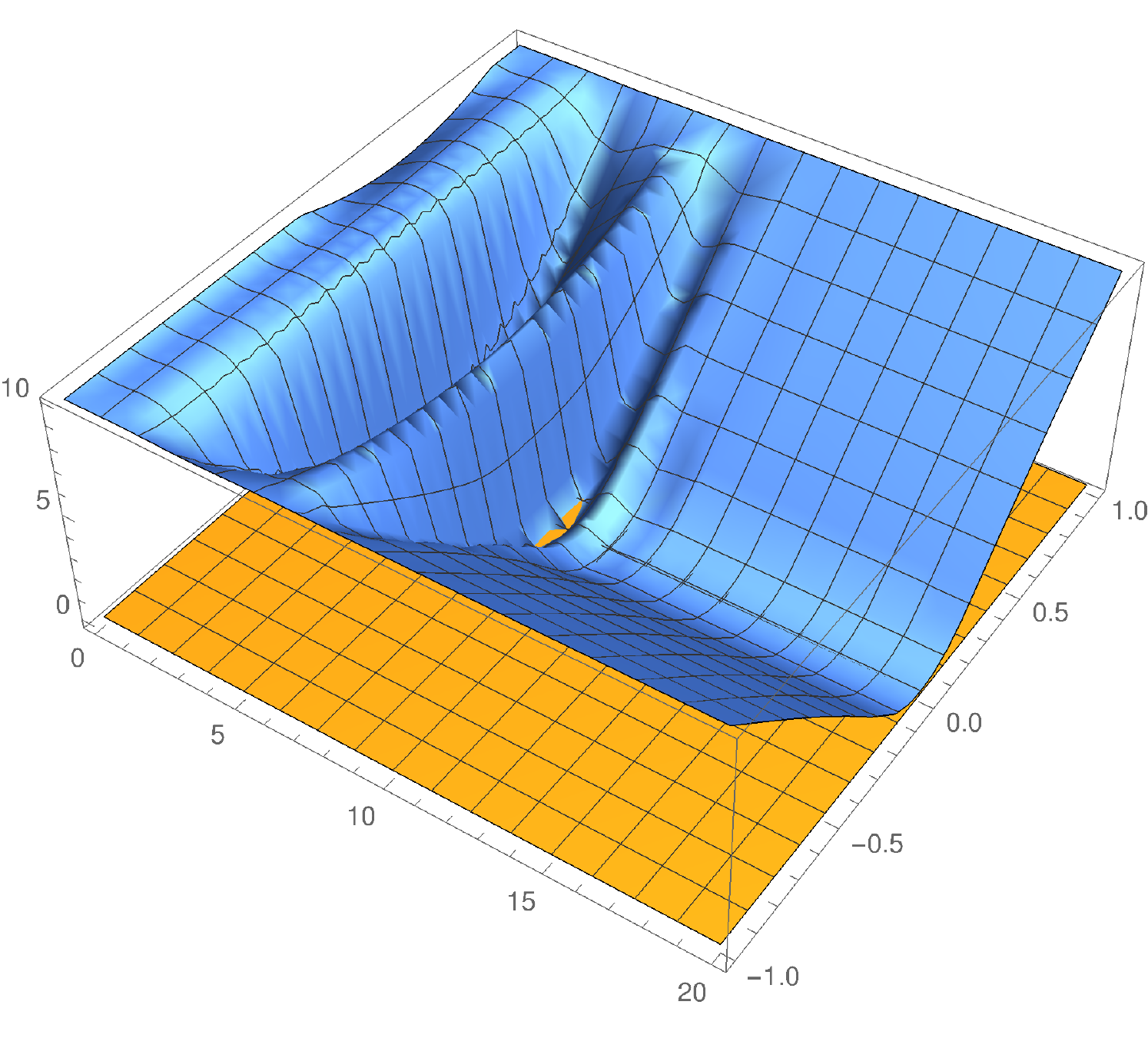}\hspace*{-1.3em}&
\includegraphics[width=0.27\textwidth]{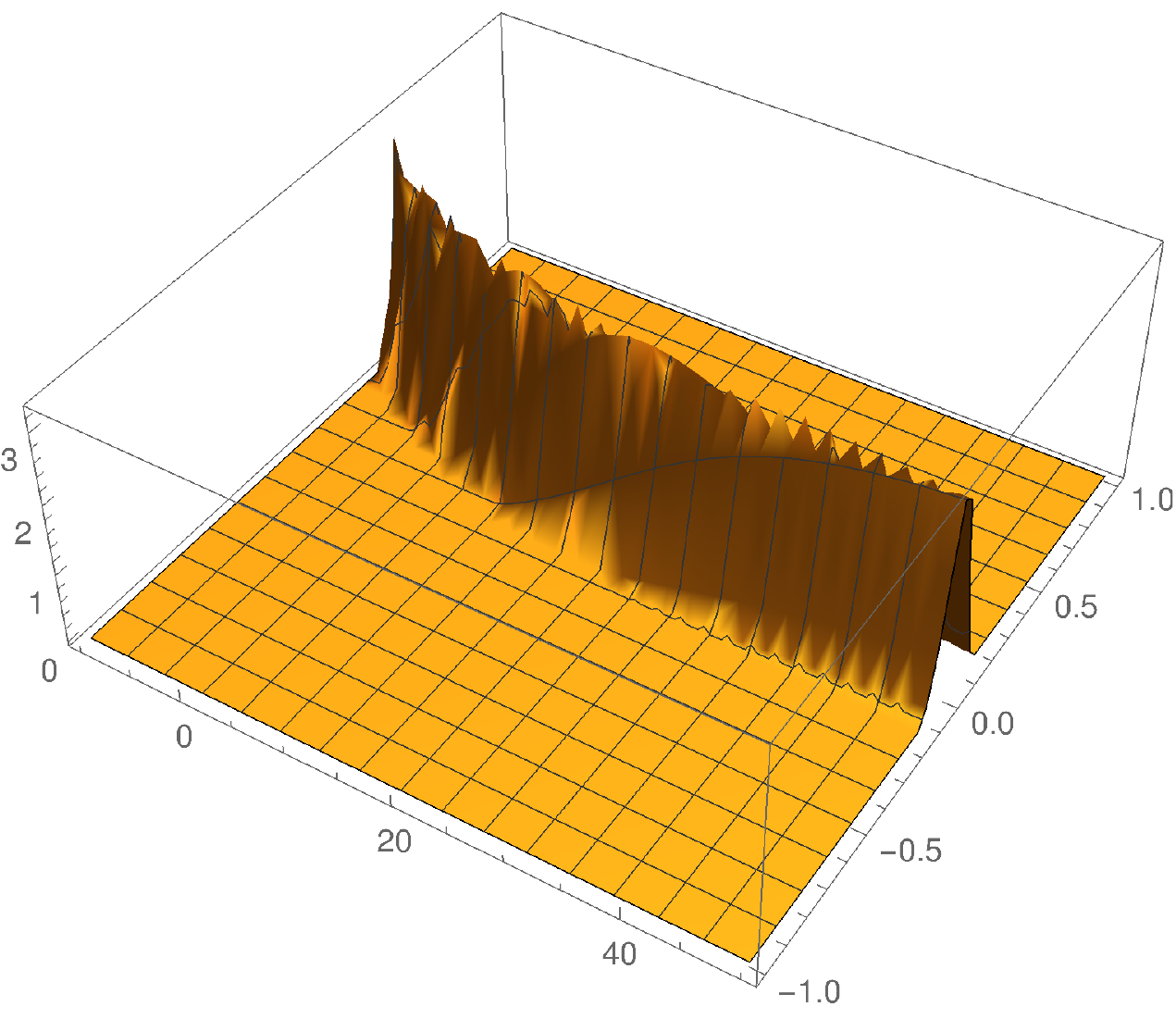}\hspace*{-1.3em}& 
\includegraphics[width=0.22\textwidth]{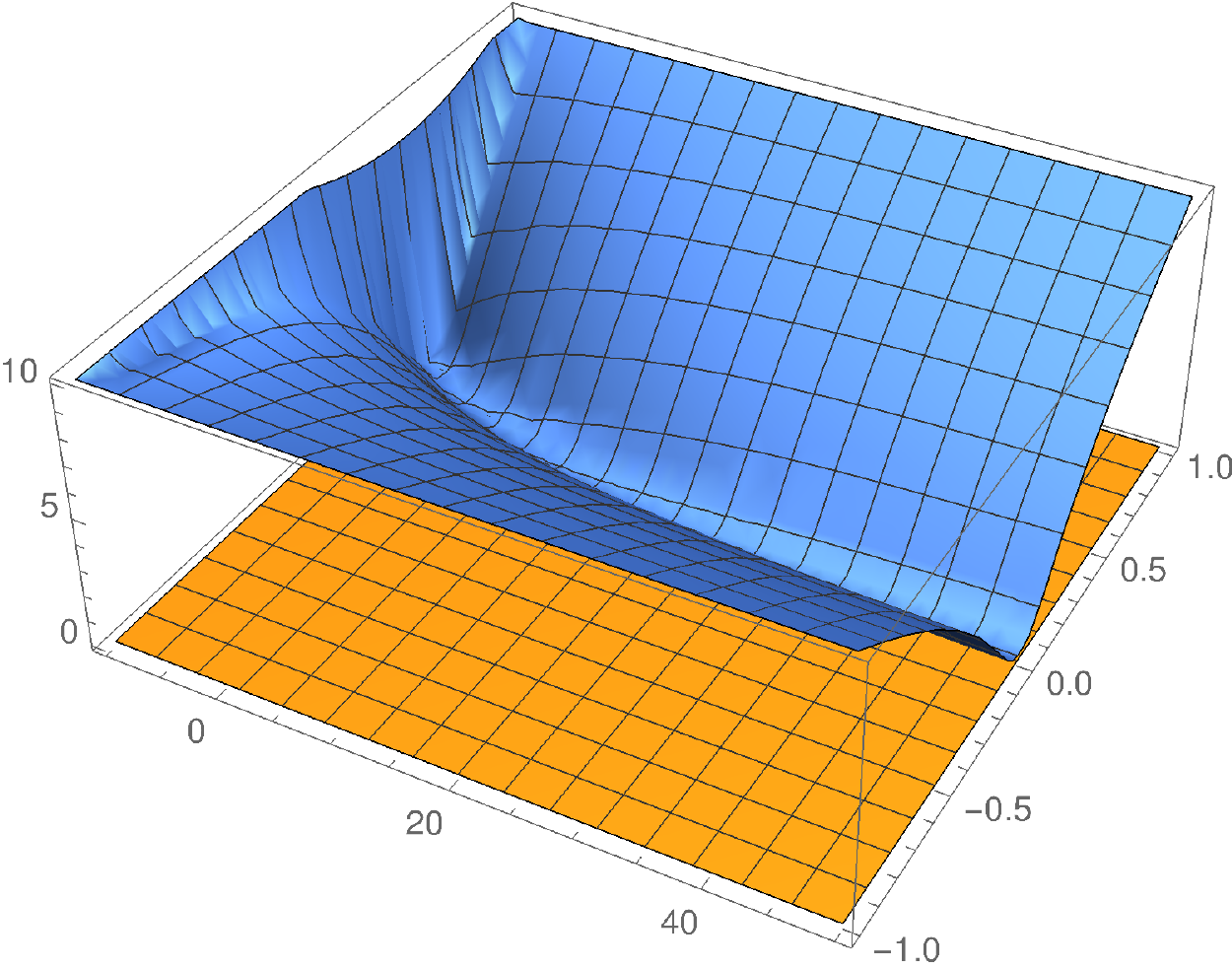}
\\
\multicolumn{2}{c}{\small$\kappa=0.16,\ v_\infty=0.6$}&
\multicolumn{2}{c}{\small$\kappa=0.004,\ v_\infty=0.2$}
\end{tabular}
\par
\caption{\small\sl Simulation of the solution $\theta$ (left) and
  $\pi=\varPi(\sigma,\theta)$ (right) for \eqref{eq:SM}.
  Convergence to a steady state can be observed in both cases.}
\label{fig:SM.converge}
\end{figure}

Figure \ref{fig:SM.converge} displays two simulation results featuring convergence
into steady state.     

\begin{figure}[h]
\centering 
\begin{tikzpicture}
\node at (0,3)
{\includegraphics[width=0.97\textwidth,trim=150 0 150 0,clip=true,angle=2]
  {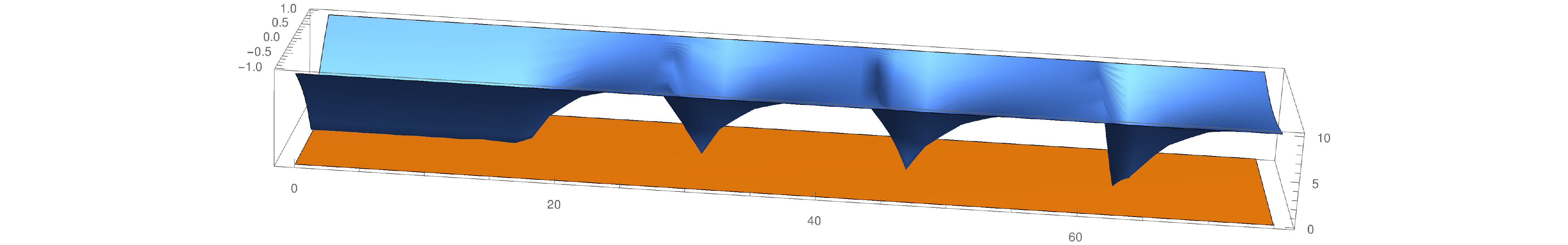}}; 
\node at (-1.3,1.8) { time $t$};
\node at (7.5,3.8) {$\theta(t,x)$};

\node at (0,0)
 {\includegraphics[width=0.85\textwidth,angle=4.3]{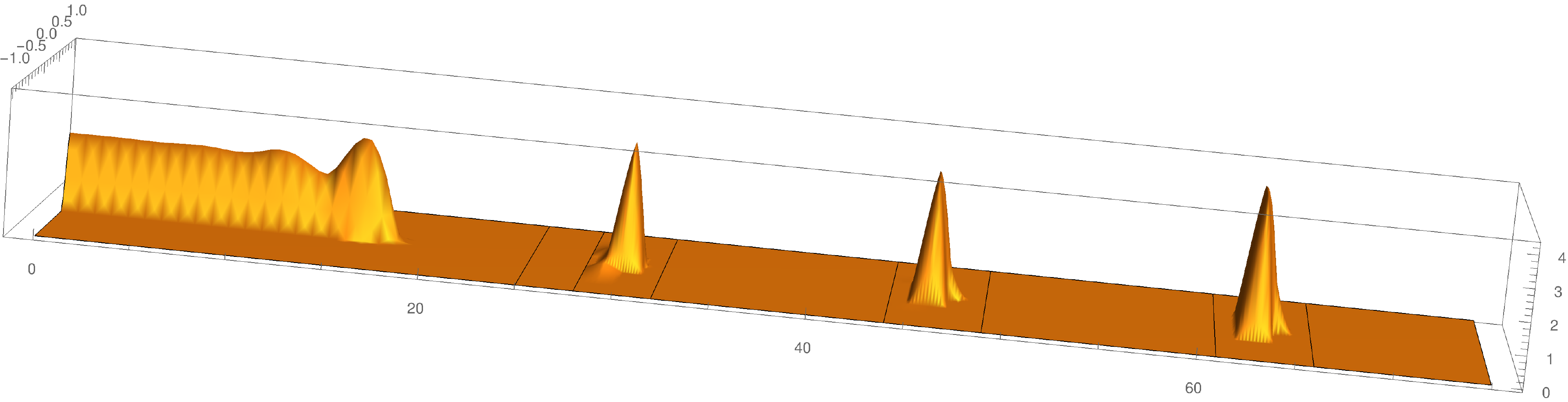}}; 
\node at (-1.3,-1.3) { time $t$};
\node at (7.5,0.4) {$\pi(t,x)$};
\end{tikzpicture}
\vspace{-2em}

\caption{\small\sl Simulation of the solution $\theta$ (top) and
  $\pi=\varPi(\sigma,\theta)$ (bottom) for \eqref{eq:SM} with $\kappa=0.04$ and
  $v_\infty=0.15$. Convergence to a periodic behavior where $\pi$ is
  localized in space and time can be observed.}
\label{fig:SM.osc}
\end{figure}
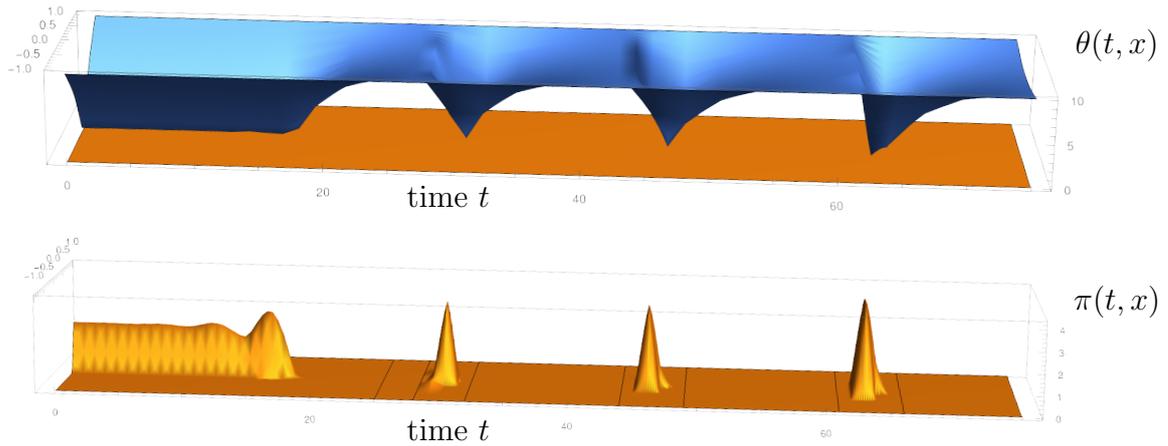

In the case $\kappa=0.04$ and the smaller shear rate $v_\infty=0.15$ one
observes oscillatory behavior. In fact, we start the solution very close to the
steady state and the solution needs some time to develop the instability but
then it switches quickly into a periodically looking regime,
see Figure \ref{fig:SM.osc}.

\bigskip\bigskip\bigskip

\paragraph*{Acknowledgments.}
A.M.\ was partially supported by DFG via the Priority Program 
SPP\,2256 
\emph{Variational Methods for Predicting Complex Phenomena in 
Engineering Structures and Materials} (project no.\,441470105, subproject Mi 459/9-1 \emph{Analysis for thermo-mechanical models with internal
  variables}). 
T.R.\ is thankful for the hospitality of the Weierstra\ss{}--Institut
Berlin and also acknowledges the support of the M\v SMT \v CR
(Ministry of Education of the Czech Republic) project
CZ.02.1.01/0.0/0.0/15-003/0000493, and the institutional support RVO:
61388998 (\v CR).  

\end{sloppypar}
\end{document}